\documentclass[reqno]{amsproc}
\usepackage{amssymb}
\usepackage{amsmath}
\usepackage{amsfonts}
\usepackage[a4paper,left=25mm,right=25mm,top=25mm,bottom=25mm,footskip=10mm]{geometry}

\makeatletter
\def\th@boldplain{  \th@plain
  \thm@headfont{\normalfont\bfseries}}
\makeatother
\theoremstyle{plain}
\theoremstyle{boldplain}
\newtheorem{theorem}{Theorem}[section]
\theoremstyle{plain}

\theoremstyle{boldplain}
\newtheorem{corollary}[theorem]{Corollary}
\theoremstyle{plain}

\theoremstyle{boldplain}
\newtheorem{definition}[theorem]{Definition}
\theoremstyle{plain}
\theoremstyle{boldplain}
\newtheorem{example}[theorem]{Example}
\theoremstyle{plain}

\theoremstyle{boldplain}
\newtheorem{lemma}[theorem]{Lemma}
\theoremstyle{plain}

\theoremstyle{boldplain}
\newtheorem{proposition}[theorem]{Proposition}
\theoremstyle{plain}
\theoremstyle{boldplain}
\newtheorem{remark}[theorem]{Remark}
\theoremstyle{plain}

\numberwithin{equation}{section}
\IfFileExists{tcilatex.tex}{\input{tcilatex}}{  \providecommand{\dsum}{\displaystyle\sum}  \providecommand{\dprod}{\displaystyle\prod}  \providecommand{\dint}{\displaystyle\int}}
\providecommand{\func}[1]{\operatorname{#1}}
\begin{document}
\title[Stationary scattering theory]{Stationary scattering theory for
commuting $d$-tuples of unitaries in von Neumann algebras}
\author{Qihui Li}
\address{School of Science, East China University of Science and Technology,
Shanghai, 200237, China}
\email{qihui\_li@126.com}
\thanks{The author was partly supported by the NSFC (Grant No. 12671153) and
the Natural Science Foundation General Project of Shanghai (24ZR1415600).}
\subjclass[2010]{Primary 47A40; Secondary 47A55,47C15}
\keywords{Mathematical scattering theory, wave operators, Kato--Rosenblum
theorem, von Neumann algebras. }

\begin{abstract}
In this paper, we develop scattering theory for commuting $d$-tuples of
unitaries in a von Neumann algebra $\mathcal{M}$. Our approach is based on
the stationary method. We introduce jointly smooth operators and define the
joint norm absolutely continuous support $P_{ac}^{\infty }\left( \boldsymbol{%
U},\mathcal{M}\right) $ for a commuting $d$-tuple $\boldsymbol{U=}\left(
U_{1},\cdots ,U_{d}\right) $ in $\mathcal{M},$ which extends the notion of
the absolutely continuous projection $P_{ac}\left( \boldsymbol{U}\right) $
to the setting of von Neumann algebras. For $d\geq 2,$ we consider two $d$%
-tuples of commuting unitaries $\boldsymbol{U}$ and $\boldsymbol{U}_{0}$ and
an identification operator $J\in \mathcal{M}.$ Under the factorizations $%
JU_{0,k}-U_{k}J=G_{k}^{\ast }G_{0,k}$ with $G_{k}$ jointly $\boldsymbol{U}$%
-smooth and $G_{0,k}$ jointly $\boldsymbol{U}_{0}$-smooth, we use products
of resolvents to construct sector stationary wave operators $w_{\boldsymbol{%
\sigma }}=w_{\boldsymbol{\sigma }}\left( \boldsymbol{U},\boldsymbol{U}_{0},J,%
\mathcal{M}\right) $ and relate them to the strong wave operator $W=W\left( 
\boldsymbol{U},\boldsymbol{U}_{0},J,\mathcal{M}\right) $. We first prove
that $W$ exists, coincides with $w_{\boldsymbol{\sigma }}$ for every sector $%
\boldsymbol{\sigma }\in \{+1,-1\}^{d}$ and satisfies $WU_{0,k}=U_{k}W$ for $%
k=1,\cdots ,d.$ Second, when the identification operator $J$ is unitary, we
show that $W^{\ast }W=P_{ac}^{\infty }\left( \boldsymbol{U}_{0},\mathcal{M}%
\right) $ and $WW^{\ast }=P_{ac}^{\infty }\left( \boldsymbol{U},\mathcal{M}%
\right) .$ For $d=1,$ let $\mathcal{M}$ be countably decomposable and
properly infinite, with a faithful normal semifinite tracial weight $\tau $
and let $J$ be a unitary in $\mathcal{M}.$ We prove that $W_{\pm }=W_{\pm
}\left( U,U_{0},J,\mathcal{M}\right) $ exists and $W_{\pm }=w_{\pm }\left(
U,U_{0},J,\mathcal{M}\right) $ provided $JU_{0}-UJ\in \mathcal{L}^{1}\left( 
\mathcal{M},\tau \right) \mathcal{\cap M}$ where $w_{\pm }\left( U,U_{0},J,%
\mathcal{M}\right) $ is the corresponding stationary wave operator. We also
establish the identities $W_{\pm }^{\ast }W_{\pm }=P_{ac}^{\infty }\left(
U_{0},\mathcal{M}\right) $ and $W_{\pm }W_{\pm }^{\ast }=P_{ac}^{\infty
}\left( U,\mathcal{M}\right) $. Consequently, no unitary $U$ in $\mathcal{M}$
with $P_{ac}^{\infty }\left( U,\mathcal{M}\right) \neq 0$ can be written as $%
U=D+K,$ where $D\in \mathcal{M}$ is a diagonal unitary and $K\in \mathcal{L}%
^{1}\left( \mathcal{M},\tau \right) \mathcal{\cap M}.$
\end{abstract}

\maketitle

\section{Introduction}

\bigskip\ Mathematical scattering theory provides a rigorous framework for
describing particle collisions in quantum mechanics and the interaction of
waves with obstacles or inhomogeneous media. It studies the long-time
behavior of dynamical systems by comparing a perturbed evolution with a
simpler reference evolution. Its central objects are wave operators, which
relate the asymptotic states of the two systems, and the scattering
operator, which connects incoming and outgoing states. Fundamental questions
concern the existence and completeness of wave operators, as well as their
relation to the spectral properties of the underlying operators. There are
two principal approaches to studying the existence of wave operators: the
time-dependent method, based on limits of the evolutions, and the stationary
method, based on resolvent boundary values and spectral representations.

Let $\mathcal{H}$ be a complex Hilbert space and let $\mathcal{B(H)}$ be the
algebra of all bounded linear operators on $\mathcal{H}.$ Mathematical
scattering theory is an area of perturbation theory (\cite{Kato4}). In \cite%
{Y}, the stationary scattering theory of self-adjoint operators is presented
in detail. Yafaev also emphasizes that the methods used for the scattering
theory of self-adjoint operators are subdivided into two classes, i.e., the
smooth method and the trace class method. The smooth and trace class
approaches do not overlap one another, and both have useful applications.
For two densely defined self-adjoint operators $H_{0}$ and $H$, if $H_{0}-H$
is in the trace class, then the Kato-Rosenblum theorem (\cite{Kato},\cite{Ro}%
) asserts that the wave operator $W_{\pm }\left( H,H_{0}\right) $ of $H$ and 
$H_{0}$ in $\mathcal{B(H)}$ exists and consequently the absolutely
continuous parts of $H$ and $H_{0}$ are unitarily equivalent. Voiculescu
developed diagonalization criteria for commuting self-adjoint tuples modulo
normed ideals using quasicentral approximate units. He also applied
quasitriangularity to trace estimates for self-commutators and related
perturbation invariants to singular integrals, hyperbolicity, and dynamical
entropy (\cite{Voi2}-\cite{Voi9}). In \cite{B}, Bhutani and Voiculescu
discuss a sharp Kato-Rosenblum type theorem for unbounded $d$-tuples of
commuting self-adjoint operators. The method used in \cite{B} can be
regarded as a time-dependent approach.

For a pair of unitary operators $U_{0}$ in $\mathcal{B(H}_{0}\mathcal{)}$
and $U$ in $\mathcal{B(H)}$ where $\mathcal{H}_{0}$ is a complex Hilbert
space, stationary scattering theory constructs wave operators of $U_{0}$ and 
$U$ through boundary values of resolvents on the unit circle. A systematic
abstract treatment of them was given by R. Tiedra de Aldecoa \cite{TR} in
2020. For the factorization $JU_{0}-UJ=G^{\ast }G_{0}$ where $J\in \mathcal{%
B(H}_{0},\mathcal{\mathcal{H})}$ is an identification operator, that work
provides sufficient smoothness and resolvent boundary-value conditions for
the existence of stationary wave operators $w_{\pm }\left( U,U_{0},J\right) $%
, establishes their coincidence with the strong wave operators $W_{\pm
}\left( U,U_{0},J\right) $, and derives representation formulas for the wave
operators and the scattering matrix. When $J$ is a unitary operator from $%
\mathcal{H}$ to $\mathcal{H}_{1},$ the wave operators implement unitary
equivalences between the absolutely continuous parts of $U_{0}$ and $U$.

Let $\mathcal{M\subseteq B(H)}$ be a von Neumann algebra with unit $1$. A
von Neumann algebra is called \textquotedblleft
semifinite\textquotedblright\ if it has a faithful normal semifinite tracial
weight. (For general knowledge about von Neumann algebras, the reader is
referred to \cite{MurrayVonNeuman1}-\cite{MurrayVonNeuman5} or \cite{Kadison}%
). In \cite{Li}, we introduce the concept of a generalized wave operator $%
W_{\pm }$ based on the notion of norm absolutely continuous projections. An
analogue of the Kato-Rosenblum theorem in a semifinite von Neumann algebra $%
\mathcal{M}$ with a faithful normal semifinite tracial weight $\tau $ is
obtained by showing the existence of the generalized wave operator $W_{\pm
}. $ To be more precise, we proved that for self-adjoint $H$ and $H_{0}$
affiliated with $\mathcal{M}$ satisfying $H-H_{0}\in \mathcal{L}^{1}\left( 
\mathcal{M},\tau \right) \mathcal{\cap M},$ the generalized wave operator $%
W_{\pm }\left( H,H_{0}\right) $ exists and then the norm absolutely
continuous parts of $H$ and $H_{0}$ are unitarily equivalent.

To the best of our knowledge, the existence of wave operators for $d$-tuples
of commuting unitaries in von Neumann algebras has not yet been
investigated. In this paper, our purpose is to fill in this gap by
developing scattering theory for commuting $d$-tuples of unitaries in a von
Neumann algebra $\mathcal{M}$ by applying stationary approaches.

One important contribution of this paper is the introduction of the joint
norm absolutely continuous support $P_{ac}^{\infty }\left( \boldsymbol{U},%
\mathcal{M}\right) $ for a commuting $d$-tuple $\boldsymbol{U}=\left(
U_{1},\cdots ,U_{d}\right) $ of unitaries in a von Neumann algebra $\mathcal{%
M}$. Motivated by the norm absolutely continuous support for a single
self-adjoint operator introduced by Li, Shen, Shi, and Wang in \cite{Li}, we
define this joint support as the supremum of all jointly smooth projections
in $\mathcal{M}$. The definition depends only on the tuple and the ambient
algebra, independently of the identification operator or the perturbation
under consideration. In particular, it belongs to the center of the relative
commutant in $\mathcal{M}$ of the algebra generated by the tuple, and its
range therefore reduces every component. When $\mathcal{M}$ is the full
algebra $\mathcal{B(H)}$, $P_{ac}^{\infty }\left( \boldsymbol{U},\mathcal{%
B(H)}\right) $ coincides with the usual joint absolutely continuous
projection. For a general von Neumann algebra, however, it can be strictly
smaller, as illustrated by an example. The scattering significance of this
support is established by our completeness theorem. Under the joint
smoothness and factorization assumptions of the theorem, and with a unitary
identification operator, we prove that the initial projection of the strong
wave operator is exactly the joint norm absolutely continuous support of the
reference tuple, while its final projection is exactly that of the perturbed
tuple. The wave operator consequently restricts to a unitary map between
these two subspaces and simultaneously intertwines all components of the
tuples. This establishes scattering completeness relative to the joint norm
absolutely continuous supports and gives these independently defined
projections a direct interpretation as the initial and final projections of
the wave operator.

Let 
\begin{equation*}
\boldsymbol{U}=\left( U_{1},\cdots ,U_{d}\right) \text{ and }\boldsymbol{U}%
_{0}=\left( U_{0,1},\cdots ,U_{0,d}\right)
\end{equation*}%
be two $d$-tuples of commuting unitaries in a von Neumann algebra $\mathcal{M%
}$ and let $J\in \mathcal{M}$. For $d\geq 2$, a central difficulty is to
control the joint $\mathbb{Z}^{d}$-action uniformly as $\left\Vert \mathbf{m}%
\right\Vert _{\infty }\rightarrow \infty $. This limit allows some
coordinates to remain bounded and therefore requires estimates that are
uniform in the remaining coordinates. Separate applications of one-operator
scattering results do not, by themselves, provide this control or identify
the limits associated with different sign sectors. We address these issues
through joint smoothness and products of resolvents. A telescoping expansion
of the product-resolvent difference incorporates the coordinate
factorizations $JU_{0,k}-U_{k}J=G_{k}^{\ast }G_{0,k}$, together with joint
smoothness assumptions on the factors. Using products of resolvents and
their boundary values under a common radial regularization, we construct
sector stationary wave operators $w_{\boldsymbol{\sigma }}:=w_{\boldsymbol{%
\sigma }}\left( \boldsymbol{U},\boldsymbol{U}_{0},J,\mathcal{M}\right) $ and
relate them to the strong wave operator $W:=W\left( \boldsymbol{U},%
\boldsymbol{U}_{0},J,\mathcal{M}\right) $. We establish two main results.
First, under a certain condition, we show that the strong wave operator $W$
and sector stationary wave operators $w_{\boldsymbol{\sigma }}$ exist, and 
\begin{equation*}
W=w_{\boldsymbol{\sigma }}\text{ for every sector }\boldsymbol{\sigma }\in
\left\{ +1,-1\right\} ^{d}.
\end{equation*}%
Moreover, 
\begin{equation*}
WU_{0,k}=U_{k}W.
\end{equation*}%
Second, when the identification operator $J$ is unitary in $\mathcal{M}$, we
prove that%
\begin{equation*}
W^{\ast }W=P_{ac}^{\infty }\left( \boldsymbol{U}_{0},\mathcal{M}\right) 
\text{ and }WW^{\ast }=P_{ac}^{\infty }\left( \boldsymbol{U},\mathcal{M}%
\right) .
\end{equation*}

For $d=1,$ we prove an analogue of the Kato-Rosenblum theorem by a
stationary method. More precisely, if $U$, $U_{0}$ are unitaries in a
countably decomposable, properly infinite von Neumann algebra $\mathcal{M}$
with a faithful normal semifinite tracial weight $\tau $ and $J$ is a
unitary in $\mathcal{M}$ satisfying $JU_{0}-UJ\in \mathcal{L}^{1}\left( 
\mathcal{M},\tau \right) \mathcal{\cap M},$ then the strong wave operators $%
W_{\pm }:=W_{\pm }\left( U,U_{0},J,\mathcal{M}\right) $ exist and 
\begin{equation*}
W_{\pm }=w_{\pm }\left( U,U_{0},J,\mathcal{M}\right)
\end{equation*}%
where $w_{\pm }\left( U,U_{0},J,\mathcal{M}\right) $ are the corresponding
stationary wave operators. The two signs are treated separately. We also
establish the identities 
\begin{equation*}
W_{\pm }^{\ast }W_{\pm }=P_{ac}^{\infty }\left( U_{0},\mathcal{M}\right) 
\text{ and }W_{\pm }W_{\pm }^{\ast }=P_{ac}^{\infty }\left( U,\mathcal{M}%
\right) .
\end{equation*}
Consequently, no unitary $U$ in $\mathcal{M}$ with $P_{ac}^{\infty }\left( U,%
\mathcal{M}\right) \neq 0$ can be written as $U=D+K,$ where $D\in \mathcal{M}
$ is a diagonal unitary and $K\in \mathcal{L}^{1}\left( \mathcal{M},\tau
\right) \mathcal{\cap M}.$

The paper begins with Section 2. In this section, we introduce the concept
of sector product resolvents of $d$-tuples of commuting unitaries. We also
recall the notions of the absolutely continuous spectrum and the normed
ideal $\mathcal{L}^{1}\left( \mathcal{M},\tau \right) \cap \mathcal{M}$ of a
von Neumann algebra $\mathcal{M}$. In Section 3, motivated by the smoothness
defined in \cite{AMA}, we introduce the notion of jointly smooth operators
associated with $d$-tuples of commuting unitaries in von Neumann algebras
and establish equivalent characterizations of joint smoothness. Building on
this notion, we define the joint norm absolutely continuous support $%
P_{ac}^{\infty }\left( \boldsymbol{U},\mathcal{M}\right) ,$ which extends
the notion of the absolutely continuous projection $P_{ac}\left( \boldsymbol{%
U}\right) $ in the setting of von Neumann algebras and establish its basic
properties. We also provide an example showing that $P_{ac}^{\infty }\left( 
\boldsymbol{U},\mathcal{M}\right) $ can be strictly smaller than $%
P_{ac}\left( \boldsymbol{U}\right) $. In Section 4, we introduce the notions
of strong and stationary wave operators for two $d$-tuples of commuting
unitaries in a von Neumann algebra. For $d\geq 2$, under joint smoothness
assumptions on the factorizations $JU_{0,k}-U_{k}J=G_{k}^{\ast }G_{0,k}$, we
use the stationary method to establish the existence of the strong wave
operator and its coincidence with the corresponding stationary wave
operators. Section 5 is devoted to the case $d=1.$ In this section, an
analogue of the Kato-Rosenblum theorem for a pair of unitaries $U$ and $%
U_{0} $ in $\mathcal{M}$ is given by applying a stationary method.

\section{ Preliminaries and notation}

Let $\mathcal{H}$ be a complex Hilbert space and let $\mathcal{B(H)}$ be the
algebra of all bounded linear operators on $\mathcal{H}.$

\subsection{The sector product resolvent of $\boldsymbol{U}$}

We denote the one-dimensional torus $\mathbb{R}/2\pi \mathbb{Z}$ by $\mathbb{%
T}$ and identify it with the half-open interval $[0,2\pi )$. The unit circle
and the open unit disk of the complex plane are respectively denoted by $%
\mathbb{S}$ and $\mathbb{D}.$ If $\Theta $ is a subset of $\mathbb{T},$ then
the corresponding subarc of $\mathbb{S},$ $\exp \left( i\Theta \right) ,$ is
written $\widehat{\Theta }$ and any function $\Phi $ defined on $\mathbb{S}$
is associated in a unique manner to the function $\phi ,$ defined on $%
\mathbb{T}$ by $\phi \left( \cdot \right) =\Phi \left( e^{i\cdot }\right) .$

Let $U$ be a unitary operator on $\mathcal{H}.$ By the spectral theorem,
there is a spectral measure $F_{U}$ defined on the Borel subsets of $\mathbb{%
S}$ such that 
\begin{equation*}
U=\int_{\mathbb{S}}zdF_{U}\left( z\right) .
\end{equation*}%
If we set $E_{U}\left( t\right) =F_{U}\left( \delta _{t}\right) $ where $%
\delta _{t}=\left\{ z=\exp \left( i\tau \right) :0\leq \tau \leq t\right\} ,$
then 
\begin{equation*}
U=\int_{\mathbb{T}}e^{it}dE_{U}\left( t\right) .
\end{equation*}%
The non-decreasing function $E_{U}\left( t\right) $ can be extended to $%
\mathbb{R}$ by setting $E_{U}\left( t\right) =0$ as $t<0$ and $E_{U}\left(
t\right) =I$ as $t\geq 2\pi .$ Throughout the remainder of the paper, we
always assume $E_{U}$ is defined on $\mathbb{R}.$ If $\varphi $ is a
function, then $\varphi (U)=\int_{\mathbb{S}}\varphi \left( z\right)
dF_{U}\left( z\right) =\int_{\mathbb{T}}\varphi \left( e^{it}\right)
dE_{U}\left( t\right) .$

The resolvent of $U$ is defined as 
\begin{equation*}
R\left( U,z\right) :=\left( 1-zU^{\ast }\right) ^{-1},\text{ \ \ \ \ }z\in 
\mathbb{C}\backslash \mathbb{S}
\end{equation*}%
and then we define 
\begin{equation*}
\delta \left( U,z\right) :=\frac{1-\left\vert z\right\vert ^{2}}{2\pi }%
R\left( U,z\right) ^{\ast }R\left( U,z\right) ,\text{ \ \ }z\in \mathbb{C}%
\backslash \mathbb{S.}
\end{equation*}%
It can be verified that 
\begin{eqnarray*}
R\left( U,re^{i\theta }\right) -R\left( U,r^{-1}e^{i\theta }\right)
&=&\left( 1-r^{2}\right) R\left( U,re^{i\theta }\right) ^{\ast }R\left(
U,re^{i\theta }\right) \\
&=&(1-r^{2})\left\vert R\left( U,re^{i\theta }\right) \right\vert ^{2}\text{
for }r\in (0,\infty )\backslash \left\{ 1\right\} .
\end{eqnarray*}%
In the remainder of the paper, we denote $R\left( U,re^{i\theta }\right) $
by $R^{+}\left( U,r,\theta \right) $ and $R\left( U,r^{-1}e^{i\theta
}\right) $ by $R^{-}\left( U,r,\theta \right) $ for $r\in \left( 0,1\right)
. $ Then 
\begin{equation*}
R^{+}\left( U,r,\theta \right) =\sum_{m=0}^{\infty }r^{m}e^{im\theta }U^{-m}.
\end{equation*}%
and%
\begin{equation*}
R^{-}\left( U,r,\theta \right) =-\sum_{n=1}^{\infty }r^{n}e^{-in\theta
}U^{n}.\text{ }
\end{equation*}%
Let%
\begin{equation}
\gamma _{+}\left( r\right) :=\frac{1-r^{2}}{2\pi }>0\text{ and }\gamma
_{-}\left( r\right) :=\frac{r^{-2}-1}{2\pi }  \label{p9}
\end{equation}%
for $r\in \left( 0,1\right) .$ If we denote $\delta \left( U,re^{i\theta
}\right) $ by $\delta \left( U,r,\theta \right) $ for $r\in \left(
0,1\right) ,$ then 
\begin{eqnarray*}
\delta \left( U,r,\theta \right) &=&\frac{1}{2\pi }\left( R^{+}\left(
U,r,\theta \right) -R^{-}\left( U,r,\theta \right) \right) \\
&=&\gamma _{\pm }\left( r\right) R^{\pm }\left( U,r,\theta \right) ^{\ast
}R^{\pm }\left( U,r,\theta \right) .
\end{eqnarray*}%
In particular, for $r\in (0,1)$ and $\theta $ in $\mathbb{T}$, we have 
\begin{equation}
\delta \left( U,r,\theta \right) =\frac{1}{2\pi }\dsum\limits_{n\in \mathbb{Z%
}}r^{\left\vert n\right\vert }e^{in\theta }U^{\ast n}=\frac{1}{2\pi }%
\dint\limits_{\mathbb{T}}P_{r}\left( \theta -\theta ^{\prime }\right)
dE_{U}\left( \theta ^{\prime }\right)  \label{p8}
\end{equation}%
where the Poisson kernel $P_{r}\left( \cdot \right) $ is defined as follows 
\begin{eqnarray*}
P_{r}\left( \theta -\theta ^{\prime }\right) &=&\frac{1-r^{2}}{\left\vert
1-re^{i\left( \theta -\theta ^{\prime }\right) }\right\vert ^{2}}%
=\dsum\limits_{n\in \mathbb{Z}}r^{\left\vert n\right\vert }e^{in\left(
\theta -\theta ^{\prime }\right) } \\
&=&\frac{1-r^{2}}{1-2r\cos \left( \theta -\theta ^{\prime }\right) +r^{2}}>0
\end{eqnarray*}%
for all $\theta ^{\prime }\in \mathbb{T}$ and 
\begin{equation*}
\int_{0}^{2\pi }\delta \left( U,r,\theta \right) d\theta =I
\end{equation*}%
for $r\in (0,1)$ by applying the results in Chapter 2 of \cite{SS}$.$ We
also recall from (2.5) in \cite{TR} that for any $\varphi $ and $\psi $ in $%
\mathcal{H}$ we have the relations 
\begin{equation*}
\lim_{\varepsilon \searrow 0}\left\langle \delta \left( U,\left(
1-\varepsilon \right) ^{\pm },\theta \right) \varphi ,\psi \right\rangle
=\pm \frac{d}{d\theta }\left\langle E_{U}\left( (0,\theta ]\right) \varphi
,\psi \right\rangle
\end{equation*}%
for a.e. $\theta \in \lbrack 0,2\pi ).$

For 
\begin{equation*}
\boldsymbol{z}=\left( z_{1},\cdots ,z_{d}\right) ,\text{ }\boldsymbol{\omega 
}=\left( \omega _{1},\cdots ,\omega _{d}\right) \text{ in }\mathbb{C}^{d},
\end{equation*}%
we define $\left\langle \boldsymbol{z},\boldsymbol{\omega }\right\rangle
=z_{1}\overline{\omega }_{1}+\cdots +z_{d}\overline{\omega }_{d}$ where $%
\overline{\omega }_{k}$ is the complex conjugate of $\omega _{k}.$ We also
write $\left\vert \boldsymbol{z}\right\vert =\sqrt{\left\langle \boldsymbol{z%
},\boldsymbol{z}\right\rangle }.$

The open unit polydisc in $\mathbb{C}^{d}$ is the set 
\begin{equation*}
\mathbb{D}^{d}=\left\{ \boldsymbol{z}\in \mathbb{C}^{d}:\left\vert
z_{j}\right\vert <1,j=1,\cdots ,d\right\} .
\end{equation*}%
Its distinguished boundary is 
\begin{equation*}
\mathbb{S}^{d}=\left\{ \boldsymbol{z}\in \mathbb{C}^{d}:\left\vert
z_{j}\right\vert =1\right\} ,
\end{equation*}%
which is identified with 
\begin{equation*}
\mathbb{T}^{d}=[0,2\pi )^{d}.
\end{equation*}%
If $\boldsymbol{t},\boldsymbol{n}\in \mathbb{Z}^{d}$ with $\boldsymbol{t}%
=\left( t_{1},\cdots ,t_{d}\right) $ and $\boldsymbol{n=}\left( n_{1},\cdots
,n_{d}\right) ,$ then 
\begin{equation*}
\boldsymbol{tn=}\left( t_{1}n_{1},\cdots ,t_{d}n_{d}\right) .
\end{equation*}%
For a commuting $d$-tuple $\boldsymbol{U}=\left( U_{1},\cdots ,U_{d}\right) $
of unitaries in $\mathcal{B}\left( \mathcal{H}\right) ,$ the spectral
measure of $\boldsymbol{U}$ is denoted by $E_{\boldsymbol{U}}:\mathcal{B}%
\left( \mathbb{T}^{d}\right) \rightarrow \mathcal{P}\left( \mathcal{B}\left( 
\mathcal{H}\right) \right) ,$ where $\mathcal{B}\left( \mathbb{T}^{d}\right) 
$ represents the $\sigma $-algebra of all Borel sets in $\mathbb{T}^{d}.$ We
define a Borel function $f_{\theta _{i}}:\mathbb{T}^{d}\rightarrow \mathbb{C}
$ for $i=1,\cdots ,d$ such that for $\boldsymbol{\theta =}\left( \theta
_{1},\cdots ,\theta _{d}\right) \in \mathbb{T}^{d}$ 
\begin{equation*}
f_{\theta _{i}}\left( \boldsymbol{\theta }\right) =e^{\theta _{i}}.
\end{equation*}%
Then $U_{i}=\int_{\mathbb{T}^{d}}f_{\theta _{i}}\left( \boldsymbol{\theta }%
\right) dE_{\boldsymbol{U}}\left( \boldsymbol{\theta }\right) $ for $%
i=1,\cdots ,d$. We also define 
\begin{equation*}
\boldsymbol{U}^{\boldsymbol{n}}=U_{1}^{n_{1}}\cdots U_{d}^{n_{d}}\text{ for }%
\boldsymbol{n}\in \mathbb{Z}^{d}.
\end{equation*}

Hence, we define%
\begin{equation*}
\Delta \left( \boldsymbol{U},\boldsymbol{z}\right)
:=\dprod\limits_{j=1}^{d}\delta \left( U_{j},z_{j}\right)
\end{equation*}%
for $\boldsymbol{z=}\left( z_{1},\cdots ,z_{d}\right) .$ If $%
z_{j}=r_{j}e^{i\theta _{j}}$ and $r_{j}\in \left( 0,1\right) ,$ then 
\begin{equation*}
\Delta \left( \boldsymbol{U},\boldsymbol{r},\boldsymbol{\theta }\right)
:=\Delta \left( \boldsymbol{U},\boldsymbol{re}^{i\boldsymbol{\theta }%
}\right) =\dprod\limits_{j=1}^{d}\delta \left( U_{j},r_{j},\theta
_{j}\right) .
\end{equation*}%
Let $R_{j}^{+}=R^{+}\left( U_{j},r_{j},\theta _{j}\right) $ and $%
R_{j}^{-}=R^{-}\left( U_{j},r_{j},\theta _{j}\right) $ for $r_{j}\in \left(
0,1\right) .$ In \cite{O}, 
\begin{eqnarray*}
P_{\boldsymbol{r}}\left( \boldsymbol{\theta -t}\right)
&:&=\dprod\limits_{j=1}^{d}P_{r_{j}}\left( \theta _{j}-t_{j}\right)
=\dprod\limits_{j=1}^{d}\frac{1-r_{j}^{2}}{\left\vert 1-r_{j}e^{i\left(
\theta _{j}-t_{j}\right) }\right\vert ^{2}} \\
&=&\sum_{\boldsymbol{m}\in \mathbb{Z}^{d}}\boldsymbol{r}^{\left\vert 
\boldsymbol{m}\right\vert }e^{i\left\langle \boldsymbol{m},\boldsymbol{%
\theta -t}\right\rangle }
\end{eqnarray*}%
for $\boldsymbol{r=}\left( r_{1},\cdots ,r_{d}\right) \in \left( 0,1\right)
^{d}$, $\boldsymbol{\theta =}\left( \theta _{1},\cdots ,\theta _{d}\right) $
and $\boldsymbol{t}=\left( t_{1},\cdots ,t_{d}\right) $ where $\boldsymbol{r}%
^{\left\vert \boldsymbol{m}\right\vert
}=\dprod\limits_{j=1}^{d}r_{j}^{\left\vert m_{j}\right\vert }$ (see \cite{O})%
$.$ Then, by \cite{G} and \cite{P}, for $F\in \mathcal{L}^{2}\left( \mathbb{T%
}^{d},\mathcal{H}\right) ,$ we have 
\begin{eqnarray}
P_{\boldsymbol{r}}\ast F\left( \boldsymbol{\theta }\right) &=&\frac{1}{%
\left( 2\pi \right) ^{d}}\int_{\mathbb{T}^{d}}\dprod%
\limits_{j=1}^{d}P_{r_{j}}\left( \theta _{j}-t_{j}\right) F\left( 
\boldsymbol{t}\right) d\boldsymbol{t}  \label{p21} \\
&=&\sum_{\boldsymbol{m}\in \mathbb{Z}^{d}}\boldsymbol{r}^{\left\vert 
\boldsymbol{m}\right\vert }e^{i\left\langle \boldsymbol{m},\boldsymbol{%
\theta }\right\rangle }\left[ \frac{1}{\left( 2\pi \right) ^{d}}\int_{%
\mathbb{T}^{d}}e^{-i\left\langle \boldsymbol{m},\boldsymbol{t}\right\rangle
}F\left( \boldsymbol{t}\right) d\boldsymbol{t}\right]  \notag \\
&=&\sum_{\boldsymbol{m}\in \mathbb{Z}^{d}}\boldsymbol{r}^{\left\vert 
\boldsymbol{m}\right\vert }e^{i\left\langle \boldsymbol{m},\boldsymbol{%
\theta }\right\rangle }\widehat{F}\left( \boldsymbol{m}\right)  \notag
\end{eqnarray}%
and 
\begin{equation*}
\widehat{P_{\boldsymbol{r}}\ast F}\left( \boldsymbol{m}\right) =\boldsymbol{r%
}^{\left\vert \boldsymbol{m}\right\vert }\widehat{F}\left( \boldsymbol{m}%
\right)
\end{equation*}%
where $\widehat{F}\left( \boldsymbol{m}\right) =\frac{1}{\left( 2\pi \right)
^{d}}\int_{\mathbb{T}^{d}}e^{-i\left\langle \boldsymbol{m},\boldsymbol{t}%
\right\rangle }F\left( \boldsymbol{t}\right) d\boldsymbol{t}.$ It implies
that 
\begin{equation}
\Delta \left( \boldsymbol{U},\boldsymbol{r},\boldsymbol{\theta }\right)
=\left( \frac{1}{2\pi }\right) ^{d}\dint\limits_{\mathbb{T}^{d}}P_{%
\boldsymbol{r}}\left( \boldsymbol{\theta }\mathbf{-}\boldsymbol{\lambda }%
\right) dE_{\boldsymbol{U}}\left( \boldsymbol{\lambda }\right)  \label{p17}
\end{equation}
and 
\begin{equation}
\int_{\mathbb{T}^{d}}\Delta \left( \boldsymbol{U},\boldsymbol{r},\boldsymbol{%
\theta }\right) d\boldsymbol{\theta }=I.  \label{p1}
\end{equation}

Now, we define the sector product resolvent 
\begin{equation*}
R_{\boldsymbol{\sigma }}\left( \boldsymbol{U},\boldsymbol{r},\boldsymbol{%
\theta }\right) :=\dprod\limits_{j=1}^{d}R_{j}^{\sigma _{j}}
\end{equation*}%
for $\boldsymbol{\sigma }\in \left\{ +1,-1\right\} ^{d}.$ Since $\delta _{j}=%
\frac{1}{2\pi }\left( R_{j}^{+}-R_{j}^{-}\right) ,$ we conclude that for $%
\boldsymbol{r=}\left( r_{1},\cdots ,r_{d}\right) \in \left( 0,1\right) ^{d}$ 
\begin{equation*}
\Delta \left( \boldsymbol{U},\boldsymbol{r},\boldsymbol{\theta }\right)
=\dprod\limits_{j=1}^{d}\delta _{j}=\left( \frac{1}{2\pi }\right) ^{d}\sum_{%
\boldsymbol{\sigma }\in \left\{ +1,-1\right\} ^{d}}\left( -1\right)
^{N_{-}\left( \boldsymbol{\sigma }\right) }R_{\boldsymbol{\sigma }}\left( 
\boldsymbol{U},\boldsymbol{r},\boldsymbol{\theta }\right) ,
\end{equation*}%
where $N_{-}\left( \boldsymbol{\sigma }\right) $ is the number of negative
components of $\boldsymbol{\sigma }.$

The following lemma extends a classical result to the setting of functions
on $\mathbb{T}^{d}$ with values in a Hilbert space $\mathcal{H}$. For
background on Bochner integration and the space $\mathcal{L}^{2}\left( 
\mathbb{T}^{d},\mathcal{H}\right) ,$ we refer the reader to \cite{D} and 
\cite{Arendt1}.

\begin{lemma}
\label{lem0}Let $\mathcal{H}$ be a Hilbert space and $F\in \mathcal{L}%
^{2}\left( \mathbb{T}^{d},\mathcal{H}\right) $ and let 
\begin{equation*}
P^{\ast }F\left( \boldsymbol{\theta }\right) :=\sup_{\boldsymbol{r}\in
\left( 0,1\right) ^{d}}\left\Vert P_{\boldsymbol{r}}\ast F\left( \boldsymbol{%
\theta }\right) \right\Vert .
\end{equation*}%
Then there exists a constant $C_{d}>0$ such that 
\begin{equation*}
\left\Vert P^{\ast }F\right\Vert _{\mathcal{L}^{2}\left( \mathbb{T}%
^{d}\right) }\leq C_{d}\left\Vert F\right\Vert _{\mathcal{L}^{2}\left( 
\mathbb{T}^{d},\mathcal{H}\right) }.
\end{equation*}%
Moreover, 
\begin{equation*}
\lim_{\underset{j=1,\cdots ,d}{r_{j}\rightarrow 1^{-}}}\left\Vert \left( P_{%
\boldsymbol{r}}\ast F\right) \left( \boldsymbol{\theta }\right) -F\left( 
\boldsymbol{\theta }\right) \right\Vert _{\mathcal{H}}=0\text{ for a.e. }%
\boldsymbol{\theta }\in \mathbb{T}^{d},
\end{equation*}%
and 
\begin{equation*}
\lim_{\underset{j=1,\cdots ,d}{r_{j}\rightarrow 1^{-}}}\left\Vert P_{%
\boldsymbol{r}}\ast F-F\right\Vert _{\mathcal{L}^{2}\left( \mathbb{T}^{d},%
\mathcal{H}\right) }=0.
\end{equation*}
\end{lemma}

\begin{proof}
For each $j=1,\cdots ,d,$ let $M_{j}$ denote the one-dimensional
Hardy-Littlewood maximal operator acting in the $j$-th coordinate. Thus, for
a scalar function $g$ acting on $\mathbb{T}^{d},$%
\begin{equation*}
M_{j}g\left( \boldsymbol{\theta }\right) =\sup_{0\leq h\leq \pi }\frac{1}{2h}%
\int_{-h}^{h}\left\vert g\left( \theta _{1},\cdots ,\theta _{j}-t,\cdots
,\theta _{d}\right) \right\vert dt
\end{equation*}%
and 
\begin{equation}
\left\Vert M_{j}g\right\Vert _{\mathcal{L}^{2}}\leq \widetilde{C}\left\Vert
g\right\Vert _{\mathcal{L}^{2}}  \label{3}
\end{equation}%
where $\widetilde{C}$ is a constant independent of $j.$ We define the
Poisson operator acting in the $j$-th coordinate below%
\begin{equation*}
\left( P_{j,r_{j}}\ast g\right) \left( \boldsymbol{\theta }\right) :=\frac{1%
}{2\pi }\int_{\mathbb{T}}P_{r_{j}}\left( \theta _{j}-t\right) g\left( \theta
_{1},\cdots ,t,\cdots ,\theta _{d}\right) dt
\end{equation*}%
The one-dimensional Poisson maximal estimate implies that, for every
nonnegative function $g$ 
\begin{equation*}
P_{j,r_{j}}\ast g\leq CM_{j}g\text{ for }0<r_{j}<1
\end{equation*}%
where $C$ is independent of $j,r_{j}$ and $g$. If we let $g\left( 
\boldsymbol{\theta }\right) =\left\Vert F\left( \boldsymbol{\theta }\right)
\right\Vert _{\mathcal{H}}$ for a function $F\in \mathcal{L}^{2}\left( 
\mathbb{T}^{d},\mathcal{H}\right) ,$ then 
\begin{eqnarray*}
\left\Vert P_{\boldsymbol{r}}\ast F\left( \boldsymbol{\theta }\right)
\right\Vert &=&\frac{1}{\left( 2\pi \right) ^{d}}\left\Vert \int_{\mathbb{T}%
^{d}}P_{\boldsymbol{r}}\left( \boldsymbol{\theta -t}\right) F\left( 
\boldsymbol{t}\right) d\lambda ^{d}\left( \boldsymbol{t}\right) \right\Vert
\\
&\leq &\frac{1}{\left( 2\pi \right) ^{d}}\int_{\mathbb{T}^{d}}P_{\boldsymbol{%
r}}\left( \boldsymbol{\theta -t}\right) \left\Vert F\left( \boldsymbol{t}%
\right) \right\Vert d\lambda ^{d}\left( \boldsymbol{t}\right) =P_{%
\boldsymbol{r}}\ast g\left( \boldsymbol{\theta }\right) \text{ a.e. }\theta .
\end{eqnarray*}%
Applying the preceding one-dimensional estimate successively in each
coordinate yields%
\begin{equation*}
P_{\boldsymbol{r}}\ast g\left( \boldsymbol{\theta }\right) =P_{1,r_{1}}\ast
\cdots \ast P_{d,r_{d}}\ast g\leq C^{d}M_{1}\cdots M_{d}g.
\end{equation*}%
Then%
\begin{eqnarray*}
P^{\ast }F\left( \boldsymbol{\theta }\right) &=&\sup_{\boldsymbol{r}\in
\left( 0,1\right) ^{d}}\left\Vert \left( P_{\boldsymbol{r}}\ast F\right)
\left( \boldsymbol{\theta }\right) \right\Vert _{\mathcal{H}}\leq \sup_{%
\boldsymbol{r}\in \left( 0,1\right) ^{d}}P_{\boldsymbol{r}}\ast g\left( 
\boldsymbol{\theta }\right) \\
&\leq &C^{d}M_{1}\cdots M_{d}g\left( \boldsymbol{\theta }\right)
=C^{d}M_{1}\cdots M_{d}\left\Vert F\left( \boldsymbol{\theta }\right)
\right\Vert _{\mathcal{H}}.
\end{eqnarray*}%
Consequently, by (\ref{3}), there is a constant $C_{d}$ such that 
\begin{equation}
\left\Vert P^{\ast }F\right\Vert _{\mathcal{L}^{2}\left( \mathbb{T}%
^{d}\right) }\leq C_{d}\left\Vert F\right\Vert _{\mathcal{L}^{2}\left( 
\mathbb{T}^{d},\mathcal{H}\right) }.  \label{1.10}
\end{equation}%
Now we choose $\mathcal{H}$-valued trigonometric polynomials $F_{N}$ such
that 
\begin{equation}
\left\Vert F-F_{N}\right\Vert _{\mathcal{L}^{2}\left( \mathbb{T}^{d},%
\mathcal{H}\right) }\rightarrow 0\text{ as }N\rightarrow \infty .
\label{1.1}
\end{equation}%
Then for each $N,$ the trigonometric polynomial $F_{N}$ can be expressed as 
\begin{equation*}
F_{N}\left( \boldsymbol{\theta }\right) =\sum_{\boldsymbol{n}\in A_{N}}a_{%
\boldsymbol{n}}e^{i\left\langle \boldsymbol{n},\boldsymbol{\theta }%
\right\rangle }\text{ where }A_{N}\subset \mathbb{Z}^{d}\text{ is finite.}
\end{equation*}%
Hence 
\begin{equation*}
\left( P_{\boldsymbol{r}}\ast F_{N}\right) \left( \boldsymbol{\theta }%
\right) =\sum_{\boldsymbol{n}\in A_{N}}\left( \Pi
_{j=1}^{d}r_{j}^{\left\vert n_{j}\right\vert }\right) a_{\boldsymbol{n}%
}e^{i\left\langle \boldsymbol{n},\boldsymbol{\theta }\right\rangle }.
\end{equation*}%
It follows that 
\begin{equation}
\left\Vert P_{\boldsymbol{r}}\ast F_{N}-F_{N}\right\Vert _{\mathcal{L}%
^{2}\left( \mathbb{T}^{d},\mathcal{H}\right) }\rightarrow 0\text{ }
\label{1.2}
\end{equation}%
as $r_{j}\rightarrow 1^{-}$ for every $j.$ We define 
\begin{equation*}
L\left( \boldsymbol{\theta }\right) =\underset{r_{j}\rightarrow 1^{-}\text{
for }j=1,\cdots ,d}{\lim \sup }\left\Vert \left( P_{\boldsymbol{r}}\ast
F\right) \left( \boldsymbol{\theta }\right) -F\left( \boldsymbol{\theta }%
\right) \right\Vert _{\mathcal{H}}.
\end{equation*}%
Then 
\begin{equation*}
L\left( \boldsymbol{\theta }\right) \leq P^{\ast }\left( F-F_{N}\right)
\left( \boldsymbol{\theta }\right) +\left\Vert \left( F-F_{N}\right) \left( 
\boldsymbol{\theta }\right) \right\Vert _{\mathcal{H}}
\end{equation*}%
a.e. $\boldsymbol{\theta }.$ Therefore, by Chebyshev's inequality and the
maximal estimate, for every $\varepsilon >0,$%
\begin{eqnarray*}
\lambda ^{d}\left\{ \boldsymbol{\theta }:L\left( \boldsymbol{\theta }\right)
>\varepsilon \right\} &\leq &\frac{1}{\varepsilon ^{2}}\left( \left\Vert
P^{\ast }\left( F-F_{N}\right) \right\Vert _{\mathcal{L}^{2}\left( \mathbb{T}%
^{d}\right) }+\left\Vert F-F_{N}\right\Vert _{\mathcal{L}^{2}\left( \mathbb{T%
}^{d},\mathcal{H}\right) }\right) ^{2} \\
&\leq &\frac{\left( C_{d}+1\right) ^{2}}{\varepsilon ^{2}}\left\Vert
F-F_{N}\right\Vert _{\mathcal{L}^{2}\left( \mathbb{T}^{d},\mathcal{H}\right)
}^{2}.
\end{eqnarray*}%
Then we obtain that $L\left( \boldsymbol{\theta }\right) =0$ a.e. $%
\boldsymbol{\theta }.$ This proves that 
\begin{equation}
\lim_{\underset{j=1,\cdots ,d}{r_{j}\rightarrow 1^{-}}}\left\Vert \left( P_{%
\boldsymbol{r}}\ast F\right) \left( \boldsymbol{\theta }\right) -F\left( 
\boldsymbol{\theta }\right) \right\Vert _{\mathcal{H}}=0\text{ for a.e. }%
\boldsymbol{\theta }\in \mathbb{T}^{d}.  \label{1.11}
\end{equation}%
Since by (\ref{1.1}) and (\ref{1.2}), 
\begin{eqnarray*}
&&\left\Vert P_{\boldsymbol{r}}\ast F-F\right\Vert _{\mathcal{L}^{2}\left( 
\mathbb{T}^{d},\mathcal{H}\right) } \\
&=&\left\Vert P_{\boldsymbol{r}}\ast \left( F-F_{N}\right) +P_{\boldsymbol{r}%
}\ast F_{N}-F_{N}+F_{N}-F\right\Vert _{\mathcal{L}^{2}\left( \mathbb{T}^{d},%
\mathcal{H}\right) } \\
&\leq &2\left\Vert F-F_{N}\right\Vert _{\mathcal{L}^{2}\left( \mathbb{T}^{d},%
\mathcal{H}\right) }+\left\Vert P_{\boldsymbol{r}}\ast
F_{N}-F_{N}\right\Vert _{\mathcal{L}^{2}\left( \mathbb{T}^{d},\mathcal{H}%
\right) }\text{ for every }N,
\end{eqnarray*}%
we obtain that%
\begin{equation}
\lim_{\underset{j=1,\cdots ,d}{r_{j}\rightarrow 1^{-}}}\left\Vert P_{r}\ast
F-F\right\Vert _{\mathcal{L}^{2}\left( \mathbb{T}^{d},\mathcal{H}\right) }=0.
\label{1.12}
\end{equation}%
Combining (\ref{1.10}), (\ref{1.11}) and (\ref{1.12}), the proof is
completed.
\end{proof}

\subsection{Absolutely continuous subspaces}

Let $H$ be a possibly unbounded self-adjoint element. If the measure $m_{x}$
defined by $m_{x}\left( S\right) =\left\langle E_{H}\left( S\right)
x,x\right\rangle $ for any Borel set $S$ is absolutely continuous with
respect to Lebesgue measure $\lambda ,$ then $x$ is called absolutely
continuous. Let 
\begin{equation*}
\mathcal{H}_{ac}\left( H\right) =\left\{ x:m_{x}\ll \text{ }\lambda \right\}
.
\end{equation*}%
We know that $\mathcal{H}_{ac}\left( H\right) $ is a closed subspace of $%
\mathcal{H}$ and the projection $P_{ac}\left( H\right) $ from $\mathcal{H}$
onto $\mathcal{H}_{ac}\left( H\right) $ commutes with $H$ and there is a
Borel set $\Omega $ such that $P_{ac}\left( H\right) =E_{H}\left( \Omega
\right) $ (see \cite{Siknia} or \cite{TG}). If, on the other hand, the
measure $m_{x}$ is singular with respect to $\lambda $, $x$ is called
singular with respect to $H.$ The set of all singular vectors is denoted by $%
\mathcal{H}_{s}$ and the projection $P_{s}\left( H\right) $ from $\mathcal{H}
$ onto $\mathcal{H}_{s}$ commutes with $H$ as well. We also have $%
P_{s}\left( H\right) +P_{ac}\left( H\right) =I.$ If $P_{s}\left( H\right)
=I, $ we say that $H$ is purely singular.

For a unitary operator $U,$ suppose 
\begin{equation*}
U=\int_{\mathbb{S}}zdF_{U}\left( z\right) =\int_{\mathbb{T}%
}e^{it}dE_{U}\left( t\right) .
\end{equation*}%
Then a vector $x$ is called absolutely continuous if the measure $%
m_{x}\left( B\right) =\left\langle E_{U}\left( B\right) x,x\right\rangle $
defined on the Borel subsets $B$ of $\mathbb{T}$ is absolutely continuous
with respect to Lebesgue measure $\lambda .$ Let $\mathcal{H}_{ac}\left(
U\right) $ be the collection of such $x$ in $\mathcal{H}$ and $P_{ac}\left(
U\right) $ be the orthogonal projection onto $\mathcal{H}_{ac}\left(
U\right) .$

Let $\boldsymbol{U=}\left( U_{1},\cdots ,U_{d}\right) $ be a commuting $d$%
-tuple of unitaries in $\mathcal{B}\left( \mathcal{H}\right) $ and $\lambda
^{d}$ be the standard Lebesgue measure on $\mathbb{T}^{d}.$ Then the
absolutely continuous subspace is defined below 
\begin{equation*}
\mathcal{H}_{ac}\left( \boldsymbol{U}\right) =\left\{ x:m_{x}\left( S\right)
=\left\langle E_{\boldsymbol{U}}\left( S\right) x,x\right\rangle \text{ is }%
\lambda ^{d}\text{-absolutely continuous}\right\} .
\end{equation*}%
Similarly, $\mathcal{H}_{ac}\left( \boldsymbol{U}\right) $ is a closed
subspace of $\mathcal{H}$ and the projection $P_{ac}\left( \boldsymbol{U}%
\right) $ from $\mathcal{H}$ onto $\mathcal{H}_{ac}\left( \boldsymbol{U}%
\right) $ commutes with each $U_{j}.$ It can also be shown that $%
P_{ac}\left( \boldsymbol{U}\right) =E_{\boldsymbol{U}}\left( \Omega \right) $
for a Borel subset $\Omega $ by a similar argument in \cite{TG}. \ 

For each fixed pair $\varphi ,\psi \in \mathcal{H}$, the boundary
convergence theorem for product Poisson integrals of finite measures yields 
\begin{equation}
\lim_{\underset{\rho \longrightarrow 1^{-}}{\boldsymbol{r=}\rho \cdot 
\boldsymbol{1}}}\left\langle \Delta \left( \boldsymbol{U},\boldsymbol{re}^{i%
\boldsymbol{\theta }}\right) \varphi ,\psi \right\rangle =\frac{d\mu
_{\varphi ,\psi }^{ac}}{d\lambda ^{d}}\left( \boldsymbol{\theta }\right) 
\text{ for a.e. }\boldsymbol{\theta }\in \mathbb{T}^{d}  \label{1}
\end{equation}%
where $\mu _{\varphi ,\psi }^{ac}$ is the absolutely continuous part of the
measure $\mu _{\varphi ,\psi }$ with respect to $\lambda ^{d}$, defined by $%
\mu _{\varphi ,\psi }\left( \Omega \right) :=\left\langle E_{\boldsymbol{U}%
}\left( \Omega \right) \varphi ,\psi \right\rangle $ and $\boldsymbol{r=}%
\rho \cdot \boldsymbol{1=}\left( \rho ,\cdots ,\rho \right) $ with $0<\rho
<1.$

\subsection{The normed ideal $\mathcal{L}^{1}\left( \mathcal{M},\protect\tau %
\right) \cap \mathcal{M}$}

\bigskip Let $\mathcal{M\subseteq B(H)}$ be a countably decomposable,
properly infinite von Neumann algebra with a faithful normal semifinite
tracial weight $\tau $. We recall that a tracial weight $\tau $ is called
normal and semifinite if for any bounded increasing net $\left\{
X_{i}\right\} $ in $\mathcal{M}_{+},$ we have $\sup_{i}\tau \left\{
X_{i}\right\} =\tau \left( \sup_{i}X_{i}\right) $ and there is a $Y\in 
\mathcal{M}_{+}$ with $\tau \left( Y\right) <\infty $ for any non-zero
element $X\in \mathcal{M}_{+} $ such that $0<Y\leq X.$

Let $\mathcal{A}\left( \mathcal{M}\right) $ be the set of all densely
defined, closed operators that are affiliated with $\mathcal{M}.$ Let

\begin{equation*}
\mathcal{F(M},\tau \mathcal{)}=\{AEB:E=E^{\ast }=E^{2}\in \mathcal{M}\text{
with }\tau (E)<\infty \text{ and }A,B\in \mathcal{M}\}
\end{equation*}%
be the set of finite-rank operators in $\mathcal{(M},\tau \mathcal{)}.$ For $%
1\leq p<\infty ,$ the noncommutative $\mathcal{L}^{p}\left( \mathcal{M},\tau
\right) $ space is defined by $\mathcal{L}^{p}\left( \mathcal{M},\tau
\right) =\overline{\mathcal{F(M},\tau \mathcal{)}}^{\left\Vert \cdot
\right\Vert _{p}}$ with respect to the norm $\left\Vert \cdot \right\Vert
_{p}$ defined by $\left\Vert A\right\Vert _{p}=\tau \left( \left\vert
A\right\vert ^{p}\right) ^{1/p}.$ We set $\mathcal{L}^{\infty }\left( 
\mathcal{M},\tau \right) =\mathcal{M}$ equipped with the operator norm. Note
that $\mathcal{L}_{p}\left( \mathcal{M},\tau \right) $ is realized as a
space of $\mathcal{A}\left( \mathcal{M}\right) $ for $1\leq p\leq \infty .$
For more details of $\mathcal{L}^{p}\left( \mathcal{M},\tau \right) $, we
refer the reader to \cite{Pisier}. In \cite{Li2}, we know that $\mathcal{J}%
_{1}=\mathcal{L}^{1}\left( \mathcal{M},\tau \right) \cap \mathcal{M}$ is a
normed ideal of $\mathcal{M}$ with the norm defined by $\left\Vert
A\right\Vert _{\mathcal{J}_{1}}=\max \left\{ \left\Vert A\right\Vert
_{1},\left\Vert A\right\Vert \right\} $ for every $A\in $ $\mathcal{L}%
^{1}\left( \mathcal{M},\tau \right) \cap \mathcal{M}.$

By slightly modifying the proof of Lemma 2.1.1 in \cite{Li}, we obtain the
following result.

\begin{lemma}
\label{trace}Let $\mathcal{M}$ be a countably decomposable, properly
infinite von Neumann algebra with a normal semifinite faithful trace $\tau $.

\begin{enumerate}
\item For any projection $P\in \mathcal{M}$, there exists a sequence $%
\left\{ P_{n}\right\} _{n\in \mathbb{N}}$ of pairwise orthogonal projections
in $\ \mathcal{M}$ such that $\tau \left( P_{n}\right) <\infty $ for each $%
n\in \mathbb{N}$ and $\Sigma _{n\in \mathbb{N}}P_{n}=P.$

\item There exists a sequence $\left\{ x_{m}\right\} _{m\in \mathbb{N}\text{ 
}}$ of vectors in $\mathcal{H}$ such that 
\begin{equation*}
\tau \left( X^{\ast }X\right) =\Sigma _{m}\left\langle X^{\ast
}Xx_{m},x_{m}\right\rangle ,\text{ \ \ \ }\forall X\in \mathcal{M}.
\end{equation*}%
Moreover, the linear span of the set $\left\{ A^{\prime }x_{m}:A^{\prime
}\in \mathcal{M}^{\prime }\text{ and }m\in \mathbb{N}\right\} $ is dense in $%
\mathcal{H}.$
\end{enumerate}
\end{lemma}

For an operator $X\in $ $\mathcal{B(H)}$, we denote by $l\left[ X\right] $
the projection onto the subspace $\overline{X(\mathcal{H})}$, which is
called the left support of $X$, i.e., $l\left[ X\right] $ is the smallest
projection $E\in \mathcal{B(H)}$ for which $EX=X.$ Now we are ready to give
the definition of generating projections.

The following lemma follows by adapting the proof of Lemma 3.5 in \cite{Ber2}%
. We omit the proof.

\begin{lemma}
(Lemma 3.5 in \cite{Ber2}) \label{generating 2}Let $\mathcal{M}$ be a von
Neumann algebra and $\boldsymbol{U}=\left( U_{1},\cdots ,U_{d}\right) $ be a
commuting $d$-tuple of unitaries in $\mathcal{M}$. For every projection $%
Q\in \mathcal{M},$ the projection 
\begin{equation}
E_{Q}=\bigvee _{\boldsymbol{\Omega }\in \mathcal{B}\left( \mathbb{T}%
^{d}\right) }l\left[ E_{\boldsymbol{U}}\left( \boldsymbol{\Omega }\right) Q%
\right]  \label{pr5}
\end{equation}
satisfies%
\begin{equation*}
Q\leq E_{Q}\in \mathcal{M\cap A}^{\prime }
\end{equation*}%
where $\mathcal{A}$ is the von Neumann algebra generated by $\left\{
U_{1},\cdots ,U_{d}\right\} $ and $\mathcal{A}^{\prime }$ is the commutant
of $\mathcal{A}$.
\end{lemma}

\section{Joint Smoothness of Commuting $d$-Tuples of Unitaries in $\mathcal{M%
}$}

Let $\mathcal{M\subseteq B(H)}$ be a von Neumann algebra and denote its
commutant by $\mathcal{M}^{\prime }$. The von Neumann algebra $\mathcal{M}%
\cap \mathcal{M}^{\prime }$ is called the center of $\mathcal{M},$ denoted
by $\mathcal{Z}\left( \mathcal{M}\right) .$ Throughout the remainder of this
paper, we denote the strong and weak operator topologies by "sot" and "wot",
respectively.

In this section, we introduce the notion of joint smoothness. Definition \ref%
{def1} and Theorem \ref{thm1.1} below generalize Definition 2.1 and Theorem
2.2 in \cite{AMA} respectively to the setting of $d$-tuples of commuting
unitaries. The notion of joint smoothness introduced here extends the
classical notion of $H$-smoothness for a self-adjoint operator $H$ (\cite{Y}%
).\bigskip

\begin{definition}
\label{def1}Let $B\in \mathcal{B(H)}$ and let $\boldsymbol{U}=\left(
U_{1},\cdots ,U_{d}\right) $ be a commuting $d$-tuple of unitaries in $%
\mathcal{M}$. We also let 
\begin{equation*}
C_{B}:=\frac{1}{\left( 2\pi \right) ^{d}}\sup_{\left\Vert \psi \right\Vert
=1}\sum_{\mathbf{n}\in \mathbb{Z}^{d}}\left\Vert B\boldsymbol{U}^{\mathbf{n}%
}\psi \right\Vert ^{2}.
\end{equation*}%
If $C_{B}<\infty ,$ then $B$ is called jointly $\boldsymbol{U}$-smooth (or $%
U $-smooth when $d=1$) and $C_{B}$ is called the jointly $\boldsymbol{U}$%
-smoothness constant of the operator $B.$ The set of all jointly $%
\boldsymbol{U}$-smooth operators in $\mathcal{M}$ is denoted by $\mathcal{S}%
_{ac}^{\infty }\left( \boldsymbol{U},\mathcal{M}\right) .$
\end{definition}

Let 
\begin{equation*}
Q=\dprod\limits_{j=1}^{d}[a_{j},b_{j})\subseteq \mathbb{T}^{d}
\end{equation*}%
with $0\leq a_{j}<b_{j}<2\pi $ and $\lambda ^{d}\left( Q\right)
=\dprod\limits_{j=1}^{d}\left( b_{j}-a_{j}\right) $ where $\lambda ^{d}$ is
the standard Lebesgue measure on $\mathbb{T}^{d}.$ The following theorem
provides equivalent characterizations of joint smoothness.

\begin{theorem}
\label{thm1.1}Let $B\in \mathcal{B(H)}$ and $\boldsymbol{U}=\left(
U_{1},\cdots ,U_{d}\right) $ be a commuting $d$-tuple of unitaries in $%
\mathcal{M}.$ Then the following statements are equivalent

\begin{enumerate}
\item $C_{1}=\frac{1}{\left( 2\pi \right) ^{d}}\sup_{\left\Vert \psi
\right\Vert =1}\sum_{\mathbf{n}\in \mathbb{Z}^{d}}\left\Vert B\boldsymbol{U}%
^{\mathbf{n}}\psi \right\Vert ^{2}<\infty ,$

\item $C_{2}=\sup_{\left\Vert \psi \right\Vert =1,\boldsymbol{z\in }\mathbb{D%
}^{d}}\left\langle \Delta \left( \boldsymbol{U},\boldsymbol{z}\right)
B^{\ast }\psi ,B^{\ast }\psi \right\rangle <\infty ,$

\item $C_{3}=\underset{\left\Vert \psi \right\Vert =1,Q=[a_{1},b_{1})\times
\cdots \times \lbrack a_{d},b_{d}\subseteq \mathbb{T}^{d}}{\sup }\frac{%
\left\Vert BE_{\mathbf{U}}\left( Q\right) \psi \right\Vert ^{2}}{\lambda
^{d}\left( Q\right) }<\infty ,$

\item $C_{4}=\underset{\left\Vert \psi \right\Vert =1,Q=[a_{1},b_{1})\times
\cdots \times \lbrack a_{d},b_{d}\subseteq \mathbb{T}^{d}}{\sup }\frac{%
\left\Vert E_{\mathbf{U}}\left( Q\right) B^{\ast }\psi \right\Vert ^{2}}{%
\lambda ^{d}\left( Q\right) }<\infty ,$

\item $C_{5}=\frac{1}{\left( 2\pi \right) ^{d}}\sup_{\left\Vert \psi
\right\Vert =1,\boldsymbol{z\in }\mathbb{D}^{d}}\dprod\limits_{j=1}^{d}%
\left( 1-\left\vert z_{j}\right\vert ^{2}\right) \left\Vert \left(
\dprod\limits_{j=1}^{d}R\left( U_{j},z_{j}\right) \right) B^{\ast }\psi
\right\Vert ^{2}<\infty .$
\end{enumerate}

If these conditions hold, then 
\begin{equation*}
C_{1}=C_{2}=C_{3}=C_{4}=C_{5}=C_{B}.
\end{equation*}
\end{theorem}

For clarity, we organize the proof of Theorem \ref{thm1.1} into several
lemmas. All the lemmas below are stated under the hypotheses of Theorem \ref%
{thm1.1}. The first lemma is obvious, so we omit its proof.

\begin{lemma}
\label{lem1.1}$C_{3}=C_{4}.$
\end{lemma}

\begin{lemma}
\label{lem1.2}If at least one of $C_{1}$ and $C_{4}$ is finite, then $%
C_{1}=C_{4}.$
\end{lemma}

\begin{proof}
First, we assume that $C_{4}$ is finite. For $y\in \mathcal{H}$ and any
Borel set in $\mathbb{T}^{d}$, we define the positive scalar measure 
\begin{equation*}
\mu _{y}\left( \Omega \right) =\left\langle E_{\boldsymbol{U}}\left( \Omega
\right) B^{\ast }y,B^{\ast }y\right\rangle =\left\Vert E_{\boldsymbol{U}%
}\left( \Omega \right) B^{\ast }y\right\Vert ^{2}.
\end{equation*}%
From the definition of $C_{4},$ for $Q=\dprod\limits_{j=1}^{d}[a_{j},b_{j})$
with $0\leq a_{i}<b_{i}<2\pi ,$ 
\begin{equation*}
\mu _{y}\left( Q\right) \leq C_{4}\lambda ^{d}\left( Q\right) \left\Vert
y\right\Vert ^{2}.
\end{equation*}%
By covering any Borel set $\Omega \subseteq \mathbb{T}^{d}$ with countably
many sets in the form of $Q$ whose total measure is arbitrarily close to $%
\lambda ^{d}\left( \Omega \right) ,$ we obtain 
\begin{equation*}
\mu _{y}\left( \Omega \right) \leq C_{4}\lambda ^{d}\left( \Omega \right)
\left\Vert y\right\Vert ^{2}.
\end{equation*}%
Now we define simple vector-valued functions 
\begin{equation*}
g(\boldsymbol{t})=\sum_{k=1}^{N}\chi _{\Omega _{k}}\left( \boldsymbol{t}%
\right) y_{k}
\end{equation*}%
in $L^{2}\left( \mathbb{T}^{d},\mathcal{H}\right) $ with pairwise disjoint
Borel sets $\Omega _{k}$ where $\chi _{\Omega _{k}}$ is the characteristic
function of $\Omega _{k},$ and define $V$ such that 
\begin{equation*}
Vg:=\sum_{k}^{N}E_{\boldsymbol{U}}\left( \Omega _{k}\right) B^{\ast }y_{k}.
\end{equation*}%
Passing to a common refinement shows that this definition is independent of
the representation of $g$. Since the ranges of $E_{\boldsymbol{U}}\left(
\Omega _{k}\right) $ are pairwise orthogonal, 
\begin{eqnarray*}
\left\Vert Vg\right\Vert ^{2} &=&\sum_{k=1}^{N}\left\Vert E_{\boldsymbol{U}%
}\left( \Omega _{k}\right) B^{\ast }y_{k}\right\Vert ^{2} \\
&\leq &C_{4}\sum_{k=1}^{N}\lambda ^{d}\left( \Omega _{k}\right) \left\Vert
y_{k}\right\Vert ^{2}=C_{4}\left\Vert g\right\Vert _{L^{2}\left( \mathbb{T}%
^{d},\mathcal{H}\right) }^{2}.
\end{eqnarray*}%
Thus $V$ extends uniquely to a bounded operator 
\begin{equation*}
V:L^{2}\left( \mathbb{T}^{d},\mathcal{H}\right) \rightarrow \mathcal{H}\text{
with }\left\Vert V\right\Vert \leq \sqrt{C_{4}}.
\end{equation*}%
Taking a sequence $\left\{ y_{\mathbf{n}}\right\} _{\mathbf{n}\in \mathbb{Z}%
^{d}}$ in $\mathcal{H}$ with only finitely many non-zero terms, we set 
\begin{equation*}
f(\boldsymbol{t}):=\sum_{\mathbf{n}}e^{-i\mathbf{n\cdot }\boldsymbol{t}}y_{%
\mathbf{n}}.
\end{equation*}%
By the joint spectral theorem, 
\begin{equation*}
V\left( e^{-i\mathbf{n\cdot }\boldsymbol{t}}y_{\mathbf{n}}\right) =\int_{%
\mathbb{T}^{d}}e^{-i\mathbf{n\cdot }\boldsymbol{t}}dE_{\boldsymbol{U}}\left( 
\boldsymbol{t}\right) B^{\ast }y_{\mathbf{n}}=\mathbf{U}^{-\mathbf{n}%
}B^{\ast }y_{\mathbf{n}}.
\end{equation*}%
Hence $Vf=\sum_{n}\mathbf{U}^{-\mathbf{n}}B^{\ast }y_{n}.$ Therefore 
\begin{eqnarray*}
\left\Vert Vf\right\Vert ^{2} &=&\left\Vert \sum_{n}\mathbf{U}^{-\mathbf{n}%
}B^{\ast }y_{n}\right\Vert ^{2}\leq C_{4}\left\Vert \sum_{n}e^{-i\mathbf{%
n\cdot }\boldsymbol{t}}y_{\mathbf{n}}\right\Vert _{L^{2}\left( \mathbb{T}%
^{d},\mathcal{H}\right) }^{2} \\
&=&C_{4}\int_{\mathbb{T}^{d}}\left\Vert \sum_{n}e^{-i\mathbf{n\cdot }%
\boldsymbol{t}}y_{\mathbf{n}}\right\Vert _{\mathcal{H}}^{2}dt.
\end{eqnarray*}%
Note 
\begin{eqnarray*}
\left\Vert f\left( \boldsymbol{t}\right) \right\Vert _{\mathcal{H}}^{2}
&=&\left\Vert \sum_{n}e^{-i\mathbf{n\cdot }\boldsymbol{t}}y_{\mathbf{n}%
}\right\Vert _{\mathcal{H}}^{2}=\left\langle \sum_{\mathbf{n}}e^{-i\mathbf{%
n\cdot }\boldsymbol{t}}y_{\mathbf{n}},\sum_{n}e^{-i\mathbf{n\cdot }%
\boldsymbol{t}}y_{\mathbf{n}}\right\rangle \\
&=&\sum_{\mathbf{n,m}}e^{-i\left( \mathbf{n-m}\right) \mathbf{\cdot }%
\boldsymbol{t}}\left\langle y_{\mathbf{n}},y_{\mathbf{m}}\right\rangle \text{
for a fixed }\boldsymbol{t.}
\end{eqnarray*}%
It follows that 
\begin{eqnarray*}
\left\Vert f\right\Vert _{L^{2}\left( \mathbb{T}^{d},\mathcal{H}\right)
}^{2} &=&\int_{\mathbb{T}^{d}}\left\Vert \sum_{n}e^{-i\mathbf{n\cdot }%
\boldsymbol{t}}y_{\mathbf{n}}\right\Vert ^{2} \\
&=&\sum_{\mathbf{n,m}}\left\langle y_{\mathbf{n}},y_{\mathbf{m}%
}\right\rangle \int_{\mathbb{T}^{d}}e^{-i\left( \mathbf{n-m}\right) \mathbf{%
\cdot }\boldsymbol{t}}=\left( 2\pi \right) ^{d}\sum_{\mathbf{n}}\left\Vert
y_{\mathbf{n}}\right\Vert ^{2}.
\end{eqnarray*}%
Thus 
\begin{equation}
\left\Vert Vf\right\Vert ^{2}=\left\Vert \sum_{\mathbf{n}}U^{-\mathbf{n}%
}B^{\ast }y_{\mathbf{n}}\right\Vert ^{2}\leq \left( 2\pi \right)
^{d}C_{4}\sum_{\mathbf{n}}\left\Vert y_{\mathbf{n}}\right\Vert ^{2}.
\label{Thm1.1(1)}
\end{equation}%
Let 
\begin{equation*}
c_{00}\left( \mathbb{Z}^{d},\mathcal{H}\right) =\left\{ \left( y_{%
\boldsymbol{n}}\right) :y_{\boldsymbol{n}}=0\text{ except for finitely many }%
\boldsymbol{n}\in \mathbb{Z}^{d}\right\} .
\end{equation*}%
For $y=\left( y_{\boldsymbol{n}}\right) \in c_{00}\left( \mathbb{Z}^{d},%
\mathcal{H}\right) ,$ define 
\begin{equation*}
Sy:=\sum_{\mathbf{n}}U^{-\mathbf{n}}B^{\ast }y_{\mathbf{n}}.
\end{equation*}%
Since 
\begin{eqnarray*}
\left\Vert Sy\right\Vert ^{2} &=&\left\Vert \sum_{\mathbf{n}}U^{-\mathbf{n}%
}B^{\ast }y_{\mathbf{n}}\right\Vert ^{2} \\
&\leq &\left( 2\pi \right) ^{d}C_{4}\sum_{\mathbf{n}}\left\Vert y_{\mathbf{n}%
}\right\Vert ^{2}=\left( 2\pi \right) ^{d}C_{4}\left\Vert y\right\Vert
_{l^{2}\left( \mathbb{Z}^{d},\mathcal{H}\right) }^{2}
\end{eqnarray*}%
by (\ref{Thm1.1(1)}). Hence $S$ can be extended to $l^{2}\left( \mathbb{Z}%
^{d},\mathcal{H}\right) .$ Therefore 
\begin{equation*}
\left\Vert S\right\Vert ^{2}\leq \left( 2\pi \right) ^{d}C_{4}.
\end{equation*}%
Now 
\begin{equation*}
S^{\ast }:\mathcal{H}\rightarrow l^{2}\left( \mathbb{Z}^{d},\mathcal{H}%
\right)
\end{equation*}%
is bounded and we can verify that $S^{\ast }\left( x\right) =\left( BU^{%
\mathbf{n}}x\right) _{\mathbf{n}\in \mathbb{Z}^{d}}\in l^{2}\left( \mathbb{Z}%
^{d},\mathcal{H}\right) .$ Hence%
\begin{equation*}
\sum_{n\in \mathbb{Z}^{d}}\left\Vert B\boldsymbol{U}^{n}x\right\Vert
^{2}\leq \left\Vert S^{\ast }\right\Vert ^{2}\left\Vert x\right\Vert
^{2}=\left\Vert S\right\Vert ^{2}\left\Vert x\right\Vert ^{2}\leq \left(
2\pi \right) ^{d}C_{4}\left\Vert x\right\Vert ^{2}.
\end{equation*}%
It follows that 
\begin{equation*}
C_{1}\leq C_{4}.
\end{equation*}%
Therefore $C_{1}$ is finite as well.

Next, we assume $C_{1}<\infty $. Then, for every $\psi \in \mathcal{H}$%
\begin{equation*}
\sum_{\boldsymbol{n}\in \mathbb{Z}^{d}}\left\Vert B\boldsymbol{U}^{%
\boldsymbol{n}}\psi \right\Vert ^{2}\leq (2\pi )^{d}C_{1}\left\Vert \psi
\right\Vert ^{2},\text{ }\psi \in \mathcal{H}.
\end{equation*}%
Let $f(\boldsymbol{\theta )}=\sum_{\boldsymbol{n}\in A}a_{\boldsymbol{n}%
}e^{i\langle \boldsymbol{n,\theta }\rangle }$, where $A\subset \mathbb{Z}%
^{d} $ is finite. By the joint functional calculus, the Cauchy-Schwarz
inequality, and Parseval's identity 
\begin{eqnarray*}
\left\Vert Bf(\boldsymbol{U})\psi \right\Vert ^{2} &=&\left\Vert \sum_{%
\boldsymbol{n}\in A}a_{\boldsymbol{n}}B\boldsymbol{U}^{\boldsymbol{n}}\psi
\right\Vert ^{2}\leq \left( \sum_{\boldsymbol{n}\in A}\left\vert a_{%
\boldsymbol{n}}\right\vert ^{2}\right) \left( \sum_{\boldsymbol{n}\in
A}\left\Vert B\boldsymbol{U}^{\boldsymbol{n}}\psi \right\Vert ^{2}\right) \\
&\leq &\frac{\left\Vert f\right\Vert _{\mathcal{L}^{2}\left( \mathbb{T}%
^{d}\right) }^{2}}{\left( 2\pi \right) ^{d}}\left( 2\pi \right)
^{d}C_{1}\left\Vert \psi \right\Vert ^{2}=C_{1}\left\Vert f\right\Vert _{%
\mathcal{L}^{2}\left( \mathbb{T}^{d}\right) }^{2}\left\Vert \psi \right\Vert
^{2}.
\end{eqnarray*}%
Consequently, 
\begin{equation*}
\left\Vert Bf(\boldsymbol{U})\right\Vert ^{2}\leq C_{1}\int_{\mathbb{T}%
^{d}}\left\vert f(\boldsymbol{\theta })\right\vert ^{2}d\lambda ^{d}\left( 
\boldsymbol{\theta }\right) .
\end{equation*}

Since trigonometric polynomials are uniformly dense in $C\left( \mathbb{T}%
^{d}\right) ,$ the same inequality holds for every $f\in C\left( \mathbb{T}%
^{d}\right) .$ For $y\in \mathcal{H},$ we define the finite positive Borel
measure 
\begin{equation*}
\mu _{y}\left( \Omega \right) :=\left\langle E_{\boldsymbol{U}}\left( \Omega
\right) B^{\ast }y,B^{\ast }y\right\rangle \text{ for any Borel set }\Omega
\subseteq \mathbb{T}^{d}.
\end{equation*}%
For every nonnegative $h\in C\left( \mathbb{T}^{d}\right) ,$ applying the
preceding inequality to $f=\sqrt{h}$ gives 
\begin{eqnarray*}
\int_{\mathbb{T}^{d}}hd\mu _{y} &=&\left\Vert h^{1/2}\left( \boldsymbol{U}%
\right) B^{\ast }y\right\Vert ^{2}\leq \left\Vert h^{1/2}\left( \boldsymbol{U%
}\right) B^{\ast }\right\Vert ^{2}\left\Vert y\right\Vert ^{2}=\left\Vert
Bh^{1/2}\left( \boldsymbol{U}\right) \right\Vert ^{2}\left\Vert y\right\Vert
^{2} \\
&\leq &C_{1}\left\Vert y\right\Vert ^{2}\int_{\mathbb{T}^{d}}hd\lambda ^{d}.
\end{eqnarray*}%
By regularity of these measures, it follows that 
\begin{equation*}
\left\Vert E_{\boldsymbol{U}}\left( Q\right) B^{\ast }y\right\Vert ^{2}=\mu
_{y}\left( Q\right) \leq C_{1}\left\Vert y\right\Vert ^{2}\lambda ^{d}\left(
Q\right)
\end{equation*}%
for $Q=\dprod\limits_{j=1}^{d}[a_{j},b_{j})\subseteq \mathbb{T}^{d}.$ Then
we obtain that $C_{4}\leq C_{1}.$
\end{proof}

\begin{lemma}
\label{lem1.3}$C_{2}=C_{5}.$
\end{lemma}

\begin{proof}
Since $U_{j}U_{i}=U_{i}U_{j}$ for any $1\leq i,j\leq d,$ we have 
\begin{eqnarray*}
\Delta \left( \boldsymbol{U},\boldsymbol{z}\right)
&=&\dprod\limits_{j=1}^{d}\delta \left( U_{j},z_{j}\right) =\frac{1}{\left(
2\pi \right) ^{d}}\dprod\limits_{j=1}^{d}\left( 1-\left\vert
z_{j}\right\vert ^{2}\right) \left\vert 1-z_{j}U_{j}^{\ast }\right\vert ^{-2}
\\
&=&\frac{1}{\left( 2\pi \right) ^{d}}\dprod\limits_{j=1}^{d}\left(
1-\left\vert z_{j}\right\vert ^{2}\right) \left\vert R_{U_{j}}\left(
z_{j}\right) \right\vert ^{2} \\
&=&\frac{1}{\left( 2\pi \right) ^{d}}\left( \dprod\limits_{j=1}^{d}\left(
1-\left\vert z_{j}\right\vert ^{2}\right) \right) \cdot \left(
\dprod\limits_{j=1}^{d}\left\vert R_{U_{j}}\left( z_{j}\right) \right\vert
^{2}\right) .
\end{eqnarray*}%
Let $X=\dprod\limits_{j}\left\vert R_{U_{j}}\left( z_{j}\right) \right\vert
. $ Then $\left\Vert XB^{\ast }\psi \right\Vert ^{2}=\left\langle XB^{\ast
}\psi ,XB^{\ast }\psi \right\rangle .$ Hence 
\begin{eqnarray*}
\left\langle \Delta \left( \mathbf{U},\boldsymbol{z}\right) B^{\ast }\psi
,B^{\ast }\psi \right\rangle &=&\frac{1}{\left( 2\pi \right) ^{d}}%
\dprod\limits_{j}\left( 1-\left\vert z_{j}\right\vert ^{2}\right)
\left\langle X^{2}B^{\ast }\psi ,B^{\ast }\psi \right\rangle \\
&=&\frac{1}{\left( 2\pi \right) ^{d}}\dprod\limits_{j}\left( 1-\left\vert
z_{j}\right\vert ^{2}\right) \left\Vert XB^{\ast }\psi \right\Vert ^{2}.
\end{eqnarray*}

It follows that 
\begin{equation*}
C_{2}=\frac{1}{\left( 2\pi \right) ^{d}}\sup_{\left\Vert \psi \right\Vert =1,%
\boldsymbol{z\in }\mathbb{D}^{d}}\dprod\limits_{j}\left( 1-\left\vert
z_{j}\right\vert ^{2}\right) \left\Vert XB^{\ast }\psi \right\Vert
^{2}=C_{5}.
\end{equation*}
\end{proof}

\begin{remark}
\label{remark2} If $B$ is jointly $\boldsymbol{U}$-smooth, then for $%
Q=\dprod\limits_{j=1}^{d}[a_{j},b_{j})\subseteq \mathbb{T}^{d}$ and $y$ with 
$\left\Vert y\right\Vert =1$ 
\begin{eqnarray*}
\left\vert \left\langle E_{\boldsymbol{U}}\left( Q\right) B^{\ast }y,B^{\ast
}y\right\rangle \right\vert &\leq &\sup_{\left\Vert y\right\Vert
=1}\left\vert \left\langle E_{\boldsymbol{U}}\left( Q\right) B^{\ast
}y,B^{\ast }y\right\rangle \right\vert =\sup_{\left\Vert y\right\Vert
=1}\left\Vert E_{\boldsymbol{U}}\left( Q\right) B^{\ast }y\right\Vert ^{2}%
\text{ } \\
&\leq &\lambda ^{d}\left( Q\right) C_{4}.
\end{eqnarray*}%
It implies that the measure $\mu _{B^{\ast }y}$ is absolutely continuous
with respect to the Lebesgue measure $\lambda ^{d}$ on $\mathbb{T}^{d}$
where $\mu _{B^{\ast }y}$ is defined by $\mu _{B^{\ast }y}\left( Q\right)
=\left\langle E_{\boldsymbol{U}}\left( Q\right) B^{\ast }y,B^{\ast
}y\right\rangle .$ Hence $\overline{R\left( B^{\ast }\right) }\subset
P_{ac}\left( \boldsymbol{U}\right) \left( \mathcal{H}\right) $ if $B$ is
jointly $\boldsymbol{U}$-smooth in $\mathcal{M}$.
\end{remark}

\begin{lemma}
\label{lem1.4}If at least one of $C_{2}$ and $C_{4}$ is finite, then $%
C_{2}=C_{4}.$
\end{lemma}

\begin{proof}
First, we assume $C_{4}<\infty .$ Then by Remark \ref{remark2}, $\frac{d\mu
_{B^{\ast }y}}{d\lambda ^{d}}$ exists for $y$ and 
\begin{equation}
\frac{d\mu _{B^{\ast }y}}{d\lambda ^{d}}\leq C_{4}\left\Vert y\right\Vert
^{2}\text{ }  \label{444}
\end{equation}%
$\lambda ^{d}$-almost everywhere. It follows that%
\begin{eqnarray*}
&&\left\langle \Delta \left( \boldsymbol{U},\boldsymbol{re}^{i\boldsymbol{%
\theta }}\right) B^{\ast }\psi ,B^{\ast }\psi \right\rangle =\frac{1}{\left(
2\pi \right) ^{d}}\dint\limits_{\mathbb{T}^{d}}P_{\boldsymbol{r}}\left( 
\boldsymbol{\theta -\lambda }\right) d\mu _{B^{\ast }\psi }\left( 
\boldsymbol{\lambda }\right) \\
&=&\frac{1}{\left( 2\pi \right) ^{d}}\dint\limits_{\mathbb{T}^{d}}P_{%
\boldsymbol{r}}\left( \boldsymbol{\theta -\lambda }\right) \frac{d\mu
_{B^{\ast }\psi }\left( \boldsymbol{\lambda }\right) }{d\lambda ^{d}}%
d\lambda ^{d} \\
&\leq &\frac{1}{\left( 2\pi \right) ^{d}}C_{4}\left\Vert \psi \right\Vert
^{2}\dint\limits_{\mathbb{T}^{d}}P_{\boldsymbol{r}}\left( \boldsymbol{\theta
-\lambda }\right) d\boldsymbol{\lambda }=C_{4}\left\Vert \psi \right\Vert
^{2}\text{ }
\end{eqnarray*}%
by (\ref{p17}), (\ref{p1}) and (\ref{444}). It implies that $C_{2}\leq
C_{4}. $

Now we are going to show that $C_{4}\leq C_{2}$ as $C_{2}<\infty .$ Since 
\begin{equation*}
\left\langle \Delta \left( \boldsymbol{U},\boldsymbol{z}\right) B^{\ast
}y,B^{\ast }y\right\rangle \leq \left\Vert B\Delta \left( \boldsymbol{U},%
\boldsymbol{z}\right) B^{\ast }\right\Vert \left\Vert y\right\Vert ^{2}\leq
C_{2}\left\Vert y\right\Vert ^{2}
\end{equation*}%
and by (\ref{p17}) 
\begin{eqnarray*}
\left\langle \Delta \left( \boldsymbol{U},\boldsymbol{r},\boldsymbol{\theta }%
\right) B^{\ast }y,B^{\ast }y\right\rangle &=&\frac{1}{\left( 2\pi \right)
^{d}}\dint\limits_{\mathbb{T}^{d}}P_{\boldsymbol{r}}\left( \boldsymbol{%
\theta -\lambda }\right) d\mu _{B^{\ast }y}\left( \boldsymbol{\lambda }%
\right) \\
&=&P_{\boldsymbol{r}}\ast \mu _{B^{\ast }y}\left( \boldsymbol{\theta }%
\right) ,
\end{eqnarray*}%
we obtain that for any continuous function $\varphi \left( \boldsymbol{%
\theta }\right) \in \mathbb{C}\left( \mathbb{T}^{d}\right) $ with $\varphi
\geq 0,$ 
\begin{equation}
\left( \varphi \cdot \left( P_{\boldsymbol{r}}\ast \mu _{B^{\ast }y}\right)
\right) \left( \boldsymbol{\theta }\right) \leq C_{2}\left\Vert y\right\Vert
^{2}\varphi \left( \boldsymbol{\theta }\right) .  \label{443}
\end{equation}%
By (\ref{443}) and Fubini's theorem, we get 
\begin{equation}
\int_{\mathbb{T}^{d}}\left\vert \varphi \cdot \left( P_{\boldsymbol{r}}\ast
\mu _{B^{\ast }y}\right) \left( \boldsymbol{\theta }\right) \right\vert d%
\boldsymbol{\theta }\leq C_{2}\left\Vert y\right\Vert ^{2}\int_{\mathbb{T}%
^{d}}\varphi \left( \boldsymbol{\theta }\right) d\boldsymbol{\theta }\mathbf{%
.}  \label{lemma3(5)}
\end{equation}%
and 
\begin{eqnarray}
&&\int_{\mathbb{T}^{d}}\left( \varphi \cdot \left( P_{\boldsymbol{r}}\ast
\mu _{B^{\ast }y}\right) \right) \left( \boldsymbol{\theta }\right) d%
\boldsymbol{\theta }\mathbf{=}\int_{\mathbb{T}^{d}}\varphi \left( 
\boldsymbol{\theta }\right) \left( \int_{\mathbb{T}^{d}}\frac{1}{\left( 2\pi
\right) ^{d}}P_{\boldsymbol{r}}\left( \boldsymbol{\theta }\mathbf{-}%
\boldsymbol{\lambda }\right) d\mu _{B^{\ast }y}\left( \boldsymbol{\lambda }%
\right) \right) d\boldsymbol{\theta }  \notag \\
&=&\frac{1}{\left( 2\pi \right) ^{d}}\int_{\mathbb{T}^{d}}\left( \int_{%
\mathbb{T}^{d}}P_{\boldsymbol{r}}\left( \boldsymbol{\theta }\mathbf{-}%
\boldsymbol{\lambda }\right) \varphi \left( \boldsymbol{\theta }\right) d%
\boldsymbol{\theta }\right) d\mu _{B^{\ast }y}\left( \boldsymbol{\lambda }%
\right) =\int_{\mathbb{T}^{d}}\left( P_{\boldsymbol{r}}\ast \varphi \right)
\left( \boldsymbol{\lambda }\right) d\mu _{B^{\ast }y}\left( \boldsymbol{%
\lambda }\right) .  \label{lemma 3(3)}
\end{eqnarray}%
Since 
\begin{equation*}
\lim_{\underset{\boldsymbol{r}\in \left( 0,1\right) ^{d}}{\boldsymbol{r}%
\rightarrow \boldsymbol{1}}}\left\Vert \left( P_{\boldsymbol{r}}\ast \varphi
\right) -\varphi \right\Vert _{\infty }=0
\end{equation*}%
by Theorem 3.1 in \cite{O}$,$ we have 
\begin{equation*}
\left\vert \int_{\mathbb{T}^{d}}P_{\boldsymbol{r}}\ast \varphi d\mu
_{B^{\ast }y}\left( \boldsymbol{\lambda }\right) -\int_{\mathbb{T}%
^{d}}\varphi \left( \boldsymbol{\lambda }\right) d\mu _{B^{\ast }y}\left( 
\boldsymbol{\lambda }\right) \right\vert \leq \left\vert \mu _{B^{\ast
}y}\left( \mathbb{T}^{d}\right) \right\vert \left\Vert P_{\boldsymbol{r}%
}\ast \varphi -\varphi \right\Vert _{\infty }\rightarrow 0.
\end{equation*}%
It implies that 
\begin{equation}
\lim_{\underset{\boldsymbol{r}\in \left( 0,1\right) ^{d}}{\boldsymbol{r}%
\rightarrow \boldsymbol{1}}}\int_{\mathbb{T}^{d}}P_{\boldsymbol{r}}\ast
\varphi \left( \boldsymbol{\lambda }\right) d\mu _{B^{\ast }y}\left( 
\boldsymbol{\lambda }\right) \boldsymbol{=}\int_{\mathbb{T}^{d}}\varphi
\left( \boldsymbol{\lambda }\right) d\mu _{B^{\ast }y}\left( \boldsymbol{%
\lambda }\right) .  \label{442}
\end{equation}%
Therefore we obtain that 
\begin{equation*}
\int_{\mathbb{T}^{d}}\varphi \left( \boldsymbol{\lambda }\right) d\mu
_{B^{\ast }y}\left( \boldsymbol{\lambda }\right) \leq C_{2}\left\Vert
y\right\Vert ^{2}\int_{\mathbb{T}^{d}}\varphi \left( \boldsymbol{\theta }%
\right) d\boldsymbol{\theta }\mathbf{.}
\end{equation*}%
by (\ref{lemma 3(3)}) and (\ref{442}). Since $\varphi \geq 0$ in $\mathbb{C}%
\left( \mathbb{T}^{d}\right) ,$ by a standard measure-theoretic argument, we
obtain 
\begin{equation*}
\left\Vert E\left( Q\right) B^{\ast }y\right\Vert ^{2}\leq C_{2}\lambda
^{d}\left( Q\right) \left\Vert y\right\Vert ^{2}
\end{equation*}%
for $Q=\dprod\limits_{j=1}^{d}[a_{j},b_{j})\subseteq \mathbb{T}^{d}.$ Hence $%
C_{4}\leq C_{2}.$
\end{proof}

\bigskip

\textbf{Proof of Theorem \ref{thm1.1}}\bigskip : The conclusion follows from
Lemmas \ref{lem1.1}-Lemma \ref{lem1.4}.

Now we are ready to define the projection $P_{ac}^{\infty }\left( 
\boldsymbol{U},\mathcal{M}\right) $ which generalizes the absolutely
continuous projection $P_{ac}\left( \boldsymbol{U}\right) $ to the setting
of von Neumann algebras$\mathcal{.}$

\begin{definition}
\label{li4}Let $\boldsymbol{U}\mathbf{=}\left( U_{1},\cdots ,U_{d}\right) $
be a $d$-tuple of commuting unitary elements in $\mathcal{M}$. Recall that $%
\mathcal{S}_{ac}^{\infty }\left( \boldsymbol{U},\mathcal{M}\right) $ is the
collection of all jointly $\boldsymbol{U}$-smooth operators in $\mathcal{M}$%
. We set $\mathcal{P}_{ac}^{\infty }\left( U,\mathcal{M}\right) $ to be the
collection of all projections in $\mathcal{S}_{ac}^{\infty }\left( U,%
\mathcal{M}\right) $. The norm absolutely continuous support $P_{ac}^{\infty
}\left( \boldsymbol{U},\mathcal{M}\right) $ is defined below 
\begin{equation*}
P_{ac}^{\infty }\left( \boldsymbol{U},\mathcal{M}\right) \triangleq \bigvee
_{P\in \mathcal{P}_{ac}^{\infty }\left( \boldsymbol{U},\mathcal{M}\right) }P.
\end{equation*}
\end{definition}

The following proposition characterizes $P_{ac}^{\infty }\left( \boldsymbol{U%
},\mathcal{M}\right) .$

\begin{proposition}
\label{pro11}Let $\boldsymbol{U}\mathbf{=}\left( U_{1},\cdots ,U_{d}\right) $
be a $d$-tuple of commuting unitary elements in $\mathcal{M}$. Then 
\begin{equation*}
P_{ac}^{\infty }\left( \boldsymbol{U},\mathcal{M}\right) =\bigvee _{P\in 
\mathcal{P}_{ac}^{\infty }\left( \boldsymbol{U},\mathcal{M}\right) }E_{P}\in 
\mathcal{Z}\left( \mathcal{M\cap A}^{\prime }\right)
\end{equation*}%
where $E_{P}$ is defined in (\ref{pr5}) and $\mathcal{A}$ is the von Neumann
algebra generated by $U_{1},\cdots ,U_{d}$.
\end{proposition}

\begin{proof}
For any projection $P\in \mathcal{P}_{ac}^{\infty }\left( \boldsymbol{U},%
\mathcal{M}\right) $ and $j\in \left\{ 1,\cdots ,d\right\} ,$ we set $%
B=U_{j}PU_{j}^{\ast }.$ Then 
\begin{eqnarray*}
\dsum\limits_{\mathbf{n\in }\mathbb{Z}^{d}}\left\Vert B\boldsymbol{U}^{%
\mathbf{n}}x\right\Vert ^{2} &=&\dsum\limits_{\mathbf{n\in }\mathbb{Z}%
^{d}}\left\Vert PU_{j}^{\ast }\boldsymbol{U}^{\mathbf{n}}x\right\Vert
^{2}=\dsum\limits_{\mathbf{n\in }\mathbb{Z}^{d}}\left\Vert P\mathbf{U}^{%
\mathbf{n}-e_{j}}x\right\Vert \\
&=&\dsum\limits_{\mathbf{m\in }\mathbb{Z}^{d}}\left\Vert P\boldsymbol{U}^{%
\mathbf{m}}x\right\Vert ^{2}\leq \left( 2\pi \right) ^{d}C_{P}\left\Vert
x\right\Vert ^{2}
\end{eqnarray*}%
where $e_{j}$ is the $j$-th standard basis vector of $\mathbb{Z}^{d}$ and $%
C_{P}$ is the jointly $\boldsymbol{U}$-smoothness constant of $P$. Thus $%
U_{j}PU_{j}^{\ast }$ is also a jointly $\boldsymbol{U}$-smooth projection.
Therefore $U_{j}P_{ac}^{\infty }\left( \boldsymbol{U},\mathcal{M}\right)
U_{j}^{\ast }=P_{ac}^{\infty }\left( \boldsymbol{U},\mathcal{M}\right) $ for 
$j=1,\cdots ,d.$ It implies that 
\begin{equation}
P_{ac}^{\infty }\left( \boldsymbol{U},\mathcal{M}\right) \in \mathcal{A}%
^{\prime }.  \label{Proposition 3(1)}
\end{equation}%
Now for any unitary $V\in \mathcal{M}\cap \mathcal{A}^{\prime },$ 
\begin{eqnarray*}
\frac{\left\Vert VPV^{\ast }E_{\boldsymbol{U}}\left( Q\right) \psi
\right\Vert ^{2}}{\lambda ^{d}\left( Q\right) } &=&\frac{\left\Vert PE_{%
\boldsymbol{U}}\left( Q\right) V^{\ast }\psi \right\Vert ^{2}}{\lambda
^{d}\left( Q\right) } \\
&\leq &C_{P}\left\Vert V^{\ast }\psi \right\Vert ^{2}=C_{P}\left\Vert \psi
\right\Vert ^{2}
\end{eqnarray*}%
for every $Q=\dprod\limits_{j=1}^{d}[a_{j},b_{j})\subseteq \mathbb{T}^{d}$.
Then $VPV^{\ast }$ is jointly $\boldsymbol{U}$-smooth if $P$ is jointly $%
\boldsymbol{U}$-smooth. It implies that 
\begin{equation}
VP_{ac}^{\infty }\left( \boldsymbol{U},\mathcal{M}\right) V^{\ast
}=P_{ac}^{\infty }\left( \boldsymbol{U},\mathcal{M}\right) .
\label{Proposition 3(2)}
\end{equation}%
Combining \ref{Proposition 3(1)} and \ref{Proposition 3(2)}, we conclude
that 
\begin{equation*}
P_{ac}^{\infty }\left( \boldsymbol{U},\mathcal{M}\right) \in \mathcal{Z}%
\left( \mathcal{M}\cap \mathcal{A}^{\prime }\right) .
\end{equation*}

Finally, for each $P\in \mathcal{P}_{ac}^{\infty }\left( \boldsymbol{U},%
\mathcal{M}\right) ,$ the projection $E_{P}$ is the smallest projection in $%
\mathcal{M}\cap \mathcal{A}^{\prime }$ dominating $P\in \mathcal{P}%
_{ac}^{\infty }\left( \boldsymbol{U},\mathcal{M}\right) .$ Since $P\in 
\mathcal{P}_{ac}^{\infty }\left( \boldsymbol{U},\mathcal{M}\right) $ and $%
P_{ac}^{\infty }\left( \boldsymbol{U},\mathcal{M}\right) \in \mathcal{M}\cap 
\mathcal{A}^{\prime },$ we have that $E_{P}\leq P_{ac}^{\infty }\left( 
\boldsymbol{U},\mathcal{M}\right) .$ It follows that $\bigvee _{P\in 
\mathcal{P}_{ac}^{\infty }\left( \boldsymbol{U},\mathcal{M}\right)
}E_{P}\leq P_{ac}^{\infty }\left( \boldsymbol{U},\mathcal{M}\right) .$
Conversely, $P\leq E_{P}$ for every such $P$, so we conclude that 
\begin{equation*}
P_{ac}^{\infty }\left( \boldsymbol{U},\mathcal{M}\right) =\bigvee _{P\in 
\mathcal{P}_{ac}^{\infty }\left( \boldsymbol{U},\mathcal{M}\right) }E_{P}\in 
\mathcal{Z}\left( \mathcal{M\cap A}^{\prime }\right) .
\end{equation*}
\end{proof}

\begin{corollary}
\label{cor2}Let $\boldsymbol{U}\mathbf{=}\left( U_{1},\cdots ,U_{d}\right) $
be a $d$-tuple of commuting unitary elements in $\mathcal{M}$. Suppose $A\in 
\mathcal{S}_{ac}^{\infty }\left( U,\mathcal{M}\right) .$ Then 
\begin{equation*}
\overline{Ran\left( A^{\ast }\right) }\subseteq P_{ac}^{\infty }\left( 
\boldsymbol{U},\mathcal{M}\right) \left( \mathcal{H}\right)
\end{equation*}%
and 
\begin{equation*}
P_{ac}^{\infty }\left( \boldsymbol{U},\mathcal{M}\right) \leq P_{ac}\left( 
\boldsymbol{U}\right) .
\end{equation*}
\end{corollary}

\begin{proof}
Since $A$ is jointly $\boldsymbol{U}$-smooth, there exists a constant $%
C_{A}<\infty $ such that 
\begin{equation*}
\sum_{\boldsymbol{n}\in \mathbb{Z}^{d}}\left\Vert A\boldsymbol{U}^{%
\boldsymbol{n}}x\right\Vert \leq \left( 2\pi \right) ^{d}C_{A}\left\Vert
x\right\Vert ^{2},\text{ }x\in \mathcal{H}.
\end{equation*}%
For $\varepsilon >0,$ we define the spectral projection 
\begin{equation*}
P_{\varepsilon }:=E_{A^{\ast }A}\left( [\varepsilon ,\infty \right) )\in 
\mathcal{M.}
\end{equation*}%
By functional calculus, $P_{\varepsilon }\leq \frac{A^{\ast }A}{\varepsilon }
$ and then 
\begin{equation*}
\left\Vert P_{\varepsilon }y\right\Vert ^{2}=\left\langle P_{\varepsilon
}y,y\right\rangle \leq \varepsilon ^{-1}\left\Vert Ay\right\Vert ^{2}\text{
for }y\in \mathcal{H}.
\end{equation*}%
Consequently%
\begin{equation*}
\sum_{\boldsymbol{n}\in \mathbb{Z}^{d}}\left\Vert P_{\varepsilon }%
\boldsymbol{U}^{\boldsymbol{n}}x\right\Vert ^{2}\leq \varepsilon \sum_{%
\boldsymbol{n}\in \mathbb{Z}^{d}}\left\Vert A\boldsymbol{U}^{\boldsymbol{n}%
}x\right\Vert ^{2}\leq \varepsilon ^{-1}\left( 2\pi \right)
^{d}C_{A}\left\Vert x\right\Vert ^{2}.
\end{equation*}%
Thus $P_{\varepsilon }$ is a jointly $\boldsymbol{U}$-smooth projection for
every $\varepsilon >0$. By the definition of $P_{ac}^{\infty }\left( 
\boldsymbol{U},\mathcal{M}\right) ,$%
\begin{equation*}
P_{\varepsilon }\leq P_{ac}^{\infty }\left( \boldsymbol{U},\mathcal{M}%
\right) .
\end{equation*}%
Letting $\varepsilon \rightarrow 0,$ we obtain 
\begin{equation*}
l\left( A^{\ast }A\right) =\bigvee_{\varepsilon >0}P_{\varepsilon }\leq
P_{ac}^{\infty }\left( \boldsymbol{U},\mathcal{M}\right) .
\end{equation*}%
Hence 
\begin{equation*}
\overline{Ran\left( A^{\ast }\right) }\subseteq P_{ac}^{\infty }\left( 
\boldsymbol{U},\mathcal{M}\right) \left( \mathcal{H}\right) \text{ for every 
}A\in \mathcal{S}_{ac}^{\infty }\left( U,\mathcal{M}\right) .
\end{equation*}%
Combining with Remark \ref{remark2} and the definition of $P_{ac}^{\infty
}\left( \boldsymbol{U},\mathcal{M}\right) $, we have 
\begin{equation*}
P_{ac}^{\infty }\left( \boldsymbol{U},\mathcal{M}\right) \leq P_{ac}\left( 
\boldsymbol{U}\right) .
\end{equation*}
\end{proof}

Note $\mathcal{B(H)}$ is a type $I$ von Neumann algebra. The following
corollary shows that the norm absolutely continuous support coincides with
the usual absolutely continuous projection when $\mathcal{M}=\mathcal{B(H)}$.

\begin{corollary}
\bigskip Let $\boldsymbol{U}\mathbf{=}\left( U_{1},\cdots ,U_{d}\right) $ be
a $d$-tuple of commuting unitary elements in $\mathcal{B(H)}$. Then 
\begin{equation*}
P_{ac}^{\infty }\left( \boldsymbol{U},\mathcal{B(H)}\right) =P_{ac}\left( 
\boldsymbol{U}\right) .
\end{equation*}
\end{corollary}

\begin{proof}
If $P$ is jointly $\boldsymbol{U}$-smooth, then by Remark \ref{remark2} we
have $\overline{R\left( P\right) }\subset P_{ac}\left( \boldsymbol{U}\right) 
$. Hence 
\begin{equation}
P_{ac}^{\infty }\left( \boldsymbol{U},\mathcal{B(H)}\right) \leq
P_{ac}\left( \boldsymbol{U}\right) .  \label{cor1(1)}
\end{equation}

To prove the reverse inequality, for a fixed $x\in \mathcal{H}_{ac}\left( 
\boldsymbol{U}\right) ,$ we have 
\begin{equation*}
\left\langle E\left( \Omega \right) x,x\right\rangle =\int_{\Omega
}w_{x}\left( \boldsymbol{t}\right) d\lambda ^{d}\left( t\right)
\end{equation*}%
for a Borel set $\Omega $ and a function $0\leq w_{x}\in \mathcal{L}%
^{1}\left( \mathbb{T}^{d}\right) .$ For $N\geq 1,$ we define 
\begin{equation*}
\Omega _{N}:=\left\{ \boldsymbol{z}\in \mathbb{T}^{d}:w_{x}\left( 
\boldsymbol{z}\right) \leq N\right\}
\end{equation*}%
and the truncated vector 
\begin{equation*}
x_{N}:=E_{\boldsymbol{U}}\left( \Omega _{N}\right) x.
\end{equation*}%
Because $E_{\boldsymbol{U}}\left( \Omega _{N}\right) \overset{sot}{%
\rightarrow }I$ as $N\rightarrow \infty ,$%
\begin{equation}
\left\Vert x-x_{N}\right\Vert ^{2}=\left\langle E_{\boldsymbol{U}}\left( 
\mathbb{T}^{d}\backslash \Omega _{N}\right) x,x\right\rangle \rightarrow 0%
\text{ as }N\rightarrow \infty .  \label{corollary1(1)}
\end{equation}%
Since 
\begin{eqnarray*}
&&\underset{\left\Vert \psi \right\Vert =1,Q=[a_{1},b_{1})\times \cdots
\times \lbrack a_{d},b_{d}\subseteq \mathbb{T}^{d}}{\sup }\frac{\left\Vert
\left( x_{N}\otimes x_{N}\right) E_{\boldsymbol{U}}\left( Q\right) \psi
\right\Vert ^{2}}{\lambda ^{d}\left( Q\right) } \\
&=&\underset{\left\Vert \psi \right\Vert =1,Q=[a_{1},b_{1})\times \cdots
\times \lbrack a_{d},b_{d}\subseteq \mathbb{T}^{d}}{\sup }\frac{\left\langle
E_{\boldsymbol{U}}\left( Q\right) \psi ,x_{N}\right\rangle ^{2}\left\Vert
x_{N}\right\Vert ^{2}}{\lambda ^{d}\left( Q\right) } \\
&\leq &\underset{\left\Vert \psi \right\Vert =1,Q=[a_{1},b_{1})\times \cdots
\times \lbrack a_{d},b_{d}\subseteq \mathbb{T}^{d}}{\sup }\frac{\left\Vert
\psi \right\Vert ^{2}\left\Vert E_{\boldsymbol{U}}\left( Q\right)
x_{N}\right\Vert ^{2}\left\Vert x_{N}\right\Vert ^{2}}{\lambda ^{d}\left(
Q\right) } \\
&=&\underset{Q=[a_{1},b_{1})\times \cdots \times \lbrack
a_{d},b_{d}\subseteq \mathbb{T}^{d}}{\sup }\frac{\left\langle E_{\boldsymbol{%
U}}\left[ Q\cap \Omega _{N}\right] x,x\right\rangle \left\Vert
x_{N}\right\Vert ^{2}}{\lambda ^{d}\left( Q\right) } \\
&\leq &\underset{Q=[a_{1},b_{1})\times \cdots \times \lbrack
a_{d},b_{d}\subseteq \mathbb{T}^{d}}{\sup }\frac{\left( \dint\limits_{Q\cap
\Omega _{N}}w_{x}\left( \boldsymbol{t}\right) d\lambda ^{d}\left( t\right)
\right) \boldsymbol{\cdot }\left\Vert x_{N}\right\Vert ^{2}}{\lambda
^{d}\left( Q\right) }\leq N\left\Vert x_{N}\right\Vert ^{2}.
\end{eqnarray*}%
Hence $x_{N}\otimes x_{N}\in \mathcal{S}_{ac}^{\infty }\left( \boldsymbol{U},%
\mathcal{B(H)}\right) $ for the fixed number $N$. It implies that 
\begin{equation*}
x_{N}\in \overline{Ran\left( x_{N}\otimes x_{N}\right) }\subseteq
P_{ac}^{\infty }\left( \boldsymbol{U},\mathcal{B(H)}\right) \left( \mathcal{H%
}\right)
\end{equation*}%
for every $N$ by Corollary \ref{cor2} and then $x\in P_{ac}^{\infty }\left( 
\boldsymbol{U},\mathcal{B(H)}\right) \left( \mathcal{H}\right) $ by (\ref%
{corollary1(1)}). Therefore 
\begin{equation}
P_{ac}^{\infty }\left( \boldsymbol{U},\mathcal{B(H)}\right) \geq
P_{ac}\left( \boldsymbol{U}\right) .  \label{cor1(2)}
\end{equation}%
Combining (\ref{cor1(1)}) and (\ref{cor1(2)}), we obtain that 
\begin{equation*}
P_{ac}^{\infty }\left( \boldsymbol{U},\mathcal{B(H)}\right) =P_{ac}\left( 
\boldsymbol{U}\right) .
\end{equation*}
\end{proof}

\begin{example}
This example shows that $P_{ac}^{\infty }\left( \boldsymbol{U},\mathcal{M}%
\right) $ can be strictly smaller than $P_{ac}\left( \boldsymbol{U}\right) $
for a $d$-tuple of commuting unitaries $\boldsymbol{U}$ in $\mathcal{M.}$

Let $\mathcal{K}=\mathcal{L}^{2}\left( \mathbb{T}^{d},\lambda ^{d}\right) $
and $\mathcal{D=}\left\{ M_{f}:f\in \mathcal{L}^{\infty }\left( \mathbb{T}%
^{d},\lambda \right) \right\} $ and $V_{j}=M_{z_{j}}$ where $M_{f}$ denotes
multiplication by $f.$ Then it is not hard to check that $P_{ac}\left( 
\boldsymbol{V}\right) =I$ where $\boldsymbol{V}=\left( V_{1},\cdots
,V_{d}\right) .$ Since any non-zero projection $P$ in $\mathcal{D}$
satisfies 
\begin{equation*}
\sup_{\left\Vert x\right\Vert =1}\left\Vert P\boldsymbol{V}^{\boldsymbol{n}%
}x\right\Vert =\left\Vert \boldsymbol{V}^{\boldsymbol{n}}Px\right\Vert =1%
\text{ for every }\boldsymbol{n}\in \mathbb{Z}^{d},
\end{equation*}%
we have 
\begin{equation*}
\sup_{\left\Vert x\right\Vert =1}\sum_{\boldsymbol{n}\in \mathbb{Z}%
^{d}}\left\Vert P\boldsymbol{V}^{\boldsymbol{n}}x\right\Vert ^{2}=\infty .
\end{equation*}
Hence $P_{ac}^{\infty }\left( \boldsymbol{V},\mathcal{D}\right)
=0<P_{ac}\left( \boldsymbol{V}\right) .$ Let $\mathcal{H}=\mathcal{K\oplus K}
$, $\mathcal{M}=\mathcal{B}\left( \mathcal{K}\right) \mathcal{\oplus D}$ and 
$\boldsymbol{U=V\oplus V}.$ Then $\mathcal{M}\neq \mathcal{B}\left( \mathcal{%
H}\right) $ and $P_{ac}^{\infty }\left( \boldsymbol{U},\mathcal{M}\right)
=I\oplus 0<P_{ac}\left( \boldsymbol{U}\right) =I\oplus I.$
\end{example}

\section{Wave operators of $\boldsymbol{U}$ and $\boldsymbol{U}_{0}$ in $%
\mathcal{M}$ ($d\geq 2$)}

Let $d\geq 2$ and let 
\begin{equation*}
\boldsymbol{U}=\left( U_{1},\cdots ,U_{d}\right) \text{ and }\boldsymbol{U}%
_{0}=\left( U_{0,1},\cdots ,U_{0,d}\right)
\end{equation*}%
be two $d$-tuples of commuting unitaries in a von Neumann algebra $\mathcal{M%
}$ and let $J\in \mathcal{M.}$ In this section, we develop a stationary
approach to scattering theory for commuting $d$-tuples of unitaries in $%
\mathcal{M}$. Our analysis is based on the factorizations $%
JU_{0,k}-U_{k}J=G_{k}^{\ast }G_{0,k}$, together with joint smoothness
assumptions on the factors. Using products of resolvents, we construct
sector stationary wave operators $w_{\boldsymbol{\sigma }}=w_{\boldsymbol{%
\sigma }}\left( \boldsymbol{U},\boldsymbol{U}_{0},J,\mathcal{M}\right) $ and
relate them to the strong wave operator $W=W\left( \boldsymbol{U},%
\boldsymbol{U}_{0},J,\mathcal{M}\right) $.

Before stating our first two main results in this section, we define the
strong wave operator as follows.

\begin{definition}
For $d\geq 2,$ let 
\begin{equation*}
\boldsymbol{U}=\left( U_{1},\cdots ,U_{d}\right) \text{ and }\boldsymbol{U}%
_{0}=\left( U_{0,1},\cdots ,U_{0,d}\right)
\end{equation*}%
be two $d$-tuples of commuting unitaries in a von Neumann algebra $\mathcal{M%
}$ and let $J\in \mathcal{M}.$ Let 
\begin{equation*}
W\left( \boldsymbol{U},\boldsymbol{U}_{0},J,\mathcal{M}\right) :=\underset{%
\left\Vert \boldsymbol{m}\right\Vert _{\infty }\rightarrow \infty }{sot\text{%
-}\lim }\boldsymbol{U}^{\boldsymbol{m}}J\boldsymbol{U}_{0}^{-\boldsymbol{m}%
}P_{ac}^{\infty }\left( \boldsymbol{U}_{0},\mathcal{M}\right)
\end{equation*}%
where $\boldsymbol{m=}\left( m_{1},\cdots ,m_{d}\right) \in \mathbb{Z}^{d}.$
If the "sot" limit exists, then $W\left( \boldsymbol{U},\boldsymbol{U}_{0},J,%
\mathcal{M}\right) $ is called the strong wave operator of $\boldsymbol{U},%
\boldsymbol{U}_{0}$ and $J$ in $\mathcal{M}.$ Here $\left\Vert \boldsymbol{m}%
\right\Vert _{\infty }=\max_{1\leq j\leq d}\left\vert m_{j}\right\vert .$
\end{definition}

In the rest of this section, we always assume that 
\begin{equation*}
\boldsymbol{U}=\left( U_{1},\cdots ,U_{d}\right) \text{ and }\boldsymbol{U}%
_{0}=\left( U_{01},\cdots ,U_{0d}\right)
\end{equation*}%
are two $d$-tuples of commuting unitaries in a von Neumann algebra $\mathcal{%
M}$, and $J\in \mathcal{M}.$ For $\boldsymbol{\sigma }=\left( \sigma
_{1},\cdots ,\sigma _{d}\right) \in \left\{ +1,-1\right\} ^{d}$ and $%
\boldsymbol{\theta }\in \mathbb{T}^{d},$ we write $R_{j}^{\sigma
_{j}}:=R_{j}^{\sigma _{j}}\left( U_{j},r_{j},\theta _{j}\right) .$ Then we
define the sector product resolvent 
\begin{equation*}
R_{\boldsymbol{\sigma }}\left( \boldsymbol{U},\boldsymbol{r},\boldsymbol{%
\theta }\right) :=\dprod\limits_{j=1}^{d}R_{j}^{\sigma _{j}}\left(
U_{j},r_{j},\theta _{j}\right) =\dprod\limits_{j=1}^{d}R_{j}^{\sigma _{j}}.
\end{equation*}%
For $\boldsymbol{r=}\left( r_{1},\cdots ,r_{d}\right) \in \left( 0,1\right)
^{d}$, $\boldsymbol{\sigma }=\left( \sigma _{1},\cdots ,\sigma _{d}\right)
\in \left\{ +1,-1\right\} ^{d}$ and $\boldsymbol{U}=\left( U_{1},\cdots
,U_{d}\right) ,$ we have 
\begin{eqnarray*}
\delta _{j} &=&\delta \left( U_{j},r_{j},\theta _{j}\right) =\frac{1}{2\pi }%
\left( R_{j}^{+}-R_{j}^{-}\right) =\frac{\sigma _{j}}{2\pi }\left(
R_{j}^{+\sigma _{j}}-R_{j}^{-\sigma _{j}}\right) \\
&=&\frac{\sigma _{j}\left( 1-r_{j}^{2\sigma _{j}}\right) }{2\pi }\left(
R_{j}^{\sigma _{j}}\right) ^{\ast }R_{j}^{\sigma _{j}}\geq 0.
\end{eqnarray*}%
So we define 
\begin{equation}
\gamma _{j,\sigma _{j}}\left( r_{j}\right) =\frac{\sigma _{j}\left(
1-r_{j}^{2\sigma _{j}}\right) }{2\pi }>0\text{ and }c_{\boldsymbol{\sigma }%
}\left( \boldsymbol{r}\right) =\dprod\limits_{j=1}^{d}\gamma _{j,\sigma
_{j}}\left( r_{j}\right) .  \label{p10}
\end{equation}%
Then $c_{\boldsymbol{\sigma }}\left( \boldsymbol{r}\right) $ is always
positive for $\boldsymbol{r=}\left( r_{1},\cdots ,r_{d}\right) \in \left(
0,1\right) ^{d}$. For $\boldsymbol{\sigma }=\left( \sigma _{1},\cdots
,\sigma _{d}\right) \in \left\{ +1,-1\right\} ^{d},$ we also have 
\begin{equation*}
\Delta \left( \boldsymbol{U},\boldsymbol{r},\boldsymbol{\theta }\right) =c_{%
\boldsymbol{\sigma }}\left( \boldsymbol{r}\right) R_{\boldsymbol{\sigma }%
}\left( \boldsymbol{U},\boldsymbol{r},\boldsymbol{\theta }\right) ^{\ast }R_{%
\boldsymbol{\sigma }}\left( \boldsymbol{U},\boldsymbol{r},\boldsymbol{\theta 
}\right) .
\end{equation*}

Now we are ready to define the stationary wave operators.

\begin{definition}
\label{def4}Let 
\begin{equation*}
\boldsymbol{U}=\left( U_{1},\cdots ,U_{d}\right) \text{ and }\boldsymbol{U}%
_{0}=\left( U_{0,1},\cdots ,U_{0,d}\right)
\end{equation*}%
be two $d$-tuples of commuting unitaries in a von Neumann algebra $\mathcal{M%
}$ and let $J\in \mathcal{M}.$ For any $\boldsymbol{\sigma }\in \left\{
+1,-1\right\} ^{d}$ and $\varphi _{0},\varphi $ in $\mathcal{H},$ the sector
stationary wave operator $w_{\sigma }\left( \boldsymbol{U},\boldsymbol{U}%
_{0},J\right) $ is defined as follows%
\begin{equation*}
\left\langle w_{\sigma }\left( \boldsymbol{U},\boldsymbol{U}_{0},J,\mathcal{M%
}\right) \varphi _{0},\varphi \right\rangle :=\int_{\mathbb{T}^{d}}\lim 
_{\substack{ \boldsymbol{r=}\rho \cdot \boldsymbol{1}\in \left( 0,1\right)
^{d}  \\ \rho \rightarrow 1^{-}}}\left( c_{\boldsymbol{\sigma }}\left( 
\boldsymbol{r}\right) \left\langle JR_{\boldsymbol{\sigma }}\left( 
\boldsymbol{U}_{0},\boldsymbol{r},\boldsymbol{\theta }\right) P_{ac}^{\infty
}\left( \boldsymbol{U}_{0},\mathcal{M}\right) \varphi _{0},R_{\boldsymbol{%
\sigma }}\left( \boldsymbol{U},\boldsymbol{r},\boldsymbol{\theta }\right)
\varphi \right\rangle \right) d\boldsymbol{\theta .}
\end{equation*}

where $\boldsymbol{r=}\rho \cdot \boldsymbol{1}=\left( \rho ,\cdots ,\rho
\right) \in \left( 0,1\right) ^{d}$ provided that the limit exists$.$ For $%
\boldsymbol{r=}\left( r_{1},\cdots ,r_{d}\right) \in \left( 0,1\right) ^{d},$
we define 
\begin{equation*}
\left\langle w_{\sigma }\left( \boldsymbol{U},\boldsymbol{U}_{0},J,%
\boldsymbol{r},\mathcal{M}\right) \varphi _{0},\varphi \right\rangle :=c_{%
\boldsymbol{\sigma }}\left( \boldsymbol{r}\right) \int_{\mathbb{T}%
^{d}}\left\langle JR_{\boldsymbol{\sigma }}\left( \boldsymbol{U}_{0},%
\boldsymbol{r},\boldsymbol{\theta }\right) P_{ac}^{\infty }\left( 
\boldsymbol{U}_{0},\mathcal{M}\right) \varphi _{0},R_{\boldsymbol{\sigma }%
}\left( \boldsymbol{U},\boldsymbol{r},\boldsymbol{\theta }\right) \varphi
\right\rangle d\boldsymbol{\theta .}
\end{equation*}
\end{definition}

\begin{remark}
\label{rem3}If 
\begin{eqnarray*}
&&\int_{\mathbb{T}^{d}}\lim_{\substack{ \boldsymbol{r=}\rho \cdot 
\boldsymbol{1}\in \left( 0,1\right) ^{d}  \\ \rho \rightarrow 1^{-}}}\left(
c_{\boldsymbol{\sigma }}\left( \boldsymbol{r}\right) \left\langle JR_{%
\boldsymbol{\sigma }}\left( \boldsymbol{U}_{0},\boldsymbol{r},\boldsymbol{%
\theta }\right) P_{ac}^{\infty }\left( \boldsymbol{U}_{0},\mathcal{M}\right)
\varphi _{0},R_{\boldsymbol{\sigma }}\left( \boldsymbol{U},\boldsymbol{r},%
\boldsymbol{\theta }\right) \varphi \right\rangle \right) d\boldsymbol{%
\theta } \\
&=&\lim_{\substack{ \boldsymbol{r=}\rho \cdot \boldsymbol{1}\in \left(
0,1\right) ^{d}  \\ \rho \rightarrow 1^{-}}}\int_{\mathbb{T}^{d}}\left( c_{%
\boldsymbol{\sigma }}\left( \boldsymbol{r}\right) \left\langle JR_{%
\boldsymbol{\sigma }}\left( \boldsymbol{U}_{0},\boldsymbol{r},\boldsymbol{%
\theta }\right) P_{ac}^{\infty }\left( \boldsymbol{U}_{0},\mathcal{M}\right)
\varphi _{0},R_{\boldsymbol{\sigma }}\left( \boldsymbol{U},\boldsymbol{r},%
\boldsymbol{\theta }\right) \varphi \right\rangle \right) d\boldsymbol{%
\theta ,}
\end{eqnarray*}%
then 
\begin{equation*}
\left\langle w_{\sigma }\left( \boldsymbol{U},\boldsymbol{U}_{0},J,\mathcal{M%
}\right) \varphi _{0},\varphi \right\rangle =\lim_{\substack{ \boldsymbol{r=}%
\rho \cdot \boldsymbol{1}\in \left( 0,1\right) ^{d}  \\ \rho \rightarrow
1^{-} }}\left\langle w_{\sigma }\left( \boldsymbol{U},\boldsymbol{U}_{0},J,%
\boldsymbol{r},\mathcal{M}\right) \varphi _{0},\varphi \right\rangle .
\end{equation*}%
This is the main idea behind the proof of the existence of the stationary
wave operators.
\end{remark}

Before stating the following proposition, we fix some notation. If $\sigma
=+1,$ we let 
\begin{equation*}
\mathbb{N}_{\sigma }=\mathbb{N}_{+1}=\left\{ 0,1,2,\cdots \right\} \text{.}
\end{equation*}%
If $\sigma =-1,$ we let 
\begin{equation*}
\mathbb{N_{\sigma }}=\mathbb{N}_{-1}=\left\{ 1,2,\cdots \right\} .
\end{equation*}%
Hence for any $\boldsymbol{\sigma }\in \left\{ +1,-1\right\} ^{d},$ if $%
\sigma _{j}=+1,$ $n_{j}\in \mathbb{N}_{+1}$, then 
\begin{equation}
\left( 1-r_{j}^{2}\right) r_{j}^{2\left( n_{j}-\frac{1-\sigma _{j}}{2}%
\right) }=\left( 1-r_{j}^{2}\right) r_{j}^{2n_{j}}.  \label{p11}
\end{equation}%
If $\sigma _{j}=-1,$ $n_{j}\in \mathbb{N}_{-1},$ then 
\begin{equation}
\left( 1-r_{j}^{2}\right) r_{j}^{2\left( n_{j}-\frac{1-\sigma _{j}}{2}%
\right) }=\left( 1-r_{j}^{2}\right) r_{j}^{2\left( n_{j}-1\right) }.
\label{p12}
\end{equation}

Now we characterize $w_{\boldsymbol{\sigma }}\left( \boldsymbol{U},%
\boldsymbol{U}_{0},J,\boldsymbol{r},\mathcal{M}\right) $ in the following
proposition.

\begin{proposition}
\label{pro1}Let 
\begin{equation*}
\boldsymbol{U}=\left( U_{1},\cdots ,U_{d}\right) \text{ and }\boldsymbol{U}%
_{0}=\left( U_{0,1},\cdots ,U_{0,d}\right)
\end{equation*}%
be two $d$-tuples of commuting unitaries in a von Neumann algebra $\mathcal{M%
}$ and let $J\in \mathcal{M}.$ Then for any $\boldsymbol{r}=\left(
r_{1},\cdots ,r_{d}\right) \boldsymbol{\in }\left( 0,1\right) ^{d}$ and $%
\boldsymbol{\sigma }=\left( \sigma _{1},\cdots ,\sigma _{d}\right) \in
\left\{ +1,-1\right\} ^{d},$ 
\begin{eqnarray*}
&&w_{\boldsymbol{\sigma }}\left( \boldsymbol{U},\boldsymbol{U}_{0},J,%
\boldsymbol{r},\mathcal{M}\right) \\
&=&\sum_{\boldsymbol{n}\in \mathbb{N}_{\sigma _{1}}\times \cdots \times 
\mathbb{N}_{\sigma _{d}}}\left( \left[ \Pi _{j=1}^{d}\left(
1-r_{j}^{2}\right) r_{j}^{2\left( n_{j}-\frac{1-\sigma _{j}}{2}\right) }%
\right] \cdot \boldsymbol{U}^{\boldsymbol{\sigma n}}J\boldsymbol{U}_{0}^{-%
\boldsymbol{\sigma n}}P_{ac}^{\infty }\left( \boldsymbol{U}_{0},\mathcal{M}%
\right) \right) .
\end{eqnarray*}%
Moreover, 
\begin{equation*}
\sum_{\boldsymbol{n}\in \mathbb{N}_{\sigma _{1}}\times \cdots \times \mathbb{%
N}_{\sigma _{d}}}\left[ \Pi _{j=1}^{d}\left( 1-r_{j}^{2}\right)
r_{j}^{2\left( n_{j}-\frac{1-\sigma _{j}}{2}\right) }\right] =1
\end{equation*}%
and 
\begin{equation*}
\left\Vert w_{\boldsymbol{\sigma }}\left( \boldsymbol{U},\boldsymbol{U}%
_{0},J,\boldsymbol{r},\mathcal{M}\right) \right\Vert \leq \left\Vert
J\right\Vert .
\end{equation*}
\end{proposition}

\begin{proof}
Since $\left( 1-zU^{\ast }\right) ^{-1}=\sum_{m=0}^{\infty }z^{m}U^{\ast m},$
we obtain that 
\begin{equation}
R_{0,j}^{+}:=R_{j}^{+}\left( U_{0,j},r_{j},\theta _{j}\right)
=\sum_{m_{j}=0}^{\infty }r_{j}^{m_{j}}e^{im_{j}\theta _{j}}U_{0,j}^{-m_{j}}.
\label{pr12}
\end{equation}%
Similarly, 
\begin{equation}
R_{j}^{+}:=R_{j}^{+}\left( U_{j},r_{j},\theta _{j}\right)
=\sum_{n_{j}=0}^{\infty }r_{j}^{n_{j}}e^{in_{j}\theta _{j}}U_{j}^{-n_{j}}.
\label{pr13}
\end{equation}%
We know that%
\begin{equation*}
\left( 1-zU_{j}^{\ast }\right) ^{-1}=-\sum_{n_{j=1}}^{\infty
}z^{-n_{j}}U_{j}^{n_{j}}\text{ for }\left\vert z\right\vert >1.
\end{equation*}%
Hence%
\begin{equation}
R_{j}^{-}:=R_{j}^{-}\left( U_{j},r_{j},\theta _{j}\right)
=-\sum_{m_{j}=1}^{\infty }r_{j}^{m_{j}}e^{-im_{j}\theta _{j}}U_{j}^{m_{j}}%
\text{ }  \label{pr17}
\end{equation}%
and 
\begin{equation}
R_{0,j}^{-}:=R_{j}^{-}\left( U_{0,j},r_{j},\theta _{j}\right)
=-\sum_{n_{j}=1}^{\infty }r_{j}^{n_{j}}e^{-in_{j}\theta _{j}}U_{0,j}^{n_{j}}.%
\text{ }  \label{pr18}
\end{equation}%
Then by applying (\ref{pr12}), (\ref{pr13}), (\ref{pr17}) and (\ref{pr18}),
in the following equation 
\begin{eqnarray}
&&\int_{\mathbb{T}^{d}}\left\langle JR_{\boldsymbol{\sigma }}\left( 
\boldsymbol{U}_{0},\boldsymbol{r},\boldsymbol{\theta }\right) P_{ac}^{\infty
}\left( \boldsymbol{U}_{0},\mathcal{M}\right) \varphi _{0},R_{\boldsymbol{%
\sigma }}\left( \boldsymbol{U},\boldsymbol{r}.\boldsymbol{\theta }\right)
\varphi \right\rangle d\boldsymbol{\theta }  \label{proposition2(1)} \\
&=&\int_{\mathbb{T}^{d}}\left\langle J\left(
\dprod\limits_{j=1}^{d}R_{0,j}^{\sigma _{j}}\right) P_{ac}^{\infty }\left( 
\boldsymbol{U}_{0},\mathcal{M}\right) \varphi _{0},\left(
\dprod\limits_{j=1}^{d}R_{j}^{\sigma _{j}}\right) \varphi \right\rangle d%
\boldsymbol{\theta },  \notag
\end{eqnarray}%
we have%
\begin{equation*}
R_{\boldsymbol{\sigma }}\left( \boldsymbol{U}_{0},\boldsymbol{r},\boldsymbol{%
\theta }\right) =\sum_{\boldsymbol{n}\in \mathbb{N}_{\sigma _{1}}\times
\cdots \times \mathbb{N}_{\sigma _{d}}}K_{\boldsymbol{n}}\left( \boldsymbol{U%
}_{0},\boldsymbol{\theta }\right)
\end{equation*}%
where $\left\Vert K_{\boldsymbol{n}}\left( \boldsymbol{U}_{0},\boldsymbol{%
\theta }\right) \right\Vert \leq \boldsymbol{r}^{\boldsymbol{n}}$ and 
\begin{equation*}
R_{\boldsymbol{\sigma }}\left( \boldsymbol{U},\boldsymbol{r},\boldsymbol{%
\theta }\right) =\sum_{\boldsymbol{m}\in \mathbb{N}_{\sigma _{1}}\times
\cdots \times \mathbb{N}_{\sigma _{d}}}K_{\boldsymbol{m}}\left( \boldsymbol{U%
},\boldsymbol{\theta }\right)
\end{equation*}%
where $\left\Vert K_{\boldsymbol{m}}\left( \boldsymbol{U},\boldsymbol{\theta 
}\right) \right\Vert \leq \boldsymbol{r}^{\boldsymbol{m}}.$ It follows that 
\begin{equation*}
R_{\boldsymbol{\sigma }}\left( \boldsymbol{U},\boldsymbol{r},\boldsymbol{%
\theta }\right) ^{\ast }JR_{\boldsymbol{\sigma }}\left( \boldsymbol{U}_{0},%
\boldsymbol{r},\boldsymbol{\theta }\right) P_{ac}^{\infty }\left( 
\boldsymbol{U}_{0},\mathcal{M}\right) =\sum_{\boldsymbol{m},\boldsymbol{n}%
\in \mathbb{N}_{\sigma _{1}}\times \cdots \times \mathbb{N}_{\sigma _{d}}}K_{%
\boldsymbol{m}}\left( \boldsymbol{U},\boldsymbol{\theta }\right) ^{\ast }JK_{%
\boldsymbol{n}}\left( \boldsymbol{U}_{0},\boldsymbol{\theta }\right)
P_{ac}^{\infty }\left( \boldsymbol{U}_{0},\mathcal{M}\right)
\end{equation*}%
and 
\begin{eqnarray*}
&&\sum_{\boldsymbol{m},\boldsymbol{n}\in \mathbb{N}_{\sigma _{1}}\times
\cdots \times \mathbb{N}_{\sigma _{d}}}\left\Vert K_{\boldsymbol{m}}\left( 
\boldsymbol{U},\boldsymbol{\theta }\right) ^{\ast }JK_{\boldsymbol{n}}\left( 
\boldsymbol{U}_{0},\boldsymbol{\theta }\right) P_{ac}^{\infty }\left( 
\boldsymbol{U}_{0},\mathcal{M}\right) \right\Vert \\
&\leq &\sum_{\boldsymbol{m},\boldsymbol{n}\in \mathbb{N}_{\sigma _{1}}\times
\cdots \times \mathbb{N}_{\sigma _{d}}}r^{\boldsymbol{m}+\boldsymbol{n}%
}\left\Vert J\right\Vert =\left( \sum_{\boldsymbol{m}\in \mathbb{N}_{\sigma
_{1}}\times \cdots \times \mathbb{N}_{\sigma _{d}}}r^{\boldsymbol{m}}\right)
^{2}\left\Vert J\right\Vert \\
&=&\left\Vert J\right\Vert \dprod\limits_{j=1}^{d}\frac{r_{j}^{1-\sigma _{j}}%
}{\left( 1-r_{j}\right) ^{2}}<\infty .
\end{eqnarray*}%
Since the sum of their norms is integrable over $\mathbb{T}^{d}$, Fubini's
theorem allows us to interchange the double sum and the integral. To express
(\ref{proposition2(1)}) in terms of $\boldsymbol{U}^{\boldsymbol{\sigma n}}J%
\boldsymbol{U}_{0}^{-\boldsymbol{\sigma n}},$ we should note that $%
\int_{0}^{2\pi }e^{i\left( m_{j}-l_{j}\right) \theta _{j}}d\theta _{j}=2\pi
\delta _{m_{j,}l_{j}}.$ Then by (\ref{p10}), (\ref{pr12}) and (\ref{pr13}),
if $\sigma _{j}=+1,$ then in 
\begin{equation*}
\int_{\mathbb{T}^{d}}\left\langle \gamma _{j,+1}\left( r_{j}\right)
JR_{0,j}^{+}P_{ac}^{\infty }\left( \boldsymbol{U}_{0},\mathcal{M}\right)
\varphi _{0},R_{j}^{+}\varphi \right\rangle ,
\end{equation*}%
there is a factor 
\begin{equation*}
\gamma _{j,+1}\left( r_{j}\right) r_{j}^{n_{j}}r_{j}^{n_{j}}=\frac{%
(1-r_{j}^{2})}{2\pi }r_{j}^{2\left( n_{j}-\frac{1-\sigma _{j}}{2}\right) }
\end{equation*}%
for each term involving $U^{n_{j}}JU_{0,j}^{-n_{j}}$ by (\ref{p11}) where $%
n_{j}=0,1,2,\cdots .$ Similarly for 
\begin{equation*}
\int_{\mathbb{T}^{d}}\left\langle \gamma _{j,-1}\left( r_{j}\right)
JR_{0,j}^{-}P_{ac}^{\infty }\left( \boldsymbol{U}_{0},\mathcal{M}\right)
\varphi _{0},R_{j}^{-}\varphi \right\rangle ,
\end{equation*}%
we can get a factor 
\begin{equation*}
\gamma _{j,-1}\left( r_{j}\right) r_{j}^{n_{j}}r_{j}^{n_{j}}=-\left( \frac{%
1-r_{j}^{-2}}{2\pi }\right) r_{j}^{2n_{j}}=\left( \frac{1-r_{j}^{2}}{2\pi }%
\right) r_{j}^{2\left( n_{j}-\frac{1-\sigma _{j}}{2}\right) }
\end{equation*}%
by (\ref{p12}) for the term involving $U^{-n_{j}}JU_{0,j}^{n_{j}}$. Now we
conclude that 
\begin{eqnarray*}
&&w_{\boldsymbol{\sigma }}\left( \boldsymbol{U},\boldsymbol{U}_{0},J,%
\boldsymbol{r},\mathcal{M}\right) \\
&=&\sum_{\boldsymbol{n}\in \mathbb{N}_{\sigma _{1}}\times \cdots \times 
\mathbb{N}_{\sigma _{d}}}\left( \left[ \Pi _{j=1}^{d}\left(
1-r_{j}^{2}\right) r_{j}^{2\left( n_{j}-\frac{1-\sigma _{j}}{2}\right) }%
\right] \cdot \boldsymbol{U}^{\boldsymbol{\sigma n}}J\boldsymbol{U}_{0}^{-%
\boldsymbol{\sigma n}}P_{ac}^{\infty }\left( \boldsymbol{U}_{0},\mathcal{M}%
\right) \right) .
\end{eqnarray*}%
Note 
\begin{eqnarray*}
&&\sum_{\boldsymbol{n}\in \mathbb{N}_{\sigma _{1}}\times \cdots \times 
\mathbb{N}_{\sigma _{d}}}\left[ \Pi _{j=1}^{d}\left( 1-r_{j}^{2}\right)
r_{j}^{2\left( n_{j}-\frac{1-\sigma _{j}}{2}\right) }\right] \\
&=&\Pi _{j=1}^{d}\left( \left( 1-r_{j}^{2}\right) \sum_{n_{j}=\frac{1-\sigma
_{j}}{2}}^{\infty }r_{j}^{2\left( n_{j}-\frac{1-\sigma _{j}}{2}\right)
}\right) =\Pi _{j=1}^{d}1=1,
\end{eqnarray*}%
so we also obtain that%
\begin{eqnarray*}
&&\left\Vert w_{\boldsymbol{\sigma }}\left( \boldsymbol{U},\boldsymbol{U}%
_{0},J,\boldsymbol{r},\mathcal{M}\right) \right\Vert \\
&=&\left\Vert \sum_{\boldsymbol{n}\in \mathbb{N}_{\sigma _{1}}\times \cdots
\times \mathbb{N}_{\sigma _{d}}}\left( \left[ \Pi _{j=1}^{d}\left(
1-r_{j}^{2}\right) r_{j}^{2\left( n_{j}-\frac{1-\sigma _{j}}{2}\right) }%
\right] \cdot \boldsymbol{U}^{\boldsymbol{\sigma n}}J\boldsymbol{U}_{0}^{-%
\boldsymbol{\sigma n}}P_{ac}^{\infty }\left( \boldsymbol{U}_{0},\mathcal{M}%
\right) \right) \right\Vert \\
&\leq &\sum_{\boldsymbol{n}\in \mathbb{N}_{\sigma _{1}}\times \cdots \times 
\mathbb{N}_{\sigma _{d}}}\left[ \Pi _{j=1}^{d}\left( 1-r_{j}^{2}\right)
r_{j}^{2\left( n_{j}-\frac{1-\sigma _{j}}{2}\right) }\right] \left\Vert 
\boldsymbol{U}^{\boldsymbol{\sigma n}}J\boldsymbol{U}_{0}^{-\boldsymbol{%
\sigma n}}P_{ac}^{\infty }\left( \boldsymbol{U}_{0},\mathcal{M}\right)
\right\Vert \\
&\leq &\sum_{\boldsymbol{n}\in \mathbb{N}_{\sigma _{1}}\times \cdots \times 
\mathbb{N}_{\sigma _{d}}}\left( \left[ \Pi _{j=1}^{d}\left(
1-r_{j}^{2}\right) r_{j}^{2\left( n_{j}-\frac{1-\sigma _{j}}{2}\right) }%
\right] \left\Vert J\right\Vert \right) =\left\Vert J\right\Vert .
\end{eqnarray*}%
This completes the proof.
\end{proof}

\begin{corollary}
\bigskip \label{cor1}Let $d\geq 2.$ Let 
\begin{equation*}
\boldsymbol{U}=\left( U_{1},\cdots ,U_{d}\right) \text{ and }\boldsymbol{U}%
_{0}=\left( U_{0,1},\cdots ,U_{0,d}\right)
\end{equation*}%
be two $d$-tuples of commuting unitaries in a von Neumann algebra $\mathcal{M%
}$ and let $J\in \mathcal{M}.$ Suppose $W\left( \boldsymbol{U},\boldsymbol{U}%
_{0},J,\mathcal{M}\right) $ exists. Then 
\begin{equation*}
W\left( \boldsymbol{U},\boldsymbol{U}_{0},J,\mathcal{M}\right) =sot\text{-}%
\lim_{\underset{j=1,\cdots ,d}{r_{j}\rightarrow 1^{-}}}w_{\boldsymbol{\sigma 
}}\left( \boldsymbol{U},\boldsymbol{U}_{0},J,\boldsymbol{r},\mathcal{M}%
\right)
\end{equation*}%
for any $\boldsymbol{\sigma }=\left( \sigma _{1},\cdots ,\sigma _{d}\right)
\in \left\{ +1,-1\right\} ^{d}$.
\end{corollary}

\begin{proof}
First, we define 
\begin{equation*}
W_{\boldsymbol{\sigma }}=W_{\boldsymbol{\sigma }}\left( \boldsymbol{U},%
\boldsymbol{U}_{0},J,\mathcal{M}\right) :=\underset{\overset{\boldsymbol{n}%
\in \mathbb{N}_{\sigma _{1}}\times \cdots \times \mathbb{N}_{\sigma _{d}}}{%
\boldsymbol{\min_{1\leq j\leq d}n}_{j}\rightarrow \infty }}{sot\text{-}\lim }%
\boldsymbol{U}^{\boldsymbol{\sigma n}}J\boldsymbol{U}_{0}^{-\boldsymbol{%
\sigma n}}P_{ac}^{\infty }\left( \boldsymbol{U}_{0},\mathcal{M}\right)
\end{equation*}%
for $\boldsymbol{\sigma }=\left( \sigma _{1},\cdots ,\sigma _{d}\right) \in
\left\{ +1,-1\right\} ^{d}$. Let 
\begin{equation*}
P_{\boldsymbol{r}}\left( \boldsymbol{\sigma },\boldsymbol{n}\right) =\Pi
_{j=1}^{d}\left( 1-r_{j}^{2}\right) r_{j}^{2\left( n_{j}-\frac{1-\sigma _{j}%
}{2}\right) }\cdot
\end{equation*}%
Then 
\begin{equation*}
\sum_{\boldsymbol{n}\in \mathbb{N}_{\sigma _{1}}\times \cdots \times \mathbb{%
N}_{\sigma _{d}}}P_{\boldsymbol{r}}\left( \boldsymbol{\sigma },\boldsymbol{n}%
\right) =1
\end{equation*}%
and 
\begin{equation*}
w_{\boldsymbol{\sigma }}\left( \boldsymbol{U},\boldsymbol{U}_{0},J,%
\boldsymbol{r},\mathcal{M}\right) =\sum_{\boldsymbol{n}\in \mathbb{N}%
_{\sigma _{1}}\times \cdots \times \mathbb{N}_{\sigma _{d}}}P\left( 
\boldsymbol{n}\right) \boldsymbol{U}^{\boldsymbol{\sigma n}}J\boldsymbol{U}%
_{0}^{-\boldsymbol{\sigma n}}P_{ac}^{\infty }\left( \boldsymbol{U}_{0},%
\mathcal{M}\right)
\end{equation*}%
by Proposition \ref{pro1}. Note $W_{\boldsymbol{\sigma }}\left( \boldsymbol{U%
},\boldsymbol{U}_{0},J,\mathcal{M}\right) $ exists by the assumption that $%
W\left( \boldsymbol{U},\boldsymbol{U}_{0},J,\mathcal{M}\right) $ exists. Now
for a fixed $x$ and $\varepsilon >0,$ there is an $N>0,$ such that 
\begin{equation*}
\left\Vert \left( W_{\boldsymbol{\sigma }}-\boldsymbol{U}^{\boldsymbol{%
\sigma n}}J\boldsymbol{U}_{0}^{-\boldsymbol{\sigma n}}\right) P_{ac}^{\infty
}\left( \boldsymbol{U}_{0},\mathcal{M}\right) x\right\Vert <\varepsilon
\end{equation*}%
as $n_{i}>N$ for $i=1,\cdots ,d.$ For this $N,$%
\begin{equation*}
q_{\sigma ,N}\left( \boldsymbol{r}\right) :=\sum_{\overset{\boldsymbol{n}\in 
\mathbb{N}_{\sigma _{1}}\times \cdots \times \mathbb{N}_{\sigma _{d}}}{\min
\left( n_{1},\cdots ,n_{d}\right) >N}}P_{\boldsymbol{r}}\left( \boldsymbol{%
\sigma },\boldsymbol{n}\right) =\Pi _{j=1}^{d}\left( 1-r_{j}^{2}\right)
\sum_{n_{j}>N}r_{j}^{2\left( n_{j}-\frac{1-\sigma _{j}}{2}\right) }=\Pi
_{j=1}^{d}r_{j}^{2N+1+\sigma _{j}}.
\end{equation*}%
Then 
\begin{equation*}
\lim_{\underset{j=1,\cdots ,d}{r_{j}\rightarrow 1^{-}}}q_{\sigma ,N}\left( 
\boldsymbol{r}\right) =1.
\end{equation*}%
It follows that for $\boldsymbol{r=}\left( r_{1},\cdots ,r_{d}\right) \in
\left( 0,1\right) ^{d}$ we have 
\begin{eqnarray*}
&&\left\Vert \left( w_{\boldsymbol{\sigma }}\left( \boldsymbol{U},%
\boldsymbol{U}_{0},J,\boldsymbol{r},\mathcal{M}\right) -W_{\boldsymbol{%
\sigma }}\right) x\right\Vert \\
&=&\left\Vert \sum_{\boldsymbol{n}\in \mathbb{N}_{\sigma _{1}}\times \cdots
\times \mathbb{N}_{\sigma _{d}}}P_{\boldsymbol{r}}\left( \boldsymbol{\sigma }%
,\boldsymbol{n}\right) \left( \boldsymbol{U}^{\boldsymbol{\sigma n}}J%
\boldsymbol{U}_{0}^{-\boldsymbol{\sigma n}}-W_{\boldsymbol{\sigma }}\right)
P_{ac}^{\infty }\left( \boldsymbol{U}_{0},\mathcal{M}\right) x\right\Vert \\
&\leq &q_{\sigma ,N}\left( \boldsymbol{r}\right) \sup_{\overset{\boldsymbol{n%
}\in \mathbb{N}_{\sigma _{1}}\times \cdots \times \mathbb{N}_{\sigma _{d}}}{%
\min \left( n_{1},\cdots ,n_{d}\right) >N}}\left\Vert \left( \boldsymbol{U}^{%
\boldsymbol{\sigma n}}J\boldsymbol{U}_{0}^{-\boldsymbol{\sigma n}}-W_{%
\boldsymbol{\sigma }}\right) P_{ac}^{\infty }\left( \boldsymbol{U}_{0},%
\mathcal{M}\right) x\right\Vert \\
&&+\left( 1-q_{\sigma ,N}\left( \boldsymbol{r}\right) \right) \sup_{%
\boldsymbol{n}\in \mathbb{N}_{\sigma _{1}}\times \cdots \times \mathbb{N}%
_{\sigma _{d}}}\left\Vert \left( \boldsymbol{U}^{\boldsymbol{\sigma n}}J%
\boldsymbol{U}_{0}^{-\boldsymbol{\sigma n}}-W_{\boldsymbol{\sigma }}\right)
P_{ac}^{\infty }\left( \boldsymbol{U}_{0},\mathcal{M}\right) x\right\Vert \\
&\leq &\varepsilon q_{\sigma ,N}\left( \boldsymbol{r}\right) +2\left\Vert
J\right\Vert \left\Vert x\right\Vert \left( 1-q_{\sigma ,N}\left( 
\boldsymbol{r}\right) \right) \\
&\leq &\varepsilon +2\left\Vert J\right\Vert \left\Vert x\right\Vert \left(
1-q_{\sigma ,N}\left( \boldsymbol{r}\right) \right) \\
&\rightarrow &\varepsilon \text{ as }r_{j}\rightarrow 1^{-}\text{ (}%
j=1,\cdots ,d\text{) for any }\varepsilon >0.
\end{eqnarray*}%
It follows that%
\begin{equation*}
\underset{\underset{j=1,\cdots ,d}{r_{j}\rightarrow 1^{-}}}{\lim \sup }%
\left\Vert \left( w_{\boldsymbol{\sigma }}\left( \boldsymbol{U},\boldsymbol{U%
}_{0},J,\boldsymbol{r},\mathcal{M}\right) -W_{\boldsymbol{\sigma }}\right)
x\right\Vert \leq \varepsilon .
\end{equation*}%
Therefore 
\begin{equation*}
W_{\boldsymbol{\sigma }}\left( \boldsymbol{U},\boldsymbol{U}_{0},J\right)
=sot\text{-}\lim_{\underset{j=1,\cdots ,d}{r_{j}\rightarrow 1^{-}}}w_{%
\boldsymbol{\sigma }}\left( \boldsymbol{U},\boldsymbol{U}_{0},J,\boldsymbol{r%
}\right) .
\end{equation*}%
By the fact that $W\left( \boldsymbol{U},\boldsymbol{U}_{0},J,\mathcal{M}%
\right) $ exists, we know that 
\begin{equation*}
W\left( \boldsymbol{U},\boldsymbol{U}_{0},J,\mathcal{M}\right) =W_{%
\boldsymbol{\sigma }}\left( \boldsymbol{U},\boldsymbol{U}_{0},J,\mathcal{M}%
\right) =sot\text{-}\lim_{\underset{j=1,\cdots ,d}{r_{j}\rightarrow 1^{-}}%
}w_{\boldsymbol{\sigma }}\left( \boldsymbol{U},\boldsymbol{U}_{0},J,%
\boldsymbol{r},\mathcal{M}\right) .
\end{equation*}
\end{proof}

To show the existence of the sector stationary wave operator $w_{\sigma
}\left( \boldsymbol{U},\boldsymbol{U}_{0},J,\mathcal{M}\right) ,$ we need
the following two lemmas.

\begin{lemma}
\label{lem5}Let 
\begin{equation*}
\boldsymbol{U}=\left( U_{1},\cdots ,U_{d}\right) \text{ and }\boldsymbol{U}%
_{0}=\left( U_{0,1},\cdots ,U_{0,d}\right)
\end{equation*}%
be two $d$-tuples of commuting unitaries in a von Neumann algebra $\mathcal{M%
}$ and let $J\in \mathcal{M}.$ Suppose 
\begin{equation*}
V_{k}:=JU_{0,k}-U_{k}J,\text{ }
\end{equation*}%
and 
\begin{equation*}
\text{\ }A_{j}:=R^{\sigma _{j}}\left( U_{j},r_{j},\theta _{j}\right) ,\text{
\ \ }B_{j}:=R^{\sigma _{j}}\left( U_{0,j},r_{j},\theta _{j}\right)
\end{equation*}%
for $j=1,\cdots ,d$ and $\boldsymbol{\sigma }=\left( \sigma _{1},\cdots
,\sigma _{d}\right) \in \left\{ +1,-1\right\} ^{d}.$ Then 
\begin{equation*}
R_{\boldsymbol{\sigma }}\left( \boldsymbol{U},\boldsymbol{r},\boldsymbol{%
\theta }\right) J-JR_{\boldsymbol{\sigma }}\left( \boldsymbol{U}_{0},%
\boldsymbol{r},\boldsymbol{\theta }\right) =\sum_{k=1}^{d}\left(
r_{k}^{\sigma _{k}}e^{i\theta _{k}}\right) \left( \dprod\limits_{j\leq
k}A_{j}\right) U_{k}^{\ast }V_{k}U_{0,k}^{\ast }\left( \dprod\limits_{j\geq
k}B_{j}\right) .
\end{equation*}
\end{lemma}

\begin{proof}
\bigskip Let $V_{k}:=JU_{0,k}-U_{k}J$ for $k=1,\cdots ,d$. Then 
\begin{eqnarray*}
&&R\left( U_{k},z\right) J-JR\left( U_{0,k},z\right) \\
&=&R\left( U_{k},z\right) \left[ J\left( 1-zU_{0,k}^{\ast }\right) -\left(
1-zU_{k}^{\ast }\right) J\right] R\left( U_{0,k},z\right) \\
&=&zR\left( U_{k},z\right) U_{k}^{\ast }\left( JU_{0,k}-U_{k}J\right)
U_{0,k}^{\ast }R\left( U_{0,k},z\right) \\
&=&zR\left( U_{k},z\right) U_{k}^{\ast }V_{k}U_{0,k}^{\ast }R\left(
U_{0,k},z\right) .
\end{eqnarray*}%
Therefore 
\begin{eqnarray}
&&R^{\sigma _{k}}\left( U_{k},r_{k},\theta _{k}\right) J-JR^{\sigma
_{k}}\left( U_{0,k},r_{k},\theta _{k}\right)  \label{lem10} \\
&=&\left( r_{k}^{\sigma _{k}}e^{i\theta _{k}}\right) R^{\sigma _{k}}\left(
U_{k},r_{k},\theta _{k}\right) U_{k}^{\ast }V_{k}U_{0,k}^{\ast }R^{\sigma
_{k}}\left( U_{0,k},r_{k},\theta _{k}\right) .  \notag
\end{eqnarray}%
Since 
\begin{eqnarray*}
&&A_{1}\cdots A_{d}J-JB_{1}\cdots B_{d} \\
&=&A_{1}\cdots A_{d}J-A_{1}\cdots A_{d-1}JB_{d}+A_{1}\cdots A_{d-1}JB_{d} \\
&&-A_{1}\cdots A_{d-2}JB_{d-1}B_{d}+\cdots +A_{1}JB_{2}\cdots
B_{d}-JB_{1}\cdots B_{d} \\
&=&A_{1}\cdots A_{d-1}\left( A_{d}J-JB_{d}\right) +A_{1}\cdots A_{d-2}\left(
A_{d-1}J-JB_{d-1}\right) B_{d}+ \\
&&\cdots +\left( A_{1}J-JB_{1}\right) B_{2}\cdots B_{d},
\end{eqnarray*}%
by (\ref{lem10}) 
\begin{equation*}
R_{\boldsymbol{\sigma }}\left( \boldsymbol{U},\boldsymbol{r},\boldsymbol{%
\theta }\right) J-JR_{\boldsymbol{\sigma }}\left( \boldsymbol{U}_{0},%
\boldsymbol{r},\boldsymbol{\theta }\right) =\sum_{k=1}^{d}\left(
r_{k}^{\sigma _{k}}e^{i\theta _{k}}\right) \left( \dprod\limits_{j\leq
k}A_{j}\right) U_{k}^{\ast }V_{k}U_{0,k}^{\ast }\left( \dprod\limits_{j\geq
k}B_{j}\right) .
\end{equation*}
\end{proof}

\begin{lemma}
\label{lem7}Let 
\begin{equation*}
\boldsymbol{U}=\left( U_{1},\cdots ,U_{d}\right) \text{ and }\boldsymbol{U}%
_{0}=\left( U_{0,1},\cdots ,U_{0,d}\right)
\end{equation*}%
be two $d$-tuples of commuting unitaries in a von Neumann algebra $\mathcal{M%
}$ and let $J\in \mathcal{M}.$ Suppose there exist a jointly $\boldsymbol{U}$%
-smooth operator $G_{k}$ and a jointly $\boldsymbol{U}_{0}$-smooth operator $%
G_{0,k}$ in $\mathcal{M}$ for $k=1,\cdots ,d$ such that 
\begin{equation*}
V_{k}=JU_{0,k}-U_{k}J=G_{k}^{\ast }G_{0,k}.
\end{equation*}%
Using the notation in Lemma \ref{lem5}, for $\boldsymbol{r}\in \left(
0,1\right) ^{d}$, $\boldsymbol{\theta }\in \mathbb{T}^{d}$, $\boldsymbol{%
\sigma }\in \left\{ +1,-1\right\} ^{d}$ and $\varphi _{0},\varphi $ in $%
\mathcal{H},$ define 
\begin{equation*}
X_{0,k,\boldsymbol{\sigma }}\left( \boldsymbol{r},\boldsymbol{\theta }%
,\varphi _{0}\right) :=G_{0,k}U_{0,k}^{\ast }\left( \dprod\limits_{j\geq
k}B_{j}\right) \varphi _{0},
\end{equation*}%
\begin{equation*}
Y_{k,\boldsymbol{\sigma }}\left( \boldsymbol{r},\boldsymbol{\theta },\varphi
\right) :=c_{\boldsymbol{\sigma }}\left( \boldsymbol{r}\right)
G_{k}U_{k}\left( \dprod\limits_{j\leq k}A_{j}\right) ^{\ast }R_{\boldsymbol{%
\sigma }}\left( \boldsymbol{U},\boldsymbol{r},\boldsymbol{\theta }\right)
\varphi .
\end{equation*}%
Then there are elements $X_{0,k,\boldsymbol{\sigma }}^{\partial }$ and $Y_{k,%
\boldsymbol{\sigma }}^{\partial }$ in $L^{2}\left( \mathbb{T}^{d},\mathcal{H}%
\right) $ for every $k$ such that 
\begin{eqnarray*}
&&\lim_{\underset{j=1,\cdots ,d}{r_{j}\boldsymbol{\rightarrow \boldsymbol{1}}%
^{-}}}\int_{\mathbb{T}^{d}}\left\langle X_{0,k,\boldsymbol{\sigma }}\left( 
\boldsymbol{r},\boldsymbol{\theta },\varphi _{0}\right) ,Y_{k,\boldsymbol{%
\sigma }}\left( \boldsymbol{r},\boldsymbol{\theta },\varphi \right)
\right\rangle d\boldsymbol{\theta } \\
&=&\int_{\mathbb{T}^{d}}\lim_{\underset{j=1,\cdots ,d}{r_{j}\boldsymbol{%
\rightarrow \boldsymbol{1}}^{-}}}\left\langle X_{0,k,\boldsymbol{\sigma }%
}\left( \boldsymbol{r},\boldsymbol{\theta },\varphi _{0}\right) ,Y_{k,%
\boldsymbol{\sigma }}\left( \boldsymbol{r},\boldsymbol{\theta },\varphi
\right) \right\rangle d\boldsymbol{\theta =}\int_{\mathbb{T}%
^{d}}\left\langle X_{0,k,\boldsymbol{\sigma }}^{\partial }\left( \boldsymbol{%
\theta }\right) ,Y_{k,\boldsymbol{\sigma }}^{\partial }\left( \boldsymbol{%
\theta }\right) \right\rangle d\boldsymbol{\theta .}
\end{eqnarray*}
\end{lemma}

\begin{proof}
Let $\mathbb{N}_{+1}=\left\{ 0,1,2,\cdots \right\} $ and $\mathbb{N}%
_{-1}=\left\{ 1,2,\cdots \right\} .$ Then by (\ref{pr12}) and (\ref{pr18}), 
\begin{equation*}
B_{j}:=R^{\sigma _{j}}\left( \boldsymbol{U}_{0},r_{j},\theta _{j}\right)
=\sigma _{j}\sum_{\sigma _{j}q_{j}\in \mathbb{N}_{\sigma
_{j}}}r_{j}^{|q_{j}|}e^{iq_{j}\theta _{j}}U_{0,j}^{-q_{j}}.
\end{equation*}%
Let 
\begin{equation*}
\Lambda _{k,\boldsymbol{\sigma }}=\left\{ q\in \mathbb{Z}^{d}:q_{j}=0\text{ }%
(j<k),\sigma _{j}q_{j}\in \mathbb{N}_{\sigma _{j}}(j\geq k)\right\} .
\end{equation*}%
Since $\boldsymbol{q}=\left( q_{1},\cdots ,q_{d}\right) $, $\boldsymbol{%
\theta =}\left( \theta _{1},\cdots ,\theta _{d}\right) $ and $\left\langle 
\boldsymbol{q},\boldsymbol{\theta }\right\rangle =q_{1}\theta _{1}+\cdots
+q_{d}\theta _{d},$ we have 
\begin{equation*}
X_{0,k,\boldsymbol{\sigma }}\left( \boldsymbol{r},\boldsymbol{\theta }%
\right) :=X_{0,k,\boldsymbol{\sigma }}\left( \boldsymbol{r},\boldsymbol{%
\theta },\varphi _{0}\right) =\sum_{\boldsymbol{q}\in \Lambda _{k,%
\boldsymbol{\sigma }}}\left( \boldsymbol{r}^{\left\vert \boldsymbol{q}%
\right\vert }\right) a_{q}e^{i\left\langle \boldsymbol{q},\boldsymbol{\theta 
}\right\rangle }
\end{equation*}%
where 
\begin{equation*}
a_{\boldsymbol{q}}=\kappa _{k,\boldsymbol{\sigma }}G_{0,k}U_{0,k}^{\ast }%
\boldsymbol{U}_{0}^{-\boldsymbol{q}}\varphi _{0}
\end{equation*}%
and%
\begin{equation*}
\kappa _{k,\boldsymbol{\sigma }}=\dprod\limits_{j\geq k}\sigma _{j}\in
\left\{ +1,-1\right\} .
\end{equation*}%
Since $U_{0,k}^{\ast }$ commutes with $\boldsymbol{U}_{0}^{-\boldsymbol{q}},$
we have 
\begin{equation}
\sum_{\boldsymbol{q}\in \Lambda _{k,\boldsymbol{\sigma }}}\left\Vert
a_{q}\right\Vert ^{2}\leq \sum_{\boldsymbol{q}\in \mathbb{Z}^{d}}\left\Vert
G_{0,k}\boldsymbol{U}_{0}^{-q}U_{0,k}^{\ast }\varphi _{0}\right\Vert
^{2}\leq \left( 2\pi \right) ^{d}C_{G_{0,k}}\left\Vert \varphi
_{0}\right\Vert ^{2}  \label{lem15}
\end{equation}%
where $C_{G_{0,k}}$ is the $\boldsymbol{U}_{0}$-smoothness constant of $%
G_{0,k}.$ It follows that 
\begin{eqnarray*}
\left\Vert X_{0,k,\boldsymbol{\sigma }}\left( \boldsymbol{r,\cdot }\right)
\right\Vert _{\mathcal{L}^{2}\left( \mathbb{T}^{d},\mathcal{H}\right) }^{2}
&=&\int_{\mathbb{T}^{d}}\left\Vert X_{0,k,\sigma }\left( \boldsymbol{r},%
\boldsymbol{\theta }\right) \right\Vert ^{2}d\boldsymbol{\theta } \\
&\boldsymbol{=}&\left( 2\pi \right) ^{2d}\sum_{q\in \Lambda _{k,\sigma
}}\left( \dprod\limits_{j\geq k}r_{j}^{2\left\vert q_{j}\right\vert }\right)
\left\Vert a_{q}\right\Vert ^{2}\leq \left( 2\pi \right)
^{2d}C_{G_{0,k}}\left\Vert \varphi _{0}\right\Vert ^{2}.
\end{eqnarray*}%
If we define $X_{0,k,\boldsymbol{\sigma }}^{\partial }\left( \boldsymbol{%
\theta }\right) =\sum_{\boldsymbol{q}\in \Lambda _{k,\sigma }}a_{\boldsymbol{%
q}}e^{i\left\langle \boldsymbol{q},\boldsymbol{\theta }\right\rangle },$
then $X_{0,k,\boldsymbol{\sigma }}^{\partial }\left( \boldsymbol{\theta }%
\right) \in \mathcal{L}^{2}\left( \mathbb{T}^{d},\mathcal{H}\right) $ by \ref%
{lem15}. We note that 
\begin{equation*}
X_{0,k,\boldsymbol{\sigma }}\left( \boldsymbol{r},\boldsymbol{\theta }%
\right) =P_{\boldsymbol{r}}\ast X_{0,k,\boldsymbol{\sigma }}^{\partial
}\left( \boldsymbol{\theta }\right)
\end{equation*}%
by (\ref{p21}). Hence 
\begin{equation}
X_{0,k,\boldsymbol{\sigma }}\left( \boldsymbol{r},\boldsymbol{\theta }%
\right) =P_{\boldsymbol{r}}\ast X_{0,k,\boldsymbol{\sigma }}^{\partial
}\rightarrow X_{0,k,\boldsymbol{\sigma }}^{\partial }\text{ a.e. }%
\boldsymbol{\theta }\text{ as }\boldsymbol{r}\rightarrow \boldsymbol{1}\text{
in }\mathcal{H}  \label{lem17}
\end{equation}%
by Lemma \ref{lem0}.

Since $A_{i}A_{j}=A_{j}A_{i}$ and $\gamma _{j,\sigma _{j}}\left(
r_{j}\right) A_{j}^{\ast }A_{j}=\delta _{j}\left( r_{j},\theta _{j}\right) $
where $\gamma _{j,\sigma _{j}}\left( r_{j}\right) $ is the number in (\ref%
{p10}) and $\sigma _{j}\in \left\{ +1,-1\right\} ,$ we have 
\begin{eqnarray*}
&&c_{\boldsymbol{\sigma }}\left( \boldsymbol{r}\right) \left(
\dprod\limits_{j\leq k}A_{j}\right) ^{\ast }R_{\boldsymbol{\sigma }}\left( 
\boldsymbol{U},\boldsymbol{r},\boldsymbol{\theta }\right) \\
&=&c_{\boldsymbol{\sigma }}\left( \boldsymbol{r}\right) \left(
\dprod\limits_{j\leq k}A_{j}^{\ast }A_{j}\right) \left(
\dprod\limits_{j>k}A_{j}\right) =\left( \dprod\limits_{j\leq k}\delta
_{j}\right) \left( \dprod\limits_{j>k}\gamma _{j,\sigma _{j}}\left(
r_{j}\right) A_{j}\right) .
\end{eqnarray*}%
It follows that 
\begin{eqnarray*}
Y_{k,\boldsymbol{\sigma }}\left( \boldsymbol{r},\boldsymbol{\theta }\right)
&:&=Y_{k,\boldsymbol{\sigma }}\left( \boldsymbol{r},\boldsymbol{\theta }%
,\varphi \right) =c_{\boldsymbol{\sigma }}\left( \boldsymbol{r}\right)
G_{k}U_{k}\left( \dprod\limits_{j\leq k}A_{j}\right) ^{\ast }R_{\boldsymbol{%
\sigma }}\left( \boldsymbol{U},\boldsymbol{r},\boldsymbol{\theta }\right)
\varphi \\
&=&\left( \dprod\limits_{j>k}\gamma _{j,\sigma _{j}}\left( r_{j}\right)
\right) G_{k}U_{k}\left( \dprod\limits_{j\leq k}\delta _{j}\right) \left(
\dprod\limits_{j>k}A_{j}\right) \varphi .
\end{eqnarray*}%
As $k=d,$ we have 
\begin{eqnarray*}
Y_{d,\boldsymbol{\sigma }}\left( \boldsymbol{r},\boldsymbol{\theta }\right)
&=&c_{\boldsymbol{\sigma }}\left( \boldsymbol{r}\right) G_{d}U_{d}\left(
\dprod\limits_{j\leq d}A_{j}\right) ^{\ast }R_{\boldsymbol{\sigma }}\left( 
\boldsymbol{U},\boldsymbol{r},\boldsymbol{\theta }\right) \varphi \\
&=&c_{\boldsymbol{\sigma }}\left( \boldsymbol{r}\right) G_{d}U_{d}R_{%
\boldsymbol{\sigma }}\left( \boldsymbol{U},\boldsymbol{r},\boldsymbol{\theta 
}\right) ^{\ast }R_{\boldsymbol{\sigma }}\left( \boldsymbol{U},\boldsymbol{r}%
,\boldsymbol{\theta }\right) \varphi \\
&=&G_{d}U_{d}\Delta \left( \boldsymbol{U},\boldsymbol{r},\boldsymbol{\theta }%
\right) \varphi .
\end{eqnarray*}

Note 
\begin{equation*}
\delta _{j}=\delta \left( U_{j},r_{j},\theta _{j}\right) =\frac{1}{2\pi }%
\sum_{q_{j}\in \mathbb{Z}}r_{j}^{\left\vert q_{j}\right\vert
}e^{iq_{j}\theta _{j}}U_{j}^{-q_{j}}
\end{equation*}%
by (\ref{p8}) and 
\begin{equation*}
A_{j}:=R^{\sigma _{j}}\left( U_{j},r_{j},\theta _{j}\right) =\sigma
_{j}\sum_{\sigma _{j}q_{j}\in \mathbb{N}_{\sigma
_{j}}}r_{j}^{|q_{j}|}e^{iq_{j}\theta _{j}}U_{j}^{-q_{j}}.
\end{equation*}%
Then 
\begin{equation*}
Y_{k,\boldsymbol{\sigma }}\left( \boldsymbol{r},\boldsymbol{\theta }\right)
=\left( \dprod\limits_{j>k}\gamma _{j,\sigma _{j}}\left( r_{j}\right)
\right) \sum_{\boldsymbol{n}\in \Omega _{Y}}\boldsymbol{r}^{\left\vert
n\right\vert }b_{\boldsymbol{n}}e^{i\left\langle \boldsymbol{n},\boldsymbol{%
\theta }\right\rangle }
\end{equation*}%
where%
\begin{equation*}
\Omega _{Y}=\mathbb{Z}^{k}\times \sigma _{k+1}N_{\sigma _{k+1}}\times \cdots
\times \sigma _{d}N_{\sigma _{d}}
\end{equation*}%
and 
\begin{equation*}
b_{\boldsymbol{n}}=\frac{\dprod\limits_{j=1}^{d}\sigma _{j}}{\left( 2\pi
\right) ^{k}}G_{k}U_{k}\boldsymbol{U}^{-\boldsymbol{n}}\varphi .
\end{equation*}%
Since $G_{k}$ is jointly $\boldsymbol{U}$-smooth, we have 
\begin{equation*}
\sum_{\boldsymbol{n}\in \mathbb{Z}^{d}}\left\Vert b_{\boldsymbol{n}%
}\right\Vert ^{2}\leq C_{k}\left\Vert \varphi \right\Vert ^{2}
\end{equation*}%
for a constant $C_{k}.$ Then $Y_{k,\boldsymbol{\sigma }}\left( \boldsymbol{r}%
,\boldsymbol{\theta }\right) \in \mathcal{L}^{2}\left( \mathbb{T}^{d},%
\mathcal{H}\right) .$ We define $F_{Y}=\sum_{\boldsymbol{n}\in \Omega
_{Y}}b_{\boldsymbol{n}}e^{i\left\langle \boldsymbol{n},\boldsymbol{\theta }%
\right\rangle }.$ Then 
\begin{equation*}
\left( \dprod\limits_{j>k}\gamma _{j,\sigma _{j}}\left( r_{j}\right) \right)
\cdot P_{\boldsymbol{r}}\ast F_{Y}=Y_{k,\boldsymbol{\sigma }}\left( 
\boldsymbol{r},\boldsymbol{\theta }\right)
\end{equation*}%
and 
\begin{equation*}
\left( \dprod\limits_{j>k}\gamma _{j,\sigma _{j}}\left( r_{j}\right) \right)
\rightarrow b_{k}^{\partial }=\left\{ 
\begin{array}{cc}
0, & k<d \\ 
1, & k=d%
\end{array}%
\right. .
\end{equation*}

By a similar argument for $X_{0,k,\sigma }^{\partial }$ in (\ref{lem17}), we
also conclude that 
\begin{equation*}
Y_{k,\boldsymbol{\sigma }}\left( \boldsymbol{r},\boldsymbol{\theta }\right)
=\left( \dprod\limits_{j>k}\gamma _{j,\sigma _{j}}\left( r_{j}\right)
\right) P_{\boldsymbol{r}}\ast F_{Y}\rightarrow Y_{k,\boldsymbol{\sigma }%
}^{\partial }\text{ a.e. }\boldsymbol{\theta }\text{ as }\boldsymbol{r}%
\rightarrow \boldsymbol{1}\text{ in }\mathcal{H}
\end{equation*}%
where 
\begin{equation*}
Y_{k,\boldsymbol{\sigma }}^{\partial }(\boldsymbol{\theta }):=\left\{ 
\begin{array}{cc}
0, & k<d, \\ 
\frac{1}{\left( 2\pi \right) ^{d}}\sum_{\boldsymbol{q\in }\mathbb{Z}%
^{d}}G_{d}U_{d}\boldsymbol{U}^{-\boldsymbol{q}}\cdot e^{i\left\langle 
\boldsymbol{q,\theta }\right\rangle }\varphi , & k=d.%
\end{array}%
\right.
\end{equation*}%
So we have proved that 
\begin{equation*}
Y_{k,\boldsymbol{\sigma }}\left( \boldsymbol{r,\theta }\right) \rightarrow
Y_{k,\boldsymbol{\sigma }}^{\partial }\text{ }
\end{equation*}%
a.e. $\boldsymbol{\theta }$ as $\boldsymbol{r}\rightarrow \boldsymbol{1}$ in 
$\mathcal{H}$ and therefore%
\begin{equation*}
\left\langle X_{0,k,\boldsymbol{\sigma }}\left( \boldsymbol{r,\theta }%
\right) ,Y_{k,\boldsymbol{\sigma }}\left( \boldsymbol{r,\theta }\right)
\right\rangle \rightarrow \left\langle X_{0,k,\boldsymbol{\sigma }%
}^{\partial },Y_{k,\boldsymbol{\sigma }}^{\partial }\right\rangle \text{
a.e. }\boldsymbol{\theta }\text{ as }\boldsymbol{r}\rightarrow \boldsymbol{1}
\end{equation*}%
as well as%
\begin{equation*}
\left\vert \left\langle X_{0,k,\boldsymbol{\sigma }}\left( \boldsymbol{%
r,\theta }\right) ,Y_{k,\boldsymbol{\sigma }}\left( \boldsymbol{r,\theta }%
\right) \right\rangle \right\vert \leq C\left( P^{\ast }X_{0,k,\boldsymbol{%
\sigma }}^{\partial }\right) \cdot \left( P^{\ast }F_{Y}\right)
\end{equation*}%
by Lemma \ref{lem0} for a constant $C$ and $C\left( P^{\ast }X_{0,k,%
\boldsymbol{\sigma }}^{\partial }\right) \cdot \left( P^{\ast }F_{Y}\right)
\in \mathcal{L}^{1}\left( \mathbb{T}^{d}\right) .$ Hence 
\begin{eqnarray*}
&&\lim_{\underset{j=1,\cdots ,d}{r_{j}\boldsymbol{\rightarrow \boldsymbol{1}}%
^{-}}}\int_{\mathbb{T}^{d}}\left\langle X_{0,k,\boldsymbol{\sigma }}\left( 
\boldsymbol{r},\boldsymbol{\theta },\varphi _{0}\right) ,Y_{k,\boldsymbol{%
\sigma }}\left( \boldsymbol{r},\boldsymbol{\theta },\varphi \right)
\right\rangle d\boldsymbol{\theta } \\
&=&\int_{\mathbb{T}^{d}}\lim_{\underset{j=1,\cdots ,d}{r_{j}\boldsymbol{%
\rightarrow \boldsymbol{1}}^{-}}}\left\langle X_{0,k,\boldsymbol{\sigma }%
}\left( \boldsymbol{r},\boldsymbol{\theta },\varphi _{0}\right) ,Y_{k,%
\boldsymbol{\sigma }}\left( \boldsymbol{r},\boldsymbol{\theta },\varphi
\right) \right\rangle d\boldsymbol{\theta =}\int_{\mathbb{T}%
^{d}}\left\langle X_{0,k,\boldsymbol{\sigma }}^{\partial }\left( \boldsymbol{%
\theta }\right) ,Y_{k,\boldsymbol{\sigma }}^{\partial }\left( \boldsymbol{%
\theta }\right) \right\rangle d\boldsymbol{\theta .}
\end{eqnarray*}%
for $\boldsymbol{r=}\left( r_{1},\cdots ,r_{d}\right) \in \left( 0,1\right)
^{d}$ by the dominated convergence theorem.
\end{proof}

The following proposition is a key ingredient in the proof of Theorem \ref%
{thm2}.

\begin{proposition}
\label{thm1}\bigskip Let 
\begin{equation*}
\boldsymbol{U}=\left( U_{1},\cdots ,U_{d}\right) \text{ and }\boldsymbol{U}%
_{0}=\left( U_{0,1},\cdots ,U_{0,d}\right)
\end{equation*}%
be two $d$-tuples of commuting unitaries in a von Neumann algebra $\mathcal{M%
}$ and let $J\in \mathcal{M}.$ Suppose there exist a jointly $\boldsymbol{U}$%
-smooth operator $G_{k}$ and a jointly $\boldsymbol{U}_{0}$-smooth operator $%
G_{0,k}$ in $\mathcal{M}$ for $k=1,\cdots ,d$ such that 
\begin{equation*}
V_{k}=JU_{0,k}-U_{k}J=G_{k}^{\ast }G_{0,k}.
\end{equation*}%
Then $w_{\boldsymbol{\sigma }}\left( \boldsymbol{U},\boldsymbol{U}_{0},J,%
\mathcal{M}\right) $ exists and 
\begin{equation*}
\left\langle w_{\sigma }\left( \boldsymbol{U},\boldsymbol{U}_{0},J,\mathcal{M%
}\right) \varphi _{0},\varphi \right\rangle =\lim_{\substack{ \boldsymbol{r=}%
\left( \rho ,\cdots ,\rho \right) \in \left( 0,1\right) ^{d}  \\ \rho
\rightarrow 1^{-}}}\left\langle w_{\sigma }\left( \boldsymbol{U},\boldsymbol{%
U}_{0},J,\boldsymbol{r},\mathcal{M}\right) \varphi _{0},\varphi \right\rangle
\end{equation*}%
for any $\boldsymbol{\sigma }=\left( \sigma _{1},\cdots ,\sigma _{d}\right)
\in \left\{ +1,-1\right\} ^{d}$ and $\left\Vert w_{\boldsymbol{\sigma }%
}\left( \boldsymbol{U},\boldsymbol{U}_{0},J,\mathcal{M}\right) \right\Vert
\leq \left\Vert J\right\Vert .$
\end{proposition}

\begin{proof}
By Lemma \ref{lem5}, we have 
\begin{eqnarray}
&&JR_{\boldsymbol{\sigma }}\left( \boldsymbol{U}_{0},\boldsymbol{r},%
\boldsymbol{\theta }\right)  \label{theorem6(1)} \\
&=&R_{\boldsymbol{\sigma }}\left( \boldsymbol{U},\boldsymbol{r},\boldsymbol{%
\theta }\right) J-\sum_{k=1}^{d}\left( r_{k}^{\sigma _{k}}e^{i\theta
_{k}}\right) \left( \dprod\limits_{j\leq k}A_{j}\right) U_{k}^{\ast
}V_{k}U_{0,k}^{\ast }\left( \dprod\limits_{j\geq k}B_{j}\right)  \notag \\
&=&R_{\boldsymbol{\sigma }}\left( \boldsymbol{U},\boldsymbol{r},\boldsymbol{%
\theta }\right) J-\sum_{k=1}^{d}\left( r_{k}^{\sigma _{k}}e^{i\theta
_{k}}\right) \left( \dprod\limits_{j\leq k}A_{j}\right) U_{k}^{\ast
}G_{k}^{\ast }G_{0,k}U_{0,k}^{\ast }\left( \dprod\limits_{j\geq
k}B_{j}\right) .  \notag
\end{eqnarray}%
In 
\begin{eqnarray*}
&&\left\langle w_{\sigma }\left( \boldsymbol{U,U}_{0},J,\boldsymbol{r},%
\mathcal{M}\right) \varphi _{0},\varphi \right\rangle \\
&=&c_{\boldsymbol{\sigma }}\left( \boldsymbol{r}\right) \int_{\mathbb{T}%
^{d}}\left\langle JR_{\boldsymbol{\sigma }}\left( \boldsymbol{U}_{0},%
\boldsymbol{r},\boldsymbol{\theta }\right) P_{ac}^{\infty }\left( 
\boldsymbol{U}_{0},\mathcal{M}\right) \varphi _{0},R_{\boldsymbol{\sigma }%
}\left( \boldsymbol{U},\boldsymbol{r},\boldsymbol{\theta }\right) \varphi
\right\rangle d\boldsymbol{\theta ,}
\end{eqnarray*}%
we have the following equation by applying (\ref{theorem6(1)}) 
\begin{eqnarray}
&&c_{\boldsymbol{\sigma }}\left( \boldsymbol{r}\right) \left\langle JR_{%
\boldsymbol{\sigma }}\left( \boldsymbol{U}_{0},\boldsymbol{r},\boldsymbol{%
\theta }\right) P_{ac}^{\infty }\left( \boldsymbol{U}_{0},\mathcal{M}\right)
\varphi _{0},R_{\boldsymbol{\sigma }}\left( \boldsymbol{U},\boldsymbol{r},%
\boldsymbol{\theta }\right) \varphi \right\rangle  \notag \\
&=&c_{\boldsymbol{\sigma }}\left( \boldsymbol{r}\right) \left\langle \left(
R_{\boldsymbol{\sigma }}J-\sum_{k=1}^{d}\left( r_{k}^{\sigma _{k}}e^{i\theta
_{k}}\right) \left( \dprod\limits_{j\leq k}A_{j}\right) U_{k}^{\ast
}V_{k}U_{0,k}^{\ast }\left( \dprod\limits_{j\geq k}B_{j}\right) \right)
P_{0}\varphi _{0},R_{\boldsymbol{\sigma }}\varphi \right\rangle
\label{theorem7(3)}
\end{eqnarray}%
where $R_{\boldsymbol{\sigma }}=R_{\boldsymbol{\sigma }}\left( \boldsymbol{U}%
,\boldsymbol{r},\boldsymbol{\theta }\right) $ and $P_{0}=P_{ac}^{\infty
}\left( \boldsymbol{U}_{0},\mathcal{M}\right) $. For the first inner product
in (\ref{theorem7(3)}), we want to show that%
\begin{eqnarray}
&&\int_{\mathbb{T}^{d}}\lim_{\substack{ \boldsymbol{r=}\left( \rho ,\cdots
,\rho \right) \in \left( 0,1\right) ^{d}  \\ \rho \rightarrow 1^{-}}}c_{%
\boldsymbol{\sigma }}\left( \boldsymbol{r}\right) \left\langle R_{%
\boldsymbol{\sigma }}JP_{0}\varphi _{0},R_{\boldsymbol{\sigma }}\varphi
\right\rangle d\boldsymbol{\theta }  \notag \\
&=&\lim_{\substack{ \boldsymbol{r=}\left( \rho ,\cdots ,\rho \right) \in
\left( 0,1\right) ^{d}  \\ \rho \rightarrow 1^{-}}}\int_{\mathbb{T}^{d}}c_{%
\boldsymbol{\sigma }}\left( \boldsymbol{r}\right) \left\langle R_{%
\boldsymbol{\sigma }}JP_{0}\varphi _{0},R_{\boldsymbol{\sigma }}\varphi
\right\rangle d\boldsymbol{\theta .}  \label{theorem7(1)}
\end{eqnarray}%
To show this equality, we first note that 
\begin{eqnarray*}
&&\int_{\mathbb{T}^{d}}\lim_{\substack{ \boldsymbol{r=}\left( \rho ,\cdots
,\rho \right) \in \left( 0,1\right) ^{d}  \\ \rho \rightarrow 1^{-}}}c_{%
\boldsymbol{\sigma }}\left( \boldsymbol{r}\right) \left\langle R_{%
\boldsymbol{\sigma }}JP_{0}\varphi _{0},R_{\boldsymbol{\sigma }}\varphi
\right\rangle d\boldsymbol{\theta } \\
&=&\int_{\mathbb{T}^{d}}\lim_{\substack{ \boldsymbol{r=}\left( \rho ,\cdots
,\rho \right) \in \left( 0,1\right) ^{d}  \\ \rho \rightarrow 1^{-}}}c_{%
\boldsymbol{\sigma }}\left( \boldsymbol{r}\right) \left\langle R_{%
\boldsymbol{\sigma }}^{\ast }R_{\boldsymbol{\sigma }}JP_{0}\varphi
_{0},\varphi \right\rangle d\boldsymbol{\theta } \\
&=&\int_{\mathbb{T}^{d}}\lim_{\substack{ \boldsymbol{r=}\left( \rho ,\cdots
,\rho \right) \in \left( 0,1\right) ^{d}  \\ \rho \rightarrow 1^{-}}}%
\left\langle \Delta \left( \boldsymbol{U},\boldsymbol{r},\boldsymbol{\theta }%
\right) JP_{0}\varphi _{0},\varphi \right\rangle d\boldsymbol{\theta =}%
\left\langle P_{ac}\left( \boldsymbol{U}\right) JP_{0}\varphi _{0},\varphi
\right\rangle
\end{eqnarray*}%
by (\ref{1}) where $P_{ac}\left( \boldsymbol{U}\right) $ is defined in
Section 2.2 and 
\begin{equation*}
\lim_{\substack{ \boldsymbol{r=}\rho \cdot \boldsymbol{1}\in \left(
0,1\right) ^{d}  \\ \rho \rightarrow 1^{-}}}\int_{\mathbb{T}%
^{d}}\left\langle \Delta \left( \boldsymbol{U},\boldsymbol{r},\boldsymbol{%
\theta }\right) JP_{0}\varphi _{0},\varphi \right\rangle d\boldsymbol{\theta
=}\left\langle JP_{0}\varphi _{0},\varphi \right\rangle
\end{equation*}%
by (\ref{p1}). Therefore, to show (\ref{theorem7(1)}), we only need to show 
\begin{equation}
P_{ac}\left( \boldsymbol{U}\right) JP_{0}=JP_{0}.  \label{theorem7(2)}
\end{equation}%
Let 
\begin{equation}
K:=\left( I-P_{ac}\left( \boldsymbol{U}\right) \right) J.
\label{theorem7(5)}
\end{equation}%
Then 
\begin{equation*}
KU_{0,k}-U_{k}K=\left( I-P_{ac}\left( \boldsymbol{U}\right) \right) \left(
JU_{0,k}-U_{k}J\right) =\left( I-P_{ac}\left( \boldsymbol{U}\right) \right)
G_{k}^{\ast }G_{0,k}=0
\end{equation*}%
by the fact that $G_{k}$ is jointly $\boldsymbol{U}$-smooth, Remark \ref%
{remark2} and Proposition \ref{pro11}. It follows that 
\begin{equation*}
KU_{0,k}=U_{k}K\text{ and }KU_{0,k}^{\ast }=U_{k}^{\ast }K
\end{equation*}%
for $k=1,\cdots ,d$. Therefore 
\begin{equation*}
KE_{\boldsymbol{U}_{0}}\left( S\right) =E_{\boldsymbol{U}}\left( S\right) K.
\end{equation*}%
For any $x\in \mathcal{H}_{ac}\left( \boldsymbol{U}_{0}\right) $ 
\begin{equation*}
\left\langle E_{\boldsymbol{U}}\left( S\right) Kx,Kx\right\rangle
=\left\Vert KE_{\boldsymbol{U}_{0}}\left( S\right) x\right\Vert ^{2}\leq
\left\Vert K\right\Vert ^{2}\left\langle E_{\boldsymbol{U}_{0}}\left(
S\right) x,x\right\rangle ,
\end{equation*}%
then $Kx\in \mathcal{H}_{ac}\left( \boldsymbol{U}\right) .$ It implies that%
\begin{equation*}
\left( I-P_{ac}\left( \boldsymbol{U}\right) \right)
JP_{0}=KP_{0}=P_{ac}\left( \boldsymbol{U}\right) KP_{0}=P_{ac}\left( 
\boldsymbol{U}\right) \left( I-P_{ac}\left( \boldsymbol{U}\right) \right)
JP_{0}=0
\end{equation*}%
by (\ref{theorem7(5)}). Hence we obtain (\ref{theorem7(2)}). Then the
equation (\ref{theorem7(1)}) holds.

Next, we consider the second inner product in (\ref{theorem7(3)}). By Lemma %
\ref{lem7}, we have 
\begin{eqnarray*}
&&\int_{\mathbb{T}^{d}}\lim_{\substack{ \boldsymbol{r=}\left( \rho ,\cdots
,\rho \right) \in \left( 0,1\right) ^{d}  \\ \rho \rightarrow 1^{-}}}c_{%
\boldsymbol{\sigma }}\left( \boldsymbol{r}\right) \left\langle
\sum_{k=1}^{d}\left( r_{k}^{\sigma _{k}}e^{i\theta _{k}}\right) \left(
\dprod\limits_{j\leq k}A_{j}\right) U_{k}^{\ast }G_{k}^{\ast
}G_{0,k}U_{0,k}^{\ast }\left( \dprod\limits_{j\geq k}B_{j}\right)
P_{0}\varphi _{0},R_{\boldsymbol{\sigma }}\varphi \right\rangle d\boldsymbol{%
\theta } \\
&=&\lim_{\substack{ \boldsymbol{r=}\left( \rho ,\cdots ,\rho \right) \in
\left( 0,1\right) ^{d}  \\ \rho \rightarrow 1^{-}}}\int_{\mathbb{T}^{d}}c_{%
\boldsymbol{\sigma }}\left( \boldsymbol{r}\right) \left\langle
\sum_{k=1}^{d}\left( r_{k}^{\sigma _{k}}e^{i\theta _{k}}\right) \left(
\dprod\limits_{j\leq k}A_{j}\right) U_{k}^{\ast }G_{k}^{\ast
}G_{0,k}U_{0,k}^{\ast }\left( \dprod\limits_{j\geq k}B_{j}\right)
P_{0}\varphi _{0},R_{\boldsymbol{\sigma }}\varphi \right\rangle d\boldsymbol{%
\theta }
\end{eqnarray*}%
Combining the above equality and (\ref{theorem7(1)}), we conclude that $%
w_{\sigma }\left( \boldsymbol{U,U}_{0},J,\mathcal{M}\right) $ exists and$\ $%
\begin{eqnarray*}
&&\lim_{\substack{ \boldsymbol{r=}\left( \rho ,\cdots ,\rho \right) \in
\left( 0,1\right) ^{d}  \\ \rho \rightarrow 1^{-}}}\left\langle w_{\sigma
}\left( \boldsymbol{U,U}_{0},J,\boldsymbol{r},\mathcal{M}\right) \varphi
_{0},\varphi \right\rangle \\
&=&\int_{\mathbb{T}^{d}}\lim_{\substack{ \boldsymbol{r=}\left( \rho ,\cdots
,\rho \right) \in \left( 0,1\right) ^{d}  \\ \rho \rightarrow 1^{-}}}c_{%
\boldsymbol{\sigma }}\left( \boldsymbol{r}\right) \left\langle JR_{%
\boldsymbol{\sigma }}\left( \boldsymbol{U}_{0},\boldsymbol{r},\boldsymbol{%
\theta }\right) P_{0}\varphi _{0},R_{\boldsymbol{\sigma }}\varphi
\right\rangle d\boldsymbol{\theta .}
\end{eqnarray*}%
Hence 
\begin{equation*}
w_{\sigma }\left( \boldsymbol{U,U}_{0},J,\mathcal{M}\right) =wot\text{-}\lim 
_{\substack{ \boldsymbol{r=}\left( \rho ,\cdots ,\rho \right) \in \left(
0,1\right) ^{d}  \\ \rho \rightarrow 1^{-}}}w_{\sigma }\left( \boldsymbol{U,U%
}_{0},J,\boldsymbol{r},\mathcal{M}\right) .
\end{equation*}%
Note 
\begin{equation*}
\left\Vert w_{\sigma }\left( \boldsymbol{U},\boldsymbol{U}_{0},J,\boldsymbol{%
r},\mathcal{M}\right) \right\Vert \leq \left\Vert J\right\Vert
\end{equation*}%
by Proposition \ref{pro1}, so we have 
\begin{equation*}
\left\Vert w_{\boldsymbol{\sigma }}\left( \boldsymbol{U},\boldsymbol{U}%
_{0},J,\mathcal{M}\right) \right\Vert \leq \left\Vert J\right\Vert .
\end{equation*}
\end{proof}

The following theorem is the first main result of this paper.

\begin{theorem}
\label{thm2}Let $d\geq 2$. Let 
\begin{equation*}
\boldsymbol{U}=\left( U_{1},\cdots ,U_{d}\right) \text{ and }\boldsymbol{U}%
_{0}=\left( U_{0,1},\cdots ,U_{0,d}\right)
\end{equation*}%
be two $d$-tuples of commuting unitaries in a von Neumann algebra $\mathcal{M%
}$ and let $J\in \mathcal{M}.$ Suppose there exist a jointly $\boldsymbol{U}$%
-smooth operator $G_{k}$ and a jointly $\boldsymbol{U}_{0}$-smooth operator $%
G_{0,k}$ in $\mathcal{M}$ for $k=1,\cdots ,d$ such that 
\begin{equation*}
V_{k}=JU_{0,k}-U_{k}J=G_{k}^{\ast }G_{0,k}.
\end{equation*}%
Then $W\left( \boldsymbol{U,U}_{0},J,\mathcal{M}\right) $ and $w_{%
\boldsymbol{\sigma }}\left( \boldsymbol{U},\boldsymbol{U}_{0},J,\mathcal{M}%
\right) $ exist, 
\begin{equation*}
W\left( \boldsymbol{U,U}_{0},J,\mathcal{M}\right) =w_{\boldsymbol{\sigma }%
}\left( \boldsymbol{U},\boldsymbol{U}_{0},J,\mathcal{M}\right)
\end{equation*}%
for any $\boldsymbol{\sigma }=\left( \sigma _{1},\cdots ,\sigma _{d}\right)
\in \left\{ +1,-1\right\} ^{d}.$ Moreover, 
\begin{equation*}
W\left( \boldsymbol{U,U}_{0},J,\mathcal{M}\right) U_{0,k}=U_{k}W\left( 
\boldsymbol{U,U}_{0},J,\mathcal{M}\right)
\end{equation*}%
for $k=1,\cdots ,d$.
\end{theorem}

\begin{proof}
If $W\left( \boldsymbol{U,U}_{0},J,\mathcal{M}\right) $ exists, then it is
clear that 
\begin{equation*}
W\left( \boldsymbol{U,U}_{0},J,\mathcal{M}\right) U_{0,k}=U_{k}W\left( 
\boldsymbol{U,U}_{0},J,\mathcal{M}\right)
\end{equation*}%
by the definition of $W\left( \boldsymbol{U,U}_{0},J,\mathcal{M}\right) .$
Applying Theorem \ref{thm1}, we know that $w_{\boldsymbol{\sigma }}\left( 
\boldsymbol{U},\boldsymbol{U}_{0},J,\mathcal{M}\right) $ exists and is
bounded for any $\boldsymbol{\sigma }\in \left\{ +1,-1\right\} ^{d}$. We
first prove that 
\begin{equation*}
W\left( \boldsymbol{U,U}_{0},J,\mathcal{M}\right) =w_{\boldsymbol{\sigma }%
}\left( \boldsymbol{U},\boldsymbol{U}_{0},J,\mathcal{M}\right)
\end{equation*}%
for any $\boldsymbol{\sigma }.$ By Proposition \ref{pro1} and Theorem \ref%
{thm1}, we know that%
\begin{eqnarray*}
&&w_{\boldsymbol{\sigma }}\left( \boldsymbol{U},\boldsymbol{U}_{0},J,%
\boldsymbol{r},\mathcal{M}\right) \\
&=&\sum_{\boldsymbol{n}\in \mathbb{N}_{\sigma _{1}}\times \cdots \times 
\mathbb{N}_{\sigma _{d}}}\left[ \Pi _{j=1}^{d}\left( 1-r_{j}^{2}\right)
r_{j}^{2\left( n_{j}-\frac{1-\sigma _{j}}{2}\right) }\right] \cdot 
\boldsymbol{U}^{\boldsymbol{\sigma n}}J\boldsymbol{U}_{0}^{-\boldsymbol{%
\sigma n}}P_{ac}^{\infty }\left( \boldsymbol{U}_{0}\right)
\end{eqnarray*}%
and 
\begin{equation*}
w_{\sigma }\left( \boldsymbol{U,U}_{0},J,\mathcal{M}\right) =wot\text{-}\lim 
_{\substack{ \boldsymbol{r=}\left( \rho ,\cdots ,\rho \right) \in \left(
0,1\right) ^{d}  \\ \rho \rightarrow 1^{-}}}w_{\sigma }\left( \boldsymbol{U,U%
}_{0},J,\boldsymbol{r},\mathcal{M}\right) .
\end{equation*}%
Let $T\left( \boldsymbol{n}\right) :=\boldsymbol{U}^{\boldsymbol{n}}J%
\boldsymbol{U}_{0}^{-\boldsymbol{n}}P_{ac}^{\infty }\left( \boldsymbol{U}%
_{0}\right) $ and $e_{k}=\left( 0,\cdots ,0,1,0,\cdots ,0\right) \in \mathbb{%
Z}^{d}$ denote the $k$-th standard basis vector of $\mathbb{Z}^{d}$. Then 
\begin{eqnarray*}
&&T\left( \boldsymbol{q}+e_{k}\right) -T\left( \boldsymbol{q}\right) \\
&=&\boldsymbol{U}^{\boldsymbol{q}}\left( U_{k}J-JU_{0,k}\right) U_{0,k}^{-1}%
\boldsymbol{U}_{0}^{-\boldsymbol{q}}=-\boldsymbol{U}^{\boldsymbol{q}%
}G_{k}^{\ast }G_{0,k}U_{0,k}^{-1}\boldsymbol{U}_{0}^{-\boldsymbol{q}%
}P_{ac}^{\infty }\left( \boldsymbol{U}_{0}\right) .
\end{eqnarray*}%
Let $\xi _{\boldsymbol{q}}=G_{0,k}U_{0,k}^{-1}\boldsymbol{U}_{0}^{-%
\boldsymbol{q}}P_{ac}^{\infty }\left( \boldsymbol{U}_{0}\right) x$ for $x\in 
\mathcal{H}.$ Then%
\begin{equation}
\left[ T\left( \boldsymbol{q}+e_{k}\right) -T\left( \boldsymbol{q}\right) %
\right] x=--\boldsymbol{U}^{\boldsymbol{q}}G_{k}^{\ast }\xi _{\boldsymbol{q}%
}.  \label{thm11}
\end{equation}

Let 
\begin{equation*}
\Omega _{M}=\left\{ \boldsymbol{m}\in \mathbb{Z}^{d}:\left\Vert \boldsymbol{m%
}\right\Vert _{\infty }\geq M\right\}
\end{equation*}%
and $\boldsymbol{n}\in \mathbb{N}_{\sigma _{1}}\times \cdots \times \mathbb{N%
}_{\sigma _{d}}.$ For $\boldsymbol{\sigma n}$ and $\boldsymbol{m}$ in $%
\Omega _{M},$ we choose an integer $L$ such that $L>\left\{ M,\left\Vert 
\boldsymbol{n}\right\Vert _{\infty },\left\Vert \boldsymbol{m}\right\Vert
_{\infty }\right\} $ and let $\boldsymbol{l}=\left( L,\cdots ,L\right) .$
For $\boldsymbol{\sigma n}$, we find an index $i$ such that $n_{i}\geq M.$
Now we take a path connecting $\boldsymbol{\sigma n}$ and $\boldsymbol{l}%
_{i}=\left( L,\cdots ,L,\sigma _{i}n_{i},L,\cdots ,L\right) $ in $\mathbb{Z}%
^{d}$ for a $\boldsymbol{\sigma }\in \left\{ +1,-1\right\} ^{d}.$ More
precisely, we let 
\begin{equation*}
\boldsymbol{v}^{(i,k)}:=\left( L,\cdots ,L,\sigma _{k+1}n_{k+1},\cdots
,\sigma _{d}n_{d}\right)
\end{equation*}%
for $k<i$ and 
\begin{equation*}
\boldsymbol{v}^{(i,k)}:=\left( L,\cdots ,L,\sigma _{i}n_{i},L,\cdots
,L,\sigma _{k+1}n_{k+1},\cdots ,\sigma _{d}n_{d}\right)
\end{equation*}%
for $k\geq i$ with $\boldsymbol{v}^{(i,0)}=\boldsymbol{\sigma n}$ and 
\begin{equation*}
\boldsymbol{v}^{(i,d)}=\boldsymbol{l}_{i}=\left( L,\cdots ,L,\sigma
_{i}n_{i},L,\cdots ,L\right) .
\end{equation*}%
Now for the index $k$ with $k\neq i$ and an integer $t$ with $\sigma
_{k}n_{k}\leq t\leq L-1$ ($k\neq i$), we write 
\begin{equation*}
\boldsymbol{q}_{k,t}^{i}=v^{\left( i,k-1\right) }+\left( t-\sigma
_{k}n_{k}\right) e_{k}\in \Omega _{M},\text{ }
\end{equation*}%
then for any $k\neq i$ 
\begin{eqnarray*}
&&\left\Vert T\left( \boldsymbol{v}^{(i,k)}\right) -T\left( \boldsymbol{v}%
^{(i,k-1)}\right) x\right\Vert =\left\Vert \sum_{t=\sigma _{k}n_{k}}^{L-1}%
\left[ T\left( \boldsymbol{q}_{k,t}^{i}+e_{k}\right) -T_{N}\left( 
\boldsymbol{q}_{k,t}^{i}\right) \right] \right\Vert \\
&=&\left\Vert \sum_{t=\sigma _{k}n_{k}}^{L-1}\boldsymbol{U}^{\boldsymbol{q}%
_{k,t}^{i}}G_{k}^{\ast }\xi _{\boldsymbol{q}_{k,t}^{i}}\right\Vert
\end{eqnarray*}%
by (\ref{thm11}) and for a unit vector $y$ 
\begin{eqnarray*}
&&\left\vert \left\langle \sum_{t=\sigma _{k}n_{k}}^{L-1}\boldsymbol{U}^{%
\boldsymbol{q}_{k,t}^{i}}G_{k}^{\ast }\xi _{\boldsymbol{q}%
_{k,t}^{i}},y\right\rangle \right\vert =\left\vert \sum_{t=\sigma
_{k}n_{k}}^{L-1}\left\langle \xi _{\boldsymbol{q}_{k,t}^{i}},G_{k}%
\boldsymbol{U}^{-\boldsymbol{q}_{k,t}^{i}}y\right\rangle \right\vert \\
&\leq &\left( \sum_{t=\sigma _{k}n_{k}}^{L-1}\left\Vert \xi _{\boldsymbol{q}%
_{k,t}^{i}}\right\Vert ^{2}\right) ^{1/2}\left( \sum_{t=\sigma
_{k}n_{k}}^{L-1}\left\Vert G_{k}\boldsymbol{U}^{-\boldsymbol{q}%
_{k,t}^{i}}y\right\Vert ^{2}\right) ^{1/2} \\
&\leq &\left( \sum_{t=\sigma _{k}n_{k}}^{L-1}\left\Vert G_{0,k}U_{0,k}^{-1}%
\boldsymbol{U}_{0}^{-\boldsymbol{q}_{k,t}^{i}}P_{ac}^{\infty }\left( 
\boldsymbol{U}_{0}\right) x\right\Vert ^{2}\right) ^{1/2}\sqrt{\left( 2\pi
\right) ^{d}C_{G_{k}}}
\end{eqnarray*}%
by the fact that $G_{k}$ is $\boldsymbol{U}$-smooth. It follows that 
\begin{eqnarray*}
&&\left\Vert T\left( \boldsymbol{v}^{(i,k)}\right) -T\left( \boldsymbol{v}%
^{(i,k-1)}\right) x\right\Vert =\left\Vert \sum_{t=\sigma _{k}n_{k}}^{L-1}%
\boldsymbol{U}^{\boldsymbol{q}_{k,t}^{i}}G_{k}^{\ast }\xi _{\boldsymbol{q}%
_{k,t}^{i}}\right\Vert \\
&\leq &\sqrt{\left( 2\pi \right) ^{d}C_{G_{k}}}\left( \sum_{t=\sigma
_{k}n_{k}}^{L-1}\left\Vert G_{0,k}U_{0,k}^{-1}\boldsymbol{U}_{0}^{-%
\boldsymbol{q}_{k,t}^{i}}P_{ac}^{\infty }\left( \boldsymbol{U}_{0}\right)
x\right\Vert ^{2}\right) ^{1/2} \\
&\leq &\sqrt{\left( 2\pi \right) ^{d}C_{G_{k}}}\cdot a_{k,M}\left( x\right)
\end{eqnarray*}%
where 
\begin{equation*}
a_{k,M}\left( x\right) =\left( \sum_{\boldsymbol{q}\in \Omega
_{M}}\left\Vert G_{0,k}U_{0,k}^{\ast }\boldsymbol{U}_{0}^{-\boldsymbol{q}%
}P_{ac}^{\infty }\left( \boldsymbol{U}_{0}\right) x\right\Vert ^{2}\right)
^{1/2}.
\end{equation*}%
Similarly, we also have 
\begin{equation*}
\left\Vert \left[ T\left( \boldsymbol{l}_{i}\right) -T\left( \boldsymbol{l}%
\right) \right] x\right\Vert \leq \sqrt{\left( 2\pi \right) ^{d}C_{G_{i}}}%
a_{i,M}\left( x\right) .
\end{equation*}%
By the following equality 
\begin{eqnarray*}
&&\left[ T\left( \boldsymbol{l}\right) -T(\boldsymbol{\sigma n})\right] x \\
&=&\sum_{k=1,k\neq i}^{d}\left[ T\left( \boldsymbol{v}^{(i,k)}\right)
-T\left( \boldsymbol{v}^{(i,k-1)}\right) \right] x+\left[ T\left( 
\boldsymbol{l}\right) -T\left( \boldsymbol{l}_{i}\right) \right] x,
\end{eqnarray*}%
we conclude that 
\begin{equation*}
\left\Vert \left[ T\left( \boldsymbol{l}\right) -T(\boldsymbol{\sigma n})%
\right] x\right\Vert \leq E_{M}\left( x\right)
\end{equation*}%
where%
\begin{equation*}
E_{M}\left( x\right) =\sum_{k=1}^{d}\sqrt{\left( 2\pi \right) ^{d}C_{G_{k}}}%
\cdot a_{k,M}\left( x\right) .
\end{equation*}%
Similarly, we also have 
\begin{equation*}
\left\Vert \left[ T\left( \boldsymbol{l}\right) -T(\boldsymbol{m})\right]
x\right\Vert \leq E_{M}\left( x\right) .
\end{equation*}%
It implies that 
\begin{equation*}
\left\Vert \left[ T\left( \boldsymbol{m}\right) -T(\boldsymbol{\sigma n})%
\right] x\right\Vert \leq 2E_{M}\left( x\right) .
\end{equation*}%
Since $G_{0,k}$ is $\boldsymbol{U}_{0}$-smooth, we conclude that 
\begin{equation*}
a_{k,M}\left( x\right) \leq \left( \sum_{\boldsymbol{q}\in \mathbb{Z}%
^{d}}\left\Vert G_{0,k}U_{0,k}^{\ast }\boldsymbol{U}_{0}^{-\boldsymbol{q}%
}P_{ac}^{\infty }\left( \boldsymbol{U}_{0}\right) x\right\Vert ^{2}\right)
^{1/2}<\infty .
\end{equation*}%
Hence $a_{k,M}\left( x\right) \rightarrow 0$ for $k=1,\cdots ,d$ as $%
M\rightarrow \infty .$ It implies that 
\begin{equation*}
E_{M}\left( x\right) \rightarrow 0\text{ as }M\rightarrow \infty .
\end{equation*}

If $\boldsymbol{m}\notin \Omega _{M},$ then 
\begin{equation*}
\left\Vert \left[ T\left( \boldsymbol{m}\right) -T\left( \boldsymbol{\sigma n%
}\right) \right] x\right\Vert \leq 2\left\Vert J\right\Vert \left\Vert
x\right\Vert .
\end{equation*}%
By Proposition \ref{pro1}, we have%
\begin{eqnarray*}
&&\left\Vert T\left( \boldsymbol{m}\right) -w_{\boldsymbol{\sigma }}\left( 
\boldsymbol{U},\boldsymbol{U}_{0},J,\boldsymbol{r},\mathcal{M}\right)
x\right\Vert \\
&=&\left\Vert \sum_{\boldsymbol{n}\in \mathbb{N}_{\sigma _{1}}\times \cdots
\times \mathbb{N}_{\sigma _{d}}}\left( \left[ \Pi _{j=1}^{d}\left(
1-r_{j}^{2}\right) r_{j}^{2\left( n_{j}-\frac{1-\sigma _{j}}{2}\right) }%
\right] \left[ T\left( \boldsymbol{m}\right) -T\left( \boldsymbol{\sigma n}%
\right) \right] \right) x\right\Vert \\
&\leq &2E_{M}\left( x\right) +2\left\Vert J\right\Vert \left\Vert
x\right\Vert \eta _{\boldsymbol{\sigma ,}M}\left( \boldsymbol{r}\right) 
\text{ for }\boldsymbol{m}\in \Omega _{M}
\end{eqnarray*}%
where 
\begin{equation*}
\eta _{\boldsymbol{\sigma ,}M}\left( \boldsymbol{r}\right) :=\sum_{\underset{%
\boldsymbol{n\notin }\Omega _{M}}{\boldsymbol{n}\in \mathbb{N}_{\sigma
_{1}}\times \cdots \times \mathbb{N}_{\sigma _{d}}}}\left[ \Pi
_{j=1}^{d}\left( 1-r_{j}^{2}\right) r_{j}^{2\left( n_{j}-\frac{1-\sigma _{j}%
}{2}\right) }\right] .
\end{equation*}%
Then $\eta _{\boldsymbol{\sigma ,}M}\left( \boldsymbol{r}\right) \rightarrow
0$ as $\boldsymbol{r=}\left( \rho ,\cdots ,\rho \right) \in \left(
0,1\right) ^{d}$ and $\rho \rightarrow 1^{-}.$ It implies that 
\begin{eqnarray*}
&&\left\vert \left\langle T\left( \boldsymbol{m}\right) x-w_{\boldsymbol{%
\sigma }}\left( \boldsymbol{U},\boldsymbol{U}_{0},J,\mathcal{M}\right)
x,y\right\rangle \right\vert \\
&=&\left\vert \lim_{\substack{ \boldsymbol{r=}\left( \rho ,\cdots ,\rho
\right) \in \left( 0,1\right) ^{d}  \\ \rho \rightarrow 1^{-}}}\left\langle
T\left( \boldsymbol{m}\right) x-w_{\boldsymbol{\sigma }}\left( \boldsymbol{U}%
,\boldsymbol{U}_{0},J,\boldsymbol{r},\mathcal{M}\right) x,y\right\rangle
\right\vert \\
&\leq &\lim_{\substack{ \boldsymbol{r=}\left( \rho ,\cdots ,\rho \right) \in
\left( 0,1\right) ^{d}  \\ \rho \rightarrow 1^{-}}}\left( 2E_{M}\left(
x\right) +2\left\Vert J\right\Vert \left\Vert x\right\Vert \eta _{%
\boldsymbol{\sigma ,}M}\left( \boldsymbol{r}\right) \right) \left\Vert
y\right\Vert \\
&=&2E_{M}\left( x\right) \left\Vert y\right\Vert .
\end{eqnarray*}%
Taking the supremum over all unit vectors $y\in \mathcal{H},$ we obtain that 
\begin{eqnarray*}
&&\lim_{\left\Vert \boldsymbol{m}\right\Vert _{\infty }\rightarrow \infty
}\left\Vert T\left( \boldsymbol{m}\right) x-w_{\boldsymbol{\sigma }}\left( 
\boldsymbol{U},\boldsymbol{U}_{0},J,\mathcal{M}\right) x\right\Vert \\
&=&\lim_{\left\Vert \boldsymbol{m}\right\Vert _{\infty }\rightarrow \infty
}\left\Vert \boldsymbol{U}^{\boldsymbol{m}}J_{N}\boldsymbol{U}_{0}^{-%
\boldsymbol{m}}P_{0,ac}^{\infty }\left( \boldsymbol{U}_{0}\right) x-w_{%
\boldsymbol{\sigma }}\left( \boldsymbol{U},\boldsymbol{U}_{0},J,\mathcal{M}%
\right) x\right\Vert \\
&=&0.
\end{eqnarray*}%
Hence the strong wave operator $W\left( \boldsymbol{U,U}_{0},J,\mathcal{M}%
\right) $ exists and 
\begin{equation*}
W\left( \boldsymbol{U,U}_{0},J,\mathcal{M}\right) =w_{\boldsymbol{\sigma }%
}\left( \boldsymbol{U},\boldsymbol{U}_{0},J,\mathcal{M}\right) .
\end{equation*}
\end{proof}

In \cite{Y}, a strong wave operator $W$ is said to be complete if it is a
partial isometry whose initial and final spaces coincide with the
corresponding absolutely continuous subspaces. The following theorem
establishes an analogous completeness property for $W\left( \boldsymbol{U,U}%
_{0},J,\mathcal{M}\right) $.

\begin{theorem}
\label{thm3}Let $d\geq 2$. Let 
\begin{equation*}
\boldsymbol{U}=\left( U_{1},\cdots ,U_{d}\right) \text{ and }\boldsymbol{U}%
_{0}=\left( U_{0,1},\cdots ,U_{0,d}\right)
\end{equation*}%
be two $d$-tuples of commuting unitaries in a von Neumann algebra $\mathcal{M%
}$ and let $J\in \mathcal{M}$ be a unitary$.$ Suppose there exist a jointly $%
\boldsymbol{U}$-smooth operator $G_{k}$ and a jointly $\boldsymbol{U}_{0}$%
-smooth operator $G_{0,k}$ in $\mathcal{M}$ for $k=1,\cdots ,d$ such that 
\begin{equation*}
V_{k}=JU_{0,k}-U_{k}J=G_{k}^{\ast }G_{0,k}.
\end{equation*}%
Then $W\left( \boldsymbol{U,U}_{0},J,\mathcal{M}\right) $ and $w_{%
\boldsymbol{\sigma }}\left( \boldsymbol{U},\boldsymbol{U}_{0},J,\mathcal{M}%
\right) $ exist, and 
\begin{equation*}
W\left( \boldsymbol{U,U}_{0},J,\mathcal{M}\right) ^{\ast }=W\left( 
\boldsymbol{U}_{0},\boldsymbol{U},J^{\ast },\mathcal{M}\right) =w_{%
\boldsymbol{\sigma }}\left( \boldsymbol{U}_{0},\boldsymbol{U},J^{\ast },%
\mathcal{M}\right) .
\end{equation*}%
Moreover, 
\begin{equation*}
W\left( \boldsymbol{U,U}_{0},J,\mathcal{M}\right) ^{\ast }W\left( 
\boldsymbol{U,U}_{0},J,\mathcal{M}\right) =P_{ac}^{\infty }\left( 
\boldsymbol{U}_{0},\mathcal{M}\right) ,
\end{equation*}%
\begin{equation*}
W\left( \boldsymbol{U,U}_{0},J,\mathcal{M}\right) W\left( \boldsymbol{U,U}%
_{0},J,\mathcal{M}\right) ^{\ast }=P_{ac}^{\infty }\left( \boldsymbol{U},%
\mathcal{M}\right) .
\end{equation*}
\end{theorem}

\begin{proof}
By Theorem \ref{thm2}, we know $W\left( \boldsymbol{U,U}_{0},J,\mathcal{M}%
\right) $ exists. To simplify the notation, we let $W=W\left( \boldsymbol{U,U%
}_{0},J,\mathcal{M}\right) $. Since $J$ is a unitary, we have 
\begin{equation*}
\left\Vert \boldsymbol{U}^{\boldsymbol{n}}J\boldsymbol{U}_{0}^{-\boldsymbol{n%
}}P_{ac}^{\infty }\left( \boldsymbol{U}_{0},\mathcal{M}\right) x\right\Vert
=\left\Vert J\boldsymbol{U}_{0}^{-\boldsymbol{n}}P_{ac}^{\infty }\left( 
\boldsymbol{U}_{0},\mathcal{M}\right) x\right\Vert =\left\Vert
P_{ac}^{\infty }\left( \boldsymbol{U}_{0},\mathcal{M}\right) x\right\Vert .
\end{equation*}%
It follows that $\left\Vert Wx\right\Vert =\left\Vert P_{ac}^{\infty }\left( 
\boldsymbol{U}_{0},\mathcal{M}\right) x\right\Vert .$ Therefore 
\begin{equation*}
\left\langle W^{\ast }Wx,x\right\rangle =\left\Vert Wx\right\Vert
^{2}=\left\Vert P_{ac}^{\infty }\left( \boldsymbol{U}_{0},\mathcal{M}\right)
x\right\Vert ^{2}=\left\langle P_{ac}^{\infty }\left( \boldsymbol{U}_{0},%
\mathcal{M}\right) x,x\right\rangle
\end{equation*}%
for all $x.$ It implies that 
\begin{equation}
W^{\ast }W=P_{ac}^{\infty }\left( \boldsymbol{U}_{0},\mathcal{M}\right) .
\label{theorem9(1)}
\end{equation}%
Let 
\begin{equation*}
V=W\left( \boldsymbol{U}_{0}\boldsymbol{,U},J^{\ast },\mathcal{M}\right) =sot%
\text{-}\lim_{\left\Vert \mathbb{m}\right\Vert _{\infty }\rightarrow \infty }%
\boldsymbol{U}_{0}^{\boldsymbol{m}}J^{\ast }\boldsymbol{U}^{-\boldsymbol{m}%
}P_{ac}^{\infty }\left( \boldsymbol{U},\mathcal{M}\right) .
\end{equation*}%
Since 
\begin{equation}
U_{0,k}^{\ast }J^{\ast }-J^{\ast }U_{k}^{\ast }=\left(
JU_{0,k}-U_{k}J\right) ^{\ast }=\left( G_{k}^{\ast }G_{0,k}\right) ^{\ast
}=G_{0,k}^{\ast }G_{k},  \label{thm9(2)}
\end{equation}%
multiplying $U_{0,k}$ from the left and $U_{k}$ from the right in (\ref%
{thm9(2)}), we get that 
\begin{equation*}
J^{\ast }U_{k}-U_{0,k}J^{\ast }=U_{0,k}G_{0,k}^{\ast }G_{k}U_{k}=\left(
G_{0,k}U_{0,k}^{\ast }\right) ^{\ast }\left( G_{k}U_{k}\right) .
\end{equation*}%
By an argument similar to that in the proof of Theorem \ref{thm2}, we
conclude that $V$ exists, $V^{\ast }V=P_{ac}^{\infty }\left( \boldsymbol{U},%
\mathcal{M}\right) $ and 
\begin{equation*}
V=w_{\boldsymbol{\sigma }}\left( \boldsymbol{U}_{0}\boldsymbol{,U},J^{\ast },%
\mathcal{M}\right) .
\end{equation*}

Let $P$ be a jointly $\boldsymbol{U}_{0}$-smooth projection and set $%
Q_{P}=WPW^{\ast }.$ Then $Q_{P}^{\ast }=Q_{P}$ and 
\begin{equation*}
Q_{P}^{2}=WPW^{\ast }WPW^{\ast }=WPP_{ac}^{\infty }\left( \boldsymbol{U}_{0},%
\mathcal{M}\right) PW^{\ast }=WPW^{\ast }=Q_{P}
\end{equation*}%
by (\ref{theorem9(1)}). Since $W^{\ast }\boldsymbol{U}^{\boldsymbol{a}}=%
\boldsymbol{U}_{0}^{\boldsymbol{a}}W^{\ast }$ for any $\boldsymbol{a}\in 
\mathbb{Z}^{d}$ by Theorem \ref{thm2} and 
\begin{equation*}
\left\Vert Q_{P}\boldsymbol{U}^{\boldsymbol{a}}y\right\Vert =\left\Vert
WPW^{\ast }\boldsymbol{U}^{\boldsymbol{a}}y\right\Vert =\left\Vert WP%
\boldsymbol{U}_{0}^{\boldsymbol{a}}W^{\ast }y\right\Vert =\left\Vert P%
\boldsymbol{U}_{0}^{\boldsymbol{a}}W^{\ast }y\right\Vert ,
\end{equation*}%
we have 
\begin{equation*}
\sum_{\boldsymbol{a}\in \mathbb{Z}^{d}}\left\Vert Q_{P}\boldsymbol{U}^{%
\boldsymbol{a}}y\right\Vert ^{2}=\sum_{\boldsymbol{a}\in \mathbb{Z}%
^{d}}\left\Vert P\boldsymbol{U}_{0}^{\boldsymbol{a}}W^{\ast }y\right\Vert
^{2}\leq \left( 2\pi \right) ^{d}C_{P}\left\Vert W^{\ast }y\right\Vert
^{2}\leq \left( 2\pi \right) ^{d}C_{P}\left\Vert y\right\Vert ^{2}
\end{equation*}%
by the fact that $P$ is jointly $\boldsymbol{U}_{0}$-smooth where $C_{P}$ is
the smoothness constant of $P.$ So we conclude that $Q_{P}$ is a jointly $%
\boldsymbol{U}$-smooth projection and then 
\begin{equation*}
Q_{P}\leq P_{ac}^{\infty }\left( \boldsymbol{U},\mathcal{M}\right) .
\end{equation*}%
Note that, by (\ref{theorem9(1)}), 
\begin{eqnarray*}
WW^{\ast }\left( \mathcal{H}\right) &=&WP_{ac}^{\infty }\left( \boldsymbol{U}%
_{0},\mathcal{M}\right) W^{\ast }\left( \mathcal{H}\right) \\
&=&\overline{span\left( \cup _{P\in \mathcal{P}_{ac}^{\infty }\left( 
\boldsymbol{U}_{0},\mathcal{M}\right) }Ran(WPW^{\ast }\right) }=\overline{%
span\cup _{P\in \mathcal{P}_{ac}^{\infty }\left( \boldsymbol{U}_{0},\mathcal{%
M}\right) }Ran\left( Q_{P}\right) } \\
&=&\left( \bigvee_{P\in \mathcal{P}_{ac}^{\infty }\left( \boldsymbol{U}_{0},%
\mathcal{M}\right) }Q_{P}\right) \left( \mathcal{H}\right) =\left(
\bigvee_{P\in \mathcal{P}_{ac}^{\infty }\left( \boldsymbol{U}_{0},\mathcal{M}%
\right) }WPW^{\ast }\right) \left( \mathcal{H}\right) \subseteq
P_{ac}^{\infty }\left( \boldsymbol{U},\mathcal{M}\right) \left( \mathcal{H}%
\right) .
\end{eqnarray*}%
It follows that 
\begin{equation}
P_{ac}^{\infty }\left( \boldsymbol{U},\mathcal{M}\right) W=W.
\label{thm3(3)}
\end{equation}%
By a similar argument, we conclude that 
\begin{equation}
P_{ac}^{\infty }\left( \boldsymbol{U}_{0},\mathcal{M}\right) V=V.
\label{thm3(4)}
\end{equation}%
Next, we are going to show $V=W^{\ast }.$ Let $T\left( \boldsymbol{n}\right)
=\boldsymbol{U}^{\boldsymbol{n}}J\boldsymbol{U}_{0}^{-\boldsymbol{n}}$. Then
by \ref{thm3(3)} and \ref{thm3(4)}, 
\begin{eqnarray*}
\left\langle Wx,y\right\rangle &=&\left\langle Wx,P_{ac}^{\infty }\left( 
\boldsymbol{U},\mathcal{M}\right) y\right\rangle =\lim_{\left\Vert 
\boldsymbol{n}\right\Vert _{\infty }\rightarrow \infty }\left\langle T\left( 
\boldsymbol{n}\right) P_{ac}^{\infty }\left( \boldsymbol{U}_{0},\mathcal{M}%
\right) x,P_{ac}^{\infty }\left( \boldsymbol{U},\mathcal{M}\right)
y\right\rangle \\
&=&\lim_{\left\Vert \boldsymbol{n}\right\Vert _{\infty }\rightarrow \infty
}\left\langle P_{ac}^{\infty }\left( \boldsymbol{U}_{0},\mathcal{M}\right)
x,B\left( \boldsymbol{n}\right) y\right\rangle =\left\langle P_{ac}^{\infty
}\left( \boldsymbol{U}_{0},\mathcal{M}\right) x,Vy\right\rangle \\
&=&\left\langle x,P_{ac}^{\infty }\left( \boldsymbol{U}_{0},\mathcal{M}%
\right) Vy\right\rangle =\left\langle x,Vy\right\rangle .
\end{eqnarray*}%
Therefore $V=W^{\ast }.$ It implies that 
\begin{equation*}
WW^{\ast }=V^{\ast }V=P_{ac}^{\infty }\left( \boldsymbol{U},\mathcal{M}%
\right) .
\end{equation*}%
The proof is completed.
\end{proof}

\section{Wave Operator of $U$ and $U_{0}$ in $\mathcal{M}$ ($d=1$)}

Let $\mathcal{M\subseteq B(H)}$ be a von Neumann algebra. In this section,
we consider the case $d=1$. In the following equations, "$\pm "$ is
understood as two independent equalities.

\begin{definition}
Let $U$ and $U_{0}$ be two unitaries and let $J\in \mathcal{M}$.

\begin{enumerate}
\item The strong wave operators of a triple $\left( U,U_{0},J\right) $ in $%
\mathcal{M}$ are defined below: 
\begin{equation*}
W_{\pm }\left( U,U_{0},J,\mathcal{M}\right) :=sot\text{-}\lim_{n\rightarrow
\pm \infty }U^{n}JU_{0}^{-n}P_{ac}^{\infty }\left( U_{0},\mathcal{M}\right)
\in \mathcal{M}
\end{equation*}%
provided that the limits exist.

\item The weak wave operators of a triple $\left( U,U_{0},J\right) $ in $%
\mathcal{M}$ are defined below:%
\begin{equation*}
\widetilde{W}_{\pm }\left( U,U_{0},J,\mathcal{M}\right) :=wot\text{-}%
\lim_{n\rightarrow \pm \infty }P_{ac}^{\infty }\left( U,\mathcal{M}\right)
U^{n}JU_{0}^{-n}P_{ac}^{\infty }\left( U_{0},\mathcal{M}\right)
\end{equation*}%
provided that the limits exist.

\item The modified stationary wave operators $\widetilde{w}_{\pm }\left(
U,U_{0},J,\mathcal{M}\right) $ of a triple $\left( U,U_{0},J\right) $ in $%
\mathcal{M}$ are defined below%
\begin{equation*}
\left\langle \widetilde{w}_{\pm }\left( U,U_{0},J,\mathcal{M}\right) \varphi
_{0},\varphi \right\rangle :=\int_{\mathbb{T}}\lim_{\underset{r\in (0,1)}{%
r\rightarrow 1}}\widetilde{\alpha }_{\pm }^{\mathcal{M}}\left( \varphi
_{0},\varphi ,r,\theta \right) d\theta
\end{equation*}%
provided that $\lim_{\underset{r\in (0,1)}{r\rightarrow 1}}\widetilde{\alpha 
}_{\pm }^{\mathcal{M}}\left( \varphi _{0},\varphi ,r,\theta \right) $ exist
where%
\begin{equation*}
\widetilde{\alpha }_{\pm }^{\mathcal{M}}\left( \varphi _{0},\varphi
,r,\theta \right) :=\gamma _{\pm }\left( r\right) \left\langle JR^{\pm
}\left( U_{0},r,\theta \right) P_{ac}^{\infty }\left( U_{0},\mathcal{M}%
\right) \varphi _{0},R^{\pm }\left( U,r,\theta \right) P_{ac}^{\infty
}\left( U,\mathcal{M}\right) \varphi \right\rangle
\end{equation*}%
for $\varphi _{0},\varphi \in \mathcal{H}$ and $\gamma _{\pm }\left(
r\right) $ are defined in (\ref{p9}) for $r\in \left( 0,1\right) .$ We also
define 
\begin{equation*}
\widetilde{w}_{+}\left( U,U_{0},J,r,\mathcal{M}\right) :=\int_{\mathbb{T}}%
\widetilde{\alpha }_{\pm }^{\mathcal{M}}\left( \varphi _{0},\varphi
,r,\theta \right) d\theta .
\end{equation*}
\end{enumerate}
\end{definition}

\begin{remark}
Define 
\begin{equation*}
\alpha _{\pm }^{\mathcal{M}}\left( \varphi _{0},\varphi ,r,\theta \right)
:=\gamma _{\pm }\left( r\right) \left\langle JR^{\pm }\left( U_{0},r,\theta
\right) P_{ac}^{\infty }\left( U_{0},\mathcal{M}\right) \varphi _{0},R^{\pm
}\left( U,r,\theta \right) \varphi \right\rangle .
\end{equation*}%
By the definition of the stationary wave operator $w_{\pm }\left( U,U_{0},J,%
\mathcal{M}\right) $ of a triple $\left( U,U_{0},J\right) $ defined in
Definition \ref{def4}$,$ we know that 
\begin{equation*}
\left\langle w_{\pm }\left( U,U_{0},J,\mathcal{M}\right) \varphi
_{0},\varphi \right\rangle :=\int_{\mathbb{T}}\lim_{\underset{r\in (0,1)}{%
r\rightarrow 1}}\alpha _{\pm }^{\mathcal{M}}\left( \varphi _{0},\varphi
,r,\theta \right) d\theta
\end{equation*}%
It is clear 
\begin{equation*}
\widetilde{\alpha }_{\pm }^{\mathcal{M}}\left( \varphi _{0},\varphi
,r,\theta \right) =\alpha _{\pm }^{\mathcal{M}}\left( \varphi
_{0},P_{ac}^{\infty }\left( U,\mathcal{M}\right) \varphi ,r,\theta \right) .
\end{equation*}%
By Proposition \ref{pro1},%
\begin{equation}
w_{+}\left( U,U_{0},J,r,\mathcal{M}\right) =\left( 1-r^{2}\right)
\sum_{n\geq 0}r^{2n}U^{n}JU_{0}^{-n}P_{ac}^{\infty }\left( U_{0},\mathcal{M}%
\right)  \label{16}
\end{equation}%
and 
\begin{equation}
w_{-}\left( U,U_{0},J,r,\mathcal{M}\right) =\left( 1-r^{2}\right)
\sum_{n\geq 1}r^{2(n-1)}U^{-n}JU_{0}^{n}P_{ac}^{\infty }\left( U_{0},%
\mathcal{M}\right) .  \label{15}
\end{equation}%
Hence we have 
\begin{equation}
\widetilde{w}_{+}\left( U,U_{0},J,r,\mathcal{M}\right) =\left(
1-r^{2}\right) \sum_{n\geq 0}r^{2n}P_{ac}^{\infty }\left( U,\mathcal{M}%
\right) U^{n}JU_{0}^{-n}P_{ac}^{\infty }\left( U_{0},\mathcal{M}\right)
\label{def3}
\end{equation}%
and 
\begin{equation}
\widetilde{w}_{-}\left( U,U_{0},J,r,\mathcal{M}\right) =\left(
1-r^{2}\right) \sum_{n\geq 1}r^{2(n-1)}P_{ac}^{\infty }\left( U,\mathcal{M}%
\right) U^{-n}JU_{0}^{n}P_{ac}^{\infty }\left( U_{0},\mathcal{M}\right) .
\label{def(2)}
\end{equation}
\end{remark}

\begin{remark}
\label{rem11}By an argument similar to that in the proof of Lemma 3.6(a) in 
\cite{TR} or Theorem 2.2.1 in \cite{Y}, we conclude that $W_{\pm }\left(
U,U_{0},J,\mathcal{M}\right) $ exists if and only if the weak wave operators 
$\widetilde{W}_{\pm }\left( U,U_{0},J,\mathcal{M}\right) $ and $\widetilde{W}%
_{\pm }\left( U_{0},U_{0},J^{\ast }J,\mathcal{M}\right) $ exist and satisfy
the chain rule 
\begin{equation*}
\widetilde{W}_{\pm }\left( U,U_{0},J,\mathcal{M}\right) ^{\ast }\widetilde{W}%
_{\pm }\left( U,U_{0},J,\mathcal{M}\right) =\widetilde{W}_{\pm }\left(
U_{0},U_{0},J^{\ast }J,\mathcal{M}\right)
\end{equation*}%
in which case we have the identity 
\begin{equation*}
W_{\pm }\left( U,U_{0},J,\mathcal{M}\right) =\widetilde{W}_{\pm }\left(
U,U_{0},J,\mathcal{M}\right) .
\end{equation*}
\end{remark}

Next we give our first main result in this section. It is a result analogous
to Theorem \ref{thm2} for the case $d=1.$

\begin{theorem}
\label{th3} Let $U$, $U_{0}$ be two unitaries in $\mathcal{M}$ and let $J\in 
\mathcal{M}.$ Suppose there are a $U_{0}$-smooth operator $G_{0}$ and a $U$%
-smooth operator $G$ such that 
\begin{equation*}
V:=JU_{0}-UJ=G^{\ast }G_{0}.
\end{equation*}%
Then $W_{\pm }\left( U,U_{0},J,\mathcal{M}\right) $ and $\omega _{\pm
}(U,U_{0},J,\mathcal{M})$ exist and 
\begin{equation*}
W_{\pm }\left( U,U_{0},J,\mathcal{M}\right) =\omega _{\pm }(U,U_{0},J,%
\mathcal{M})
\end{equation*}%
where "$\pm "$ is understood as two independent equalities. Moreover, 
\begin{equation*}
W_{\pm }\left( U,U_{0},J,\mathcal{M}\right) U_{0}=UW_{\pm }\left( U,U_{0},J,%
\mathcal{M}\right) .
\end{equation*}
\end{theorem}

\begin{proof}
By Proposition \ref{thm1}, applied with $d=1,$ the stationary wave operators 
$\omega _{\pm }(U,U_{0},J,\mathcal{M})$ exist and 
\begin{equation}
w_{\pm }\left( U,U_{0},J,\mathcal{M}\right) =wot\text{-}\lim_{\underset{r\in
\left( 0,1\right) }{r\rightarrow 1}}w_{\pm }\left( U,U_{0},J,r,\mathcal{M}%
\right) .  \label{2}
\end{equation}%
We now prove the existence of the strong wave operators and identify them
with these stationary wave operators. We define operators $A$ and $B$ from $%
\mathcal{H}$ to $l^{2}\left( \mathbb{Z},\mathcal{H}\right) $ as follows 
\begin{equation*}
\left( Ay\right) _{n}=GU^{-n}y\text{ and }\left( Bx\right)
_{n}=G_{0}U_{0}^{-n-1}x.
\end{equation*}%
The smoothness assumptions imply that $A$ and $B$ are bounded, with 
\begin{equation*}
\left\Vert A\right\Vert \leq \sqrt{2\pi C_{G}}\text{ and }\left\Vert
B\right\Vert \leq \sqrt{2\pi C_{G_{0}}}.
\end{equation*}%
For every finitely supported sequence $v=\left( v_{n}\right) ,$%
\begin{equation*}
A^{\ast }v=\sum_{n}U^{n}G^{\ast }v_{n}.
\end{equation*}%
For $E\subseteq \mathbb{Z},$ let $\Pi _{E}$ denote the orthogonal projection
in $l^{2}\left( \mathbb{Z},\mathcal{H}\right) $ onto the subspace of
sequences supported in $E.$ We write 
\begin{equation*}
\Pi _{+}=\Pi _{\left\{ n\geq 0\right\} }\text{ and }\Pi _{-}=\Pi _{\left\{
n\leq -1\right\} }
\end{equation*}%
which are the projections from $l^{2}\left( \mathbb{Z},\mathcal{H}\right) $
onto the subspaces of sequences supported in $n\geq 0$ and $n\leq -1$
respectively. For integers $a\leq b,$ we then have 
\begin{equation*}
A^{\ast }\Pi _{\lbrack a,b]}Bx=\sum_{n=a}^{b}U^{n}G^{\ast
}G_{0}U_{0}^{-n-1}x.
\end{equation*}%
Let $T_{n}=U^{n}JU_{0}^{-n}$ and then $T_{0}=J.$ Then 
\begin{equation}
T_{n+1}-T_{n}=-U^{n}G^{\ast }G_{0}U_{0}^{-n-1}  \label{35}
\end{equation}%
and 
\begin{eqnarray*}
T_{N}-J &=&\left( T_{1}-T_{0}\right) +\cdots +\left( T_{N}-T_{N-1}\right) \\
&=&-\sum_{n=0}^{N-1}U^{n}G^{\ast }G_{0}U_{0}^{-n-1}.
\end{eqnarray*}%
It follows that 
\begin{equation*}
T_{N}=J-A^{\ast }\Pi _{\lbrack 0,N-1]}B\text{ and similarly }%
T_{-N}=J+A^{\ast }\Pi _{\lbrack -N,-1]}B.
\end{equation*}%
As $N\rightarrow \infty ,$%
\begin{equation*}
\Pi _{\lbrack 0,N-1]}\rightarrow \Pi _{+}\text{ and }\Pi _{\lbrack
-N,-1]}\rightarrow \Pi _{-}
\end{equation*}%
strongly on $l^{2}\left( \mathbb{Z},\mathcal{H}\right) .$ Consequently, the
strong wave operators exist and satisfy%
\begin{equation*}
W_{+}\left( U,U_{0},J,\mathcal{M}\right) =sot\text{-}\lim_{N\rightarrow
\infty }T_{N}P_{ac}^{\infty }\left( U_{0},\mathcal{M}\right) =\left(
J-A^{\ast }\Pi _{+}B\right) P_{ac}^{\infty }\left( U_{0},\mathcal{M}\right)
\end{equation*}%
and 
\begin{equation*}
W_{-}\left( U,U_{0},J,\mathcal{M}\right) =\underset{N\rightarrow \infty }{sot%
\text{-}\lim }T_{-N}P_{ac}^{\infty }\left( U_{0},\mathcal{M}\right) =\left(
J+A^{\ast }\Pi _{-}B\right) P_{ac}^{\infty }\left( U_{0},\mathcal{M}\right) .
\end{equation*}%
By (\ref{16}) and (\ref{15}), applying the proof of Corollary \ref{cor1}, we
obtain that 
\begin{equation*}
sot\text{-}\lim_{\underset{r\in \left( 0,1\right) }{r\rightarrow 1}}w_{\pm
}\left( U,U_{0},J,r,\mathcal{M}\right) =W_{\pm }\left( U,U_{0},J,\mathcal{M}%
\right) .
\end{equation*}%
It implies that 
\begin{equation*}
w_{\pm }\left( U,U_{0},J,\mathcal{M}\right) =W_{\pm }\left( U,U_{0},J,%
\mathcal{M}\right)
\end{equation*}%
by (\ref{2}). Finally, since $P_{ac}^{\infty }\left( U_{0},\mathcal{M}%
\right) $ commutes with $U_{0},$%
\begin{equation*}
T_{n}P_{ac}^{\infty }\left( U_{0},\mathcal{M}\right)
U_{0}=UT_{n-1}P_{ac}^{\infty }\left( U_{0},\mathcal{M}\right) .
\end{equation*}%
Taking the strong limit as $n\rightarrow \infty ,$ we conclude that 
\begin{equation*}
W_{\pm }\left( U,U_{0},J,\mathcal{M}\right) U_{0}=UW_{\pm }\left( U,U_{0},J,%
\mathcal{M}\right) .
\end{equation*}%
This completes the proof.
\end{proof}

The following theorem is an analogue of Theorem \ref{thm3}. It shows that $%
W_{\pm }\left( U,U_{0},J,\mathcal{M}\right) $ are complete if $J$ is a
unitary in $\mathcal{M}$.

\begin{theorem}
\label{th4}Let $U$, $U_{0}$ be two unitaries in $\mathcal{M}$ and let $J$ be
a unitary in $\mathcal{M}.$ Suppose there are a $U_{0}$-smooth operator $%
G_{0}$ and a $U$-smooth operator $G$ such that 
\begin{equation*}
V:=JU_{0}-UJ=G^{\ast }G_{0}.
\end{equation*}%
Then $W_{\pm }\left( U,U_{0},J,\mathcal{M}\right) $ and $\omega _{\pm
}(U,U_{0},J,\mathcal{M})$ exist and%
\begin{equation*}
W_{\pm }=\omega _{\pm }(U,U_{0},J,\mathcal{M}).
\end{equation*}%
Moreover, 
\begin{equation*}
W_{\pm }^{\ast }W_{\pm }=P_{ac}^{\infty }\left( U_{0},\mathcal{M}\right) 
\text{, \ }W_{\pm }W_{\pm }^{\ast }=P_{ac}^{\infty }\left( U,\mathcal{M}%
\right)
\end{equation*}%
where "$\pm "$ is understood as two independent equalities.
\end{theorem}

\begin{proof}
\bigskip By Theorem \ref{th3} and a similar argument used in the proof of
Theorem \ref{thm3}, we obtain the desired conclusion.
\end{proof}

\begin{example}
\bigskip This example shows that the strong wave operators $W_{+}\left(
U,U_{0},J,\mathcal{M}\right) $ and $W_{-}\left( U,U_{0},J,\mathcal{M}\right) 
$ need not coincide.

Let 
\begin{equation*}
\mathcal{M=L}^{\infty }\left( \left[ 0,1\right] \right) \otimes \mathcal{B}%
\left( l^{2}\left( \mathbb{Z}\right) \right) \text{ and }\mathcal{H}=%
\mathcal{L}^{2}\left( \left[ 0,1\right] \right) \otimes l^{2}\left( \mathbb{Z%
}\right) .
\end{equation*}%
Let $S$ be the bilateral shift on $l^{2}\left( \mathbb{Z}\right) ,$ defined
by $Se_{k}=e_{k+1}$ where $\left\{ e_{k}\right\} $ is the standard
orthonormal basis in $l^{2}\left( \mathbb{Z}\right) $. We set 
\begin{equation*}
U=U_{0}=I\otimes S\text{ \ and \ }\lambda \left( t\right) =e^{2\pi it}.
\end{equation*}%
Let $f\in $ $\mathcal{L}^{\infty }\left( \left[ 0,1\right] \right) .$ Let $%
M_{f}$ denote multiplication by $f$ on $\mathcal{L}^{2}\left( \left[ 0,1%
\right] \right) .$ We define a unitary $J\in \mathcal{M}$ by%
\begin{equation*}
J(t)e_{k}=\left\{ 
\begin{array}{cc}
e_{k}, & k\geq 0, \\ 
\lambda \left( t\right) e_{k}, & k<0.%
\end{array}%
\right.
\end{equation*}%
Let $E_{ij}$ denote the matrix unit determined by $E_{ij}e_{k}=\delta
_{jk}e_{i}.$ A direct calculation gives 
\begin{equation*}
\left( JU_{0}-UJ\right) =M_{1-\lambda }\otimes E_{0,-1}=G^{\ast }G_{0}
\end{equation*}%
where $G=I\otimes E_{00}=I\otimes \left( e_{0}\otimes e_{0}\right) $ and $%
G_{0}=M_{1-\lambda }\otimes \left( e_{0}\otimes e_{-1}\right) $. So it is
easy to see that $G$ and $G_{0}$ are both in $\mathcal{M.}$ Identifying $%
\mathcal{H}$ with $\mathcal{L}^{2}\left( \left[ 0,1\right] ,l^{2}\left( 
\mathbb{Z}\right) \right) ,$ we obtain, for any $x\in \mathcal{H},$%
\begin{equation*}
\sum_{n\in \mathbb{z}}\left\Vert GU^{n}x\right\Vert ^{2}=\left\Vert
x\right\Vert ^{2},
\end{equation*}%
and 
\begin{equation*}
\sum_{n\in \mathbb{Z}}\left\Vert G_{0}U_{0}^{n}x\right\Vert
^{2}=\int_{0}^{1}\left\vert 1-\lambda \left( t\right) \right\vert
^{2}\left\Vert x\right\Vert _{l^{2}\left( \mathbb{Z}\right) }^{2}dt\leq
4\left\Vert x\right\Vert ^{2}.
\end{equation*}%
Hence $G$ is $U$-smooth and $G_{0}$ is $U_{0}$-smooth. Moreover, $I\otimes
\left( e_{k}\otimes e_{k}\right) $ is $U_{0}$-smooth for every $k\in \mathbb{%
Z}$ in $\mathcal{M}$. Hence $P_{ac}^{\infty }\left( U_{0},\mathcal{M}\right)
=I.$ Now for every $n,k\in \mathbb{Z},$ we have 
\begin{equation*}
\left( U^{n}JU_{0}^{-n}\right) \left( t\right) e_{k}=\left\{ 
\begin{array}{cc}
e_{k}, & k\geq n, \\ 
\lambda \left( t\right) e_{k}, & k<n.%
\end{array}%
\right.
\end{equation*}%
Then $W_{+}\left( U,U_{0},J,\mathcal{M}\right) =M_{\lambda }\otimes I$ and 
\begin{equation*}
W_{-}\left( U,U_{0},J,\mathcal{M}\right) =I\neq W_{+}\left( U,U_{0},J,%
\mathcal{M}\right) .
\end{equation*}
\end{example}

The proof of the following lemma is based on Theorem \ref{th4}.

\begin{lemma}
\label{lem6}Let $\mathcal{M\subseteq B}\left( \mathcal{H}\right) $ be a
countably decomposable, properly infinite von Neumann algebra with a
faithful normal semifinite tracial weight $\tau .$ Let $U$, $U_{0}$ be two
unitaries in $\mathcal{M}$ and $J$ be a unitary in $\mathcal{M}.$ If 
\begin{equation*}
V:=JU_{0}-UJ\in \mathcal{L}^{1}\left( \mathcal{M},\tau \right) \mathcal{\cap
M,}
\end{equation*}%
then the modified stationary wave operators $\widetilde{w}_{\pm }=\widetilde{%
w}_{\pm }\left( U,U_{0},J,\mathcal{M}\right) $ exist and 
\begin{equation*}
\widetilde{w}_{\pm }=wot\text{-}\lim_{\underset{r\in (0,1)}{r\rightarrow 1}%
}P_{ac}^{\infty }\left( U,\mathcal{M}\right) w_{\pm }\left( U,U_{0},r,J,%
\mathcal{M}\right) .
\end{equation*}
\end{lemma}

\begin{proof}
If $V=v\left\vert V\right\vert $ is the polar decomposition of $V,$ we set $%
G_{0}=\left\vert V\right\vert ^{1/2}$ and $G=\left\vert V\right\vert
^{1/2}v^{\ast }.$ Then $V=G^{\ast }G_{0}$ with $G,$ $G_{0}\in \mathcal{L}%
^{2}\left( \mathcal{M}\right) \mathcal{\cap M}.$ By Lemma [\ref{trace}],
there exists a sequence $\left\{ x_{l}\right\} _{l\in \mathbb{N}\text{ }}$
of vectors in $\mathcal{H}$ such that 
\begin{equation*}
\tau \left( X^{\ast }X\right) =\Sigma _{l\geq 1}\left\langle X^{\ast
}Xx_{l},x_{l}\right\rangle ,\text{ \ \ \ }\forall X\in \mathcal{M}
\end{equation*}%
and the linear span of $\left\{ A^{\prime }x_{l}:A^{\prime }\in \mathcal{M}%
^{\prime }\right\} $ is dense in $\mathcal{H}.$ Let%
\begin{equation*}
\mathcal{D}_{0}=span\left\{ PA^{\prime }x_{l}:P\in \mathcal{P}_{ac}^{\infty
}\left( U_{0},\mathcal{M}\right) ,\text{ }A^{\prime }\in \mathcal{M}^{\prime
},l\geq 1\right\} .
\end{equation*}%
Then $\mathcal{D}_{0}$ is dense in $P_{ac}^{\infty }\left( U_{0},\mathcal{M}%
\right) \left( \mathcal{H}\right) .$ Let $P\in \mathcal{P}_{ac}^{\infty
}\left( U_{0},\mathcal{M}\right) .$ Then for any $z\in \mathcal{H,}$%
\begin{equation*}
\sum_{k\in \mathbb{Z}}\left\Vert PU_{0}^{k}z\right\Vert ^{2}\leq 2\pi
C_{P}\left\Vert z\right\Vert ^{2}.
\end{equation*}%
Let $S_{\mathbb{Z}}=\sum_{k\in \mathbb{Z}}U_{0}^{-k}PU_{0}^{k}$. Then for
any $z\in \mathcal{H}$ 
\begin{equation*}
\left\langle S_{\mathbb{Z}}z,z\right\rangle =\sum_{k\in \mathbb{Z}%
}\left\langle U_{0}^{-k}PU_{0}^{k}z,z\right\rangle =\sum_{k\in \mathbb{Z}%
}\left\Vert PU_{0}^{k}z\right\Vert ^{2}\leq 2\pi C_{P}\left\Vert
z\right\Vert ^{2}.
\end{equation*}%
It follows that 
\begin{equation}
S_{\mathbb{Z}}\leq 2\pi C_{P}I.  \label{9.9}
\end{equation}%
If $x=PD^{\prime }x_{l}$ for $D^{\prime }\in \mathcal{M}^{\prime }$ and $%
P\in \mathcal{P}_{ac}^{\infty }\left( U_{0},\mathcal{M}\right) $ and $B\in 
\mathcal{M},$ then by (\ref{9.9}) we have 
\begin{eqnarray*}
\sum_{k\in \mathbb{Z}}\left\Vert G_{0}BU_{0}^{-k}x\right\Vert ^{2}
&=&\sum_{k\in \mathbb{Z}}\left\Vert G_{0}BU_{0}^{-k}PD^{\prime
}x_{l}\right\Vert ^{2}\leq \left\Vert D^{\prime }\right\Vert ^{2}\sum_{k\in 
\mathbb{Z}}\left\Vert G_{0}BU_{0}^{-k}Px_{l}\right\Vert ^{2} \\
&\leq &\left\Vert D^{\prime }\right\Vert ^{2}\sum_{k\in \mathbb{Z}%
}\left\Vert G_{0}BU_{0}^{-k}P\right\Vert _{2,\tau }^{2}=\left\Vert D^{\prime
}\right\Vert ^{2}\tau \left( G_{0}BS_{\mathbb{Z}}B^{\ast }G_{0}\right) \\
&\leq &2\pi C_{P}\left\Vert D^{\prime }\right\Vert ^{2}\tau \left(
G_{0}BB^{\ast }G_{0}^{\ast }\right) \leq 2\pi C_{P}\left\Vert D^{\prime
}\right\Vert ^{2}\left\Vert B\right\Vert ^{2}\left\Vert G_{0}\right\Vert
_{2,\tau }^{2}.
\end{eqnarray*}%
So for every $x\in \mathcal{D}_{0},$ there is a constant $c_{x}$ such that 
\begin{equation}
\left( \sum_{k\in \mathbb{Z}}\left\Vert G_{0}BU_{0}^{-k}x\right\Vert
^{2}\right) ^{1/2}\leq c_{x}\left\Vert B\right\Vert <\infty .
\label{thm5(1)}
\end{equation}%
Let $\mathcal{D}$ be a dense subspace of $P_{ac}^{\infty }\left( U,\mathcal{M%
}\right) \left( \mathcal{H}\right) $ which is defined in the same way as $%
\mathcal{D}_{0}.$ The corresponding estimate in (\ref{thm5(1)}) for $G$ and $%
U$ is obtained in the same way.

Now we define an operator $\Gamma _{0}:P_{ac}^{\infty }\left( U_{0},\mathcal{%
M}\right) \left( \mathcal{H}\right) \rightarrow l^{2}\left( \mathbb{Z,}%
\mathcal{H}\right) $ by 
\begin{equation*}
\Gamma _{0}x=\left( G_{0}U_{0}^{n}x\right) _{n\in \mathbb{Z}}\text{ }
\end{equation*}%
with the domain 
\begin{equation*}
Dom\left( \Gamma _{0}\right) =\left\{ x\in P_{ac}^{\infty }\left( U_{0},%
\mathcal{M}\right) \left( \mathcal{H}\right) :\sum_{k\in \mathbb{Z}%
}\left\Vert G_{0}U_{0}^{-k}x\right\Vert ^{2}<\infty \right\} .
\end{equation*}

By (\ref{thm5(1)}), $\Gamma _{0}$ is densely defined. We note that $\Gamma
_{0}$ is closed because convergence in $l^{2}\left( \mathbb{Z,}\mathcal{H}%
\right) $ implies convergence in every coordinate, and each $G_{0}U_{0}^{n}$
is bounded. We also have 
\begin{equation*}
\Gamma _{0}U_{0}x=\left( G_{0}U_{0}^{n+1}x\right) _{n\in \mathbb{Z}}=S\Gamma
_{0}x\text{ for }x\in Dom\left( \Gamma _{0}\right)
\end{equation*}%
where $\left( Sy\right) _{n}=y_{n+1}.$ It follows that 
\begin{equation}
\Gamma _{0}U_{0}|_{Dom\left( \Gamma _{0}\right) }=S\Gamma _{0}.
\label{thm5(2)}
\end{equation}%
For a unitary $V\in \mathcal{M}^{\prime }$ and $y=\left( y_{n}\right) \in
l^{2}\left( \mathbb{Z,}\mathcal{H}\right) ,$ we define $\widetilde{V}$ on $%
l^{2}\left( \mathbb{Z,}\mathcal{H}\right) $ by $\widetilde{V}y=\left(
Vy_{n}\right) _{n}$. Hence 
\begin{equation*}
\Gamma _{0}Vx=\left( G_{0}U_{0}^{n}Vx\right) _{n\in \mathbb{Z}}=\left(
VG_{0}U_{0}^{n}x\right) _{n\in \mathbb{Z}}=\left( \widetilde{V}\Gamma
_{0}x\right) _{n}=\widetilde{V}\Gamma _{0}x.
\end{equation*}%
It implies that 
\begin{equation*}
\Gamma _{0}V|_{Dom\left( \Gamma _{0}\right) }=\widetilde{V}\Gamma _{0}.
\end{equation*}%
By (\ref{thm5(2)}),%
\begin{equation*}
U_{0}^{\ast }\Gamma _{0}^{\ast }\Gamma _{0}U_{0}=\Gamma _{0}^{\ast }S^{\ast
}S\Gamma _{0}=\Gamma _{0}^{\ast }\Gamma _{0}.
\end{equation*}%
Then $E_{\Gamma _{0}^{\ast }\Gamma _{0}}\left( \Omega \right)
U_{0}|_{Dom\left( \Gamma _{0}\right) }=U_{0}|_{Dom\left( \Gamma _{0}\right)
}E_{\Gamma _{0}^{\ast }\Gamma _{0}}\left( \Omega \right) $ for any Borel set 
$\Omega .$ Similarly, we have 
\begin{equation*}
E_{\Gamma _{0}^{\ast }\Gamma _{0}}\left( \Omega \right) V|_{Dom\left( \Gamma
_{0}\right) }=V|_{Dom\left( \Gamma _{0}\right) }E_{\Gamma _{0}^{\ast }\Gamma
_{0}}\left( \Omega \right) .
\end{equation*}%
Let $P_{m}=E_{\Gamma _{0}^{\ast }\Gamma _{0}}\left( [0,m]\right) \oplus 0.$
Then $P_{m}$ commutes with $U_{0}$ and $P_{m}\uparrow P_{ac}^{\infty }\left(
U_{0},\mathcal{M}\right) $. Moreover, 
\begin{equation*}
\sum_{n\in \mathbb{Z}}\left\Vert G_{0}P_{m}U_{0}^{n}x\right\Vert
^{2}=\left\Vert \Gamma _{0}P_{m}x\right\Vert ^{2}\leq m\left\Vert
x\right\Vert ^{2}.
\end{equation*}%
It implies that $G_{0}P_{m}$ is $U_{0}$-smooth. The same construction for $%
G, $ $U,$ and $P_{ac}^{\infty }\left( U,\mathcal{M}\right) $ gives a
projection $Q_{m}$ such that $Q_{m}\uparrow P_{ac}^{\infty }\left( U,%
\mathcal{M}\right) $ and commutes with $U$. By a similar argument, we
conclude that $GQ_{m}$ is $U$-smooth.

For $J_{m}=Q_{m}JP_{m},$ 
\begin{equation*}
J_{m}U_{0}-UJ_{m}=\left( GQ_{m}\right) ^{\ast }\left( G_{0}P_{m}\right) .
\end{equation*}%
Then by Theorem \ref{th3}, $W_{\pm }\left( U,U_{0},J_{m},\mathcal{M}\right) $
and $w_{\pm }\left( U,U_{0},J_{m},\mathcal{M}\right) $ both exist and 
\begin{equation}
W_{\pm ,m}:=W_{\pm }\left( U,U_{0},J_{m},\mathcal{M}\right) =w_{\pm }\left(
U,U_{0},J_{m},\mathcal{M}\right) .  \label{thm5(11)}
\end{equation}%
We also have 
\begin{equation}
w_{\pm }\left( U,U_{0},J_{m},r,\mathcal{M}\right) =Q_{m}w_{\pm }\left(
U,U_{0},J,r,\mathcal{M}\right) P_{m}  \label{thm5(12)}
\end{equation}%
and 
\begin{equation*}
\left\Vert w_{\pm }\left( U,U_{0},J_{m},r,\mathcal{M}\right) \right\Vert
\leq \left\Vert J\right\Vert \text{ for every }m.
\end{equation*}%
We now prove the existence of the modified stationary wave operators. For $%
x\in P_{ac}^{\infty }\left( U_{0},\mathcal{M}\right) \left( \mathcal{H}%
\right) ,$ we set 
\begin{equation}
q_{U_{0},r}\left( x,\theta \right) =\gamma _{\pm }\left( r\right) \left\Vert
R^{\pm }\left( U_{0},r,\theta \right) x\right\Vert ^{2}=P_{r}\ast \frac{d\mu
_{U_{0},x}}{d\theta }  \label{23}
\end{equation}%
where $\mu _{U_{0},x}$ is the measure defined by $\mu _{U_{0},x}\left(
\Omega \right) =\left\langle E_{U_{0}}\left( \Omega \right) x,x\right\rangle 
$ for any Borel set $\Omega .$ Then by Fatou's theorem for Poisson
integrals, 
\begin{equation}
q_{U_{0},r}\left( x,\theta \right) \rightarrow \frac{d\mu _{U_{0},x}}{%
d\theta }\left( \theta \right) \text{ as }r\rightarrow 1^{-}\text{ a.e. }%
\theta  \label{22}
\end{equation}%
and $q_{U_{0},r}\left( x,\theta \right) $ is bounded a.e. $\theta $ with $%
\int_{\mathbb{T}}q_{U_{0},r}\left( x,\theta \right) d\theta \leq \left\Vert
x\right\Vert ^{2}.$ By a similar discussion, the same property is
established for $q_{U,r}\left( y,\theta \right) $ for $y\in P_{ac}^{\infty
}\left( U,\mathcal{M}\right) \left( \mathcal{H}\right) .$ Let $x_{j}$ be a
vector in $P_{m_{j}}\left( \mathcal{H}\right) $ and let $y_{j}$ be a vector
in $Q_{m_{j}}\left( \mathcal{H}\right) $ satisfying 
\begin{equation*}
\left\Vert x-x_{j}\right\Vert +\left\Vert y-y_{j}\right\Vert \leq 2^{-j}.
\end{equation*}%
Since $P_{m}$ and $Q_{m}$ commute with the corresponding resolvents, 
\begin{equation*}
\widetilde{\alpha }_{\pm }^{\mathcal{M}}(x_{j},y_{j},r,\theta )=\gamma _{\pm
}\left( r\right) \left\langle J_{m_{j}}R^{\pm }\left( U_{0},r,\theta \right)
x_{j},R^{\pm }\left( U,r,\theta \right) y_{j}\right\rangle .
\end{equation*}%
Then by (\ref{thm5(11)}) and Theorem \ref{thm3}, 
\begin{equation}
\lim_{\underset{r\in (0,1)}{r\rightarrow 1}}\widetilde{\alpha }_{\pm }^{%
\mathcal{M}}(x_{j},y_{j},r,\theta )  \label{5}
\end{equation}%
exists for almost every $\theta $ and every $j.$ Furthermore, 
\begin{equation*}
\sum_{j=1}^{\infty }\int_{0}^{2\pi }\frac{d\mu _{U_{0},x-x_{j}}}{d\theta }%
d\theta \leq \sum_{j=1}^{\infty }\left\Vert x-x_{j}\right\Vert ^{2}<\infty .
\end{equation*}%
Similarly, we have 
\begin{equation*}
\sum_{j=1}^{\infty }\int_{0}^{2\pi }\frac{d\mu _{U,y-y_{j}}}{d\theta }%
d\theta <\infty .
\end{equation*}

Then for almost all $\theta ,$ 
\begin{equation}
\frac{d\mu _{U_{0},x-x_{j}}}{d\theta }\rightarrow 0\text{ as }j\rightarrow
\infty .  \label{26}
\end{equation}%
and 
\begin{equation}
\frac{d\mu _{U,y-y_{j}}}{d\theta }\rightarrow 0\text{ as }j\rightarrow
\infty \text{ a.e. }\theta .  \label{36}
\end{equation}%
Since 
\begin{eqnarray*}
&&\left\vert \widetilde{\alpha }_{\pm }^{\mathcal{M}}\left( x,y,r,\theta
\right) -\widetilde{\alpha }_{\pm }^{\mathcal{M}}\left( x_{j},y_{j},r,\theta
\right) \right\vert \\
&\leq &q_{U_{0},r}\left( x-x_{j},\theta \right) ^{1/2}q_{U,r}\left( y,\theta
\right) ^{1/2}+q_{U_{0},r}\left( x_{j},\theta \right) ^{1/2}q_{U,r}\left(
y-y_{j},\theta \right) ^{1/2},
\end{eqnarray*}%
we have 
\begin{equation}
\lim_{j\rightarrow \infty }\underset{\underset{r\in (0,1)}{r\rightarrow 1}}{%
\lim \sup }\left\vert \widetilde{\alpha }_{\pm }^{\mathcal{M}}\left(
x,y,r,\theta \right) -\widetilde{\alpha }_{\pm }^{\mathcal{M}}\left(
x_{j},y_{j},r,\theta \right) \right\vert =0  \label{5.5}
\end{equation}%
for a.e. $\theta $ by (\ref{22}), (\ref{26}) and (\ref{36})$.$ For such a $%
\theta ,$%
\begin{eqnarray*}
&&\left\vert \widetilde{\alpha }_{\pm }^{\mathcal{M}}\left( x,y,r,\theta
\right) -\widetilde{\alpha }_{\pm }^{\mathcal{M}}\left( x,y,s,\theta \right)
\right\vert \\
&\leq &\left\vert \widetilde{\alpha }_{\pm }^{\mathcal{M}}\left(
x,y,r,\theta \right) -\widetilde{\alpha }_{\pm }^{\mathcal{M}}\left(
x_{j},y_{j},r,\theta \right) \right\vert +\left\vert \widetilde{\alpha }%
_{\pm }^{\mathcal{M}}\left( x_{j},y_{j},r,\theta \right) -\widetilde{\alpha }%
_{\pm }^{\mathcal{M}}\left( x_{j},y_{j},s,\theta \right) \right\vert \\
&&+\left\vert \widetilde{\alpha }_{\pm }^{\mathcal{M}}\left(
x_{j},y_{j},s,\theta \right) -\widetilde{\alpha }_{\pm }^{\mathcal{M}}\left(
x,y,s,\theta \right) \right\vert .
\end{eqnarray*}%
Then by (\ref{5}), we have%
\begin{eqnarray*}
&&\underset{\underset{r,s\in (0,1)}{r,s\rightarrow 1}}{\lim \sup }\left\vert 
\widetilde{\alpha }_{\pm }^{\mathcal{M}}\left( x,y,r,\theta \right) -%
\widetilde{\alpha }_{\pm }^{\mathcal{M}}\left( x,y,s,\theta \right)
\right\vert \\
&\leq &2\underset{\underset{r,s\in (0,1)}{r,s\rightarrow 1}}{\lim \sup }%
\left\vert \widetilde{\alpha }_{\pm }^{\mathcal{M}}\left( x,y,r,\theta
\right) -\widetilde{\alpha }_{\pm }^{\mathcal{M}}\left( x_{j},y_{j},r,\theta
\right) \right\vert \text{ for every }j.
\end{eqnarray*}%
As $j\rightarrow \infty ,$ we conclude that 
\begin{equation*}
\widetilde{\alpha }_{\pm }^{\mathcal{M}}\left( x,y,\theta \right) :=\lim_{%
\underset{r\in (0,1)}{r\rightarrow 1}}\widetilde{\alpha }_{\pm }^{\mathcal{M}%
}\left( x,y,r,\theta \right) \text{ exists a.e. }\theta
\end{equation*}%
by (\ref{5}) and (\ref{5.5}). We now verify the uniform integrability of $%
\left\{ q_{U_{0},r}\left( x,\cdot \right) :0<r<1\right\} $. For a measurable
set $E\subseteq \mathbb{T}$ and $f_{M}:=\min \left\{ \frac{d\mu _{U_{0},x}}{%
d\theta },M\right\} \in \mathcal{L}^{1}\left( \mathbb{T}\right) $ with $M>0,$
we have 
\begin{eqnarray*}
\int_{E}q_{U_{0},r}\left( x,\theta \right) d\theta &\leq &\int_{E}\left(
P_{r}\ast f_{M}\right) \left( \theta \right) d\theta +\int_{\mathbb{T}%
}\left( P_{r}\ast \left( \frac{d\mu _{U_{0},x}}{d\theta }-f_{M}\right)
\right) \left( \theta \right) d\theta \\
&\leq &M\lambda \left( E\right) +\left\Vert \frac{d\mu _{U_{0},x}}{d\theta }%
-f_{M}\right\Vert _{\mathcal{L}^{1}\left( \mathbb{T}\right) }.
\end{eqnarray*}%
This bound is independent of $r.$ Given $\varepsilon >0,$ we choose $M>0$
such that 
\begin{equation*}
\left\Vert \frac{d\mu _{U_{0},x}}{d\theta }-f_{M}\right\Vert _{\mathcal{L}%
^{1}\left( \mathbb{T}\right) }<\frac{\varepsilon }{2}.
\end{equation*}%
Then, whenever $\lambda \left( E\right) <\frac{\varepsilon }{2M},$%
\begin{equation*}
\sup_{0<r<1}\int_{E}q_{U_{0},r}\left( x,\theta \right) d\theta <\varepsilon .
\end{equation*}%
Together with $\int_{\mathbb{T}}q_{U_{0},r}\left( x,\theta \right) d\theta
\leq \left\Vert x\right\Vert ^{2},$ we conclude that $\left\{
q_{U_{0},r}\left( x,\cdot \right) :0<r<1\right\} $ is uniformly integrable.
Since 
\begin{eqnarray*}
\int_{E}\left\vert \widetilde{\alpha }_{\pm }^{\mathcal{M}}\left(
x,y,r,\theta \right) \right\vert d\theta &\leq &\left(
\int_{E}q_{U_{0},r}\left( x,\theta \right) d\theta \right) ^{1/2}\left(
\int_{E}q_{U,r}\left( y,\theta \right) d\theta \right) ^{1/2} \\
&\leq &\left\Vert y\right\Vert \left( \int_{E}q_{U_{0},r}\left( x,\theta
\right) d\theta \right) ^{1/2},
\end{eqnarray*}%
$\widetilde{\alpha }_{\pm }^{\mathcal{M}}\left( x,y,r,\cdot \right) $ is
uniformly integrable$.$ By Vitali's convergence theorem, 
\begin{equation}
\widetilde{\alpha }_{\pm }^{\mathcal{M}}\left( x,y,r,\cdot \right)
\longrightarrow \widetilde{\alpha }_{\pm }^{\mathcal{M}}\left( x,y,\cdot
\right) \text{ in }\mathcal{L}^{1}\left( \mathbb{T}\right) .  \label{24}
\end{equation}%
In particular, 
\begin{eqnarray*}
\int_{\mathbb{T}}\left\vert \widetilde{\alpha }_{\pm }^{\mathcal{M}}\left(
x,y,\theta \right) \right\vert d\theta &\leq &\left( \int_{\mathbb{T}%
}q_{U_{0},r}\left( x,\theta \right) d\theta \right) ^{1/2}\left( \int_{%
\mathbb{T}}q_{U,r}\left( y,\theta \right) d\theta \right) ^{1/2} \\
&\leq &\left\Vert x\right\Vert \left\Vert y\right\Vert .
\end{eqnarray*}%
Consequently, 
\begin{equation*}
\widetilde{w}_{\pm }\left( U,U_{0},J,\mathcal{M}\right)
\end{equation*}%
exists and 
\begin{equation*}
\left\langle \widetilde{w}_{\pm }\left( U,U_{0},J,\mathcal{M}\right)
h,k\right\rangle =\int_{\mathbb{T}}\lim_{\underset{r\in (0,1)}{r\rightarrow 1%
}}\widetilde{\alpha }_{\pm }^{\mathcal{M}}\left( h,k,r,\theta \right)
d\theta =\lim_{\underset{r\in (0,1)}{r\rightarrow 1}}\int_{\mathbb{T}}%
\widetilde{\alpha }_{\pm }^{\mathcal{M}}\left( h,k,r,\theta \right) d\theta .
\end{equation*}%
Hence 
\begin{eqnarray*}
\widetilde{w}_{\pm }\left( U,U_{0},J,\mathcal{M}\right) &=&wot\text{-}\lim_{%
\underset{r\in (0,1)}{r\rightarrow 1}}\widetilde{w}_{\pm }\left( U,U_{0},J,r,%
\mathcal{M}\right) \\
&=&wot\text{-}\lim_{\underset{r\in (0,1)}{r\rightarrow 1}}P_{ac}^{\infty
}\left( U,\mathcal{M}\right) w_{\pm }\left( U,U_{0},J,r,\mathcal{M}\right) .
\end{eqnarray*}
\end{proof}

The following lemma is an important ingredient of the last main result of
this section.

\begin{lemma}
\label{lem16}Let $\mathcal{M\subseteq B}\left( \mathcal{H}\right) $ be a
countably decomposable, properly infinite von Neumann algebra with a
faithful normal semifinite tracial weight $\tau .$ Let $U$, $U_{0}$ be two
unitaries in $\mathcal{M}$ and $J$ be a unitary in $\mathcal{M}.$ If 
\begin{equation*}
V:=JU_{0}-UJ\in \mathcal{L}^{1}\left( \mathcal{M},\tau \right) \mathcal{\cap
M,}
\end{equation*}%
then $\widetilde{w}_{\pm }=\widetilde{w}_{\pm }\left( U,U_{0},J,\mathcal{M}%
\right) $ exist and 
\begin{equation*}
\underset{r\rightarrow 1^{-}}{sot\text{-}\lim }w_{\pm }\left( U,U_{0},J,r,%
\mathcal{M}\right) =\widetilde{w}_{\pm }.
\end{equation*}%
Moreover, 
\begin{equation*}
\widetilde{w}_{\pm }^{\ast }\widetilde{w}_{\pm }=P_{ac}^{\infty }\left(
U_{0},\mathcal{M}\right) ,\text{ }U\widetilde{w}_{\pm }=\widetilde{w}_{\pm
}U_{0},\text{ }P_{ac}^{\infty }\left( U,\mathcal{M}\right) \widetilde{w}%
_{\pm }=\widetilde{w}_{\pm }
\end{equation*}
\end{lemma}

\begin{proof}
Throughout this proof, we use the same notation as in the proofs of Lemma %
\ref{lem6} and Theorem \ref{th3}. We consider first the positive sign. For $%
x\in \mathcal{D}_{0},$ we define 
\begin{equation*}
\rho _{x}\left( n\right) =\left( \sum_{k\geq n}\left\Vert
G_{0}U_{0}^{-k}x\right\Vert ^{2}\right) ^{1/2}.
\end{equation*}%
By (\ref{thm5(1)}), 
\begin{equation}
\rho _{x}\left( n\right) \rightarrow 0\text{ as }n\rightarrow \infty .
\label{19}
\end{equation}%
By (\ref{35}), 
\begin{equation*}
T_{k+1}^{\ast }-T_{k}^{\ast }=-U_{0}^{k+1}V^{\ast }U^{-k}.
\end{equation*}%
We claim that 
\begin{equation}
\left\Vert \left( T_{n}-T_{m}\right) x\right\Vert ^{2}\leq 2c_{x}\left( \rho
_{x}\left( n+1\right) +\rho _{x}\left( m+1\right) \right) \text{ for }%
n,m\geq 0.  \label{21}
\end{equation}%
To prove this, we set $F_{n,m}=T_{n}^{\ast }T_{m}.$ Note $T_{n}$ is unitary.
For fixed $n,m$ and $d=m-n,$ we obtain that 
\begin{eqnarray*}
&&F_{n+j+1},_{m+j+1}-F_{n+j},_{m+j}=T_{n+j+1}^{\ast }T_{m+j+1}-T_{n+j}^{\ast
}T_{m+j} \\
&=&\left( T_{n+j+1}^{\ast }-T_{n+j}^{\ast }\right) T_{m+j+1}+T_{n+j}^{\ast
}\left( T_{m+j+1}-T_{m+j}\right) \\
&=&-U_{0}^{n+j+1}V^{\ast }U^{d+1}JU_{0}^{-\left( m+j+1\right)
}-U_{0}^{n+j}J^{\ast }U^{d}VU_{0}^{-\left( m+j+1\right) }.
\end{eqnarray*}%
Note $G=G_{0}v^{\ast }.$ If we let $B_{1}=v^{\ast }U^{d+1}J$ and $%
B_{2}=v^{\ast }U^{-d}J,$ then $\left\Vert B_{1}\right\Vert \leq 1$ and $%
\left\Vert B_{2}\right\Vert \leq 1.$ After substituting $V=G^{\ast }G_{0},$
by (\ref{thm5(1)}) we obtain that 
\begin{eqnarray}
&&\left\vert \left\langle \left( F_{n+a,m+a}-F_{n,m}\right) x,x\right\rangle
\right\vert  \notag \\
&\leq &\left( \sum_{j=0}^{a-1}\left\Vert GU^{d+1}JU_{0}^{-\left(
m+j+1\right) }x\right\Vert ^{2}\right) ^{1/2}\left(
\sum_{j=0}^{a-1}\left\Vert G_{0}U_{0}^{-\left( n+j+1\right) }x\right\Vert
^{2}\right) ^{1/2}  \notag \\
&&+\left( \sum_{j=0}^{a-1}\left\Vert GU^{-d}JU_{0}^{-\left( n+j\right)
}x\right\Vert ^{2}\right) ^{1/2}\left( \sum_{j=0}^{a-1}\left\Vert
G_{0}U_{0}^{-\left( m+j+1\right) }x\right\Vert ^{2}\right) ^{1/2}  \notag \\
&\leq &\left( \sum_{k\in \mathbb{Z}}\left\Vert
G_{0}B_{1}U_{0}^{-k}x\right\Vert ^{2}\right) ^{1/2}\rho _{x}\left(
n+1\right) +\left( \sum_{k\in \mathbb{Z}}\left\Vert
G_{0}B_{2}U_{0}^{-k}x\right\Vert ^{2}\right) ^{1/2}\rho _{x}(m+1)  \notag \\
&\leq &c_{x}\left( \rho _{x}\left( n+1\right) +\rho _{x}\left( m+1\right)
\right) .  \label{20}
\end{eqnarray}

By (\ref{35}) and (\ref{19}), for fixed $n,m,$%
\begin{equation*}
\left\Vert \left( T_{n+a}-T_{m+a}\right) x\right\Vert \leq \left\Vert
G\right\Vert \sum_{k=\min \left( n,m\right) +a}^{\max \left( n,m\right)
+a-1}\left\Vert G_{0}U_{0}^{-k-1}x\right\Vert \rightarrow 0\text{ as }%
a\rightarrow \infty
\end{equation*}%
Since every $T_{n}$ is unitary, 
\begin{equation}
\left( F_{n+a,m+a}-I\right) x=T_{n+a}^{\ast }\left( T_{m+a}-T_{n+a}\right)
x\rightarrow 0\text{ as }a\rightarrow \infty .  \label{20.1}
\end{equation}%
By (\ref{20}) and (\ref{20.1}), letting $a\rightarrow \infty ,$ we obtain%
\begin{eqnarray}
\left\vert \left\langle \left( I-F_{n,m}\right) x,x\right\rangle \right\vert
&=&\left\vert \left\langle \left( I-F_{n+a,m+a}+F_{n+a,m+a}-F_{n,m}\right)
x,x\right\rangle \right\vert  \notag \\
&\leq &\left\vert \left\langle I-F_{n+a,m+a})x,x\right\rangle \right\vert
+\left\vert \left\langle \left( F_{n+a,m+a}-F_{n,m}\right) x,x\right\rangle
\right\vert  \notag \\
&\leq &c_{x}\left( \rho _{x}\left( n+1\right) +\rho _{x}\left( m+1\right)
\right) .  \label{33.3}
\end{eqnarray}%
It follows that 
\begin{eqnarray*}
\left\Vert \left( T_{n}-T_{m}\right) x\right\Vert ^{2} &=&2\func{Re}%
\left\langle \left( I-F_{n,m}\right) x,x\right\rangle \\
&\leq &2c_{x}\left( \rho _{x}\left( n+1\right) +\rho _{x}\left( m+1\right)
\right) \text{ for }n,m\geq 0
\end{eqnarray*}%
by (\ref{33.3}). Then the inequality (\ref{21}) holds.

We denote $w_{\pm }\left( U,U_{0},r,J,\mathcal{M}\right) $ by $w_{\pm
}\left( r\right) .$ Then 
\begin{eqnarray*}
\left\Vert w_{+}\left( r\right) x-w_{+}\left( s\right) x\right\Vert ^{2}
&\leq &\sum \left( 1-r^{2}\right) \left( 1-s^{2}\right)
r^{2n}s^{2m}\left\Vert \left( T_{n}-T_{m}\right) x\right\Vert ^{2} \\
&\leq &2c_{x}\left( k_{x}\left( r\right) +k_{x}\left( s\right) \right)
\rightarrow 0\text{ as }r,s\rightarrow 1^{-}
\end{eqnarray*}%
by (\ref{16}) and (\ref{19}) where $k_{x}\left( r\right) =\sum_{n\geq
0}\left( 1-r^{2}\right) r^{2n}\rho _{x}\left( n+1\right) $. We also have 
\begin{eqnarray*}
&&\left\vert \left\langle P_{ac}^{\infty }\left( U_{0},\mathcal{M}\right)
x,x\right\rangle -\left\langle w_{+}\left( r\right) ^{\ast }w_{+}\left(
r\right) x,x\right\rangle \right\vert \\
&=&\left\Vert x\right\Vert ^{2}-\left\Vert w_{+}\left( r\right) x\right\Vert
^{2} \\
&=&\frac{1}{2}\sum_{n,m\geq 0}\left( 1-r^{2}\right) ^{2}r^{2\left(
n+m\right) }\left\Vert \left( T_{n}-T_{m}\right) x\right\Vert ^{2} \\
&\leq &2c_{x}k_{x}\left( r\right) \rightarrow 0\text{ as }r\rightarrow 1^{-}.
\end{eqnarray*}%
Since $\mathcal{D}_{0}$ is dense in $P_{ac}^{\infty }\left( U_{0},\mathcal{M}%
\right) \mathcal{H}$ and $\left\Vert w_{+}\left( r\right) \right\Vert \leq
1, $ there is an operator $Z_{+}\in \mathcal{M}$ such that 
\begin{equation}
w_{+}\left( r\right) \overset{sot}{\rightarrow }Z_{+},\text{ }%
Z_{+}P_{ac}^{\infty }\left( U_{0},\mathcal{M}\right) =Z_{+}\text{ and }%
Z_{+}^{\ast }Z_{+}=P_{ac}^{\infty }\left( U_{0},\mathcal{M}\right) .
\label{30}
\end{equation}%
The relation $UT_{n}=T_{n+1}U_{0}$ gives%
\begin{eqnarray*}
Uw_{+}\left( r\right) &=&\left( 1-r^{2}\right) \sum_{n=0}^{\infty
}r^{2n}UT_{n}P_{ac}^{\infty }\left( U_{0},\mathcal{M}\right) =\left(
1-r^{2}\right) \sum_{n=0}^{\infty }r^{2n}T_{n+1}U_{0}P_{ac}^{\infty }\left(
U_{0},\mathcal{M}\right) \\
&=&\frac{1-r^{2}}{r^{2}}\sum_{k=1}^{\infty }r^{2k}T_{k}U_{0}P_{ac}^{\infty
}\left( U_{0},\mathcal{M}\right) =\frac{1}{r^{2}}\left( w_{+}\left( r\right)
-\left( 1-r^{2}\right) JP_{ac}^{\infty }\left( U_{0},\mathcal{M}\right)
\right) U_{0}.
\end{eqnarray*}%
Hence 
\begin{eqnarray*}
&&Uw_{+}\left( r\right) -w_{+}\left( r\right) U_{0} \\
&=&\frac{1}{r^{2}}w_{+}\left( r\right) U_{0}-\frac{\left( 1-r^{2}\right) }{%
r^{2}}JP_{ac}^{\infty }\left( U_{0},\mathcal{M}\right) U_{0}-w_{+}\left(
r\right) U_{0} \\
&=&\frac{1-r^{2}}{r^{2}}\left( w_{+}\left( r\right) -JP_{ac}^{\infty }\left(
U_{0},\mathcal{M}\right) \right) U_{0}.
\end{eqnarray*}%
It implies that 
\begin{equation}
UZ_{+}=Z_{+}U_{0}.  \label{18}
\end{equation}

Let $P\in \mathcal{P}_{ac}^{\infty }\left( U_{0},\mathcal{M}\right) $ and $%
D=PZ_{+}^{\ast }$. Then by (\ref{18}) 
\begin{equation*}
\sum_{k\in \mathbb{Z}}\left\Vert DU^{k}h\right\Vert ^{2}=\sum_{k\in \mathbb{Z%
}}\left\Vert PU_{0}^{k}Z_{+}^{\ast }h\right\Vert ^{2}\leq 2\pi
C_{P}\left\Vert h\right\Vert ^{2}
\end{equation*}%
where $C_{P}$ is the smoothness constant of $P.$ It implies that $D$ is $U$%
-smooth. Therefore by Corollary \ref{cor2}, 
\begin{equation*}
\overline{Ran\left( Z_{+}P\right) }=\overline{Ran\left( D^{\ast }\right) }%
\subseteq P_{ac}^{\infty }\left( U,\mathcal{M}\right) \left( \mathcal{H}%
\right) .
\end{equation*}%
Hence 
\begin{equation}
P_{ac}^{\infty }\left( U,\mathcal{M}\right) Z_{+}=Z_{+}.  \label{29}
\end{equation}%
It follows that 
\begin{equation*}
P_{ac}^{\infty }\left( U,\mathcal{M}\right) w_{+}\left( r\right) \overset{sot%
}{\rightarrow }Z_{+}.
\end{equation*}%
Hence by Lemma \ref{lem6}, 
\begin{equation}
Z_{+}=\widetilde{w}_{+}\left( U,U_{0},J,\mathcal{M}\right) .  \label{30.1}
\end{equation}%
Thus 
\begin{equation}
P_{ac}^{\infty }\left( U,\mathcal{M}\right) \widetilde{w}_{+}=\widetilde{w}%
_{+},\text{ }\widetilde{w}_{+}^{\ast }\widetilde{w}_{+}=P_{ac}^{\infty
}\left( U_{0},\mathcal{M}\right) ,\text{ }U\widetilde{w}_{+}=\widetilde{w}%
_{+}U_{0}  \label{46}
\end{equation}%
by (\ref{30}) and (\ref{29}).

For the negative sign, we consider the pair $U^{\ast },U_{0}^{\ast }.$ Then $%
JU_{0}^{\ast }-U^{\ast }J\in \mathcal{L}^{1}\left( \mathcal{M},\tau \right) 
\mathcal{\cap M}$ and 
\begin{equation*}
P_{ac}^{\infty }\left( U_{0}^{\ast },\mathcal{M}\right) =P_{ac}^{\infty
}\left( U_{0},\mathcal{M}\right) ,\text{ }P_{ac}^{\infty }\left( U^{\ast },%
\mathcal{M}\right) =P_{ac}^{\infty }\left( U,\mathcal{M}\right) .
\end{equation*}%
Applying the positive-sign discussion to the triple $\left( U^{\ast
},U_{0}^{\ast },J\right) ,$ we obtain that 
\begin{equation}
P_{ac}^{\infty }\left( U,\mathcal{M}\right) \widetilde{w}_{-}=\widetilde{w}%
_{-},\text{ }\widetilde{w}_{-}^{\ast }\widetilde{w}_{-}=P_{ac}^{\infty
}\left( U_{0},\mathcal{M}\right) ,\text{ }U\widetilde{w}_{-}=\widetilde{w}%
_{-}U_{0}.  \label{46.6}
\end{equation}%
Note that we have the following equality 
\begin{equation*}
w_{-}\left( U,U_{0},J,r,\mathcal{M}\right) =\frac{1}{r^{2}}\left(
w_{+}\left( U^{\ast },U_{0}^{\ast },J,r,\mathcal{M}\right) -\left(
1-r^{2}\right) JP_{ac}^{\infty }\left( U_{0},\mathcal{M}\right) \right) .
\end{equation*}%
Hence 
\begin{equation}
w_{-}\left( U,U_{0},J,r,\mathcal{M}\right) \overset{sot}{\rightarrow }%
\widetilde{w}_{-}.  \label{30.2}
\end{equation}%
Combining (\ref{30}), (\ref{30.1}), (\ref{46}), (\ref{30.2}) and (\ref{46.6}%
), the proof is completed.
\end{proof}

Throughout the following proof of the lemma, we use the same notation as in
the proofs of Lemma \ref{lem6} and Theorem \ref{th3}.

\begin{lemma}
\label{lem26}Let $\mathcal{M\subseteq B}\left( \mathcal{H}\right) $ be a
countably decomposable, properly infinite von Neumann algebra with a
faithful normal semifinite tracial weight $\tau .$ Let $U$, $U_{0}$ be two
unitaries in $\mathcal{M}$ and $J$ be a unitary in $\mathcal{M}.$ If 
\begin{equation*}
V:=JU_{0}-UJ\in \mathcal{L}^{1}\left( \mathcal{M},\tau \right) \mathcal{\cap
M,}
\end{equation*}%
then $w_{\pm }\left( U,U_{0},J,\mathcal{M}\right) $ exist and 
\begin{equation*}
w_{\pm }\left( U,U_{0},J,\mathcal{M}\right) =\widetilde{w}_{\pm }\left(
U,U_{0},J,\mathcal{M}\right) .
\end{equation*}
\end{lemma}

\begin{proof}
We write $\widetilde{w}_{\pm }=\widetilde{w}_{\pm }\left( U,U_{0},J,\mathcal{%
M}\right) .$ For $x\in P_{ac}^{\infty }\left( U_{0},\mathcal{M}\right)
\left( \mathcal{H}\right) ,$ we define 
\begin{equation}
e_{r}\left( \theta \right) =\gamma _{\pm }\left( r\right) \left\Vert \left(
J-\widetilde{w}_{\pm }\right) R^{\pm }\left( U_{0},r,\theta \right)
x\right\Vert ^{2}\geq 0.  \label{28}
\end{equation}%
By Lemma \ref{lem16},%
\begin{equation*}
\widetilde{w}_{\pm }R^{\pm }\left( U_{0},r,\theta \right) =R^{\pm }\left(
U,r,\theta \right) \widetilde{w}_{\pm }.
\end{equation*}%
Since $J$ is a unitary and $\widetilde{w}_{\pm }^{\ast }\widetilde{w}_{\pm
}=P_{ac}^{\infty }\left( U_{0},\mathcal{M}\right) $ by (\ref{23}), we have%
\begin{eqnarray}
e_{r}\left( \theta \right) &=&\gamma _{\pm }\left( r\right) \left\Vert
\left( J-\widetilde{w}_{\pm }\right) R^{\pm }\left( U_{0},r,\theta \right)
x\right\Vert ^{2}  \label{27} \\
&=&\gamma _{\pm }\left( r\right) \left\Vert JR^{\pm }\left( U_{0},r,\theta
\right) x-\widetilde{w}_{\pm }R^{\pm }\left( U_{0},r,\theta \right)
x\right\Vert  \notag \\
&=&\gamma _{\pm }\left( r\right) \left\Vert JR^{\pm }\left( U_{0},r,\theta
\right) x\right\Vert ^{2}+\gamma _{\pm }\left( r\right) \left\Vert 
\widetilde{w}_{\pm }R^{\pm }\left( U_{0},r,\theta \right) x\right\Vert ^{2} 
\notag \\
&&-2\func{Re}\gamma _{\pm }\left( r\right) \left\langle JR^{\pm }\left(
U_{0},r,\theta \right) x,\widetilde{w}_{\pm }R^{\pm }\left( U_{0},r,\theta
\right) x\right\rangle  \notag \\
&=&2q_{U_{0},r}\left( x,\theta \right) -2\func{Re}\widetilde{\alpha }_{\pm
}^{\mathcal{M}}\left( x,\widetilde{w}_{\pm }x,r,\theta \right) .  \notag
\end{eqnarray}%
Moreover, $\widetilde{w}_{\pm }x\in \mathcal{H}_{ac}\left( U\right) $ by
Lemma \ref{lem16}$.$ And we have%
\begin{equation*}
\lim_{\underset{r\in (0,1)}{r\rightarrow 1}}\widetilde{\alpha }_{\pm }^{%
\mathcal{M}}\left( x,\widetilde{w}_{\pm }x,r,\theta \right) \text{ exists
a.e. }\theta
\end{equation*}%
with 
\begin{equation*}
\int_{\mathbb{T}}\lim_{\underset{r\in (0,1)}{r\rightarrow 1}}\widetilde{%
\alpha }_{\pm }^{\mathcal{M}}\left( x,\widetilde{w}_{\pm }x,r,\theta \right)
=\left\langle \widetilde{w}_{\pm }x,\widetilde{w}_{\pm }x\right\rangle .
\end{equation*}%
Thus by (\ref{28}) and (\ref{27}), 
\begin{equation*}
e\left( \theta \right) =\lim_{\underset{r\in (0,1)}{r\rightarrow 1}%
}e_{r}\left( \theta \right) =2\frac{d\mu _{U_{0},x}}{d\theta }-2\func{Re}%
\lim_{\underset{r\in (0,1)}{r\rightarrow 1}}\widetilde{\alpha }_{\pm }^{%
\mathcal{M}}\left( x,\widetilde{w}_{\pm }x,r,\theta \right) \geq 0
\end{equation*}%
and 
\begin{equation*}
\int_{\mathbb{T}}e\left( \theta \right) d\theta =2\left\Vert x\right\Vert
^{2}-2\left\Vert \widetilde{w}_{\pm }x\right\Vert ^{2}=0.
\end{equation*}%
Hence $e\left( \theta \right) =0$ a.e. $\theta .$ Therefore 
\begin{equation}
e_{r}\left( \theta \right) =\gamma _{\pm }\left( r\right) \left\Vert \left(
J-\widetilde{w}_{\pm }\right) R^{\pm }\left( U_{0},r,\theta \right)
x\right\Vert ^{2}\rightarrow 0\text{ as }r\rightarrow 1^{-}\text{ for a.e. }%
\theta .  \label{25.1}
\end{equation}%
Let $y\in \mathcal{H}$ and set $z=\left( 1-P_{ac}^{\infty }\left( U,\mathcal{%
M}\right) \right) y.$ Then $R^{\pm }\left( U,r,\theta \right) z\in \left(
1-P_{ac}^{\infty }\left( U,\mathcal{M}\right) \right) \left( \mathcal{H}%
\right) $ and 
\begin{equation*}
\widetilde{w}_{\pm }R^{\pm }\left( U_{0},r,\theta \right) x\in
P_{ac}^{\infty }\left( U,\mathcal{M}\right) \left( \mathcal{H}\right) .
\end{equation*}%
Then $\left\langle \widetilde{w}_{\pm }R^{\pm }\left( U_{0},r,\theta \right)
x,R^{\pm }\left( U,r,\theta \right) z\right\rangle =0.$ Hence 
\begin{eqnarray*}
\left\vert \alpha _{\pm }^{\mathcal{M}}\left( x,z,r,\theta \right)
\right\vert &=&\gamma _{\pm }\left( r\right) \left\vert \left\langle \left(
J-\widetilde{w}_{\pm }\right) R^{\pm }\left( U_{0},r,\theta \right) x,R^{\pm
}\left( U,r,\theta \right) z\right\rangle \right\vert \\
&\leq &e_{r}\left( \theta \right) ^{1/2}q_{U,r}\left( z,\theta \right)
^{1/2}.
\end{eqnarray*}%
Now we write $\mu _{U,z}=f_{z}\left( \theta \right) d\theta +\mu _{U,z}^{s}$
with $f_{z}\left( \theta \right) \in \mathcal{L}_{+}^{1}\left( \mathbb{T}%
\right) $ where $\mu _{U,z}^{s}$ is singular with respect to $\lambda ^{d}.$
Hence $q_{U,r}\left( z,\theta \right) \rightarrow f_{z}\left( \theta \right) 
$ a.e. $\theta $. By (\ref{25.1}) 
\begin{equation*}
\left\vert \alpha _{\pm }^{\mathcal{M}}\left( x,z,r,\theta \right)
\right\vert \leq e_{r}\left( \theta \right) ^{1/2}q_{U,r}\left( z,\theta
\right) ^{1/2}\rightarrow 0\text{ for a.e. }\theta .
\end{equation*}%
Consequently, $\lim_{\underset{r\in (0,1)}{r\rightarrow 1}}\alpha _{\pm }^{%
\mathcal{M}}\left( x,y,r,\theta \right) $ exists and 
\begin{eqnarray*}
\lim_{\underset{r\in (0,1)}{r\rightarrow 1}}\alpha _{\pm }^{\mathcal{M}%
}\left( x,y,r,\theta \right) &=&\lim_{\underset{r\in (0,1)}{r\rightarrow 1}%
}(\alpha _{\pm }^{\mathcal{M}}\left( x,P_{ac}^{\infty }\left( U,\mathcal{M}%
\right) y,r,\theta \right) +\alpha _{\pm }^{\mathcal{M}}\left( x,\left(
1-P_{ac}^{\infty }\left( U,\mathcal{M}\right) \right) y,r,\theta \right) ) \\
&=&\lim_{\underset{r\in (0,1)}{r\rightarrow 1}}\alpha _{\pm }^{\mathcal{M}%
}\left( x,P_{ac}^{\infty }\left( U,\mathcal{M}\right) y,r,\theta \right)
=\lim_{\underset{r\in (0,1)}{r\rightarrow 1}}\widetilde{\alpha }_{\pm }^{%
\mathcal{M}}\left( x,y,r,\theta \right) .
\end{eqnarray*}%
Hence by the definitions of $w_{\pm }\left( U,U_{0},J,\mathcal{M}\right) $, $%
\widetilde{w}_{\pm }\left( U,U_{0},J,\mathcal{M}\right) $ and Lemma \ref%
{lem16}$,$ we have $w_{\pm }\left( U,U_{0},J,\mathcal{M}\right) $ exists and 
\begin{equation*}
w_{\pm }\left( U,U_{0},J,\mathcal{M}\right) =\widetilde{w}_{\pm }\left(
U,U_{0},J,\mathcal{M}\right) .
\end{equation*}%
The proof is completed.
\end{proof}

The following theorem is the last main result of this paper. It is an
analogue of the Kato-Rosenblum theorem for a pair of unitaries in a
semifinite von Neumann algebra $\mathcal{M}$.

\begin{theorem}
\label{thm16}Let $\mathcal{M\subseteq B}\left( \mathcal{H}\right) $ be a
countably decomposable, properly infinite von Neumann algebra with a
faithful normal semifinite tracial weight $\tau .$ Let $U$, $U_{0}$ be two
unitaries in $\mathcal{M}$ and let $J$ be a unitary in $\mathcal{M}.$ If 
\begin{equation*}
V:=JU_{0}-UJ\in \mathcal{L}^{1}\left( \mathcal{M},\tau \right) \mathcal{\cap
M,}
\end{equation*}%
then $W_{\pm }=W_{\pm }\left( U,U_{0},J,\mathcal{M}\right) $ and $w_{\pm
}\left( U,U_{0},J,\mathcal{M}\right) $ exist and 
\begin{equation*}
W_{\pm }=w_{\pm }\left( U,U_{0},J,\mathcal{M}\right) .
\end{equation*}%
Moreover, 
\begin{equation*}
W_{\pm }W_{\pm }^{\ast }=P_{ac}^{\infty }\left( U,\mathcal{M}\right) ,\text{ 
}W_{\pm }^{\ast }W_{\pm }=P_{ac}^{\infty }\left( U_{0},\mathcal{M}\right)
\end{equation*}%
and $UW_{\pm }=W_{\pm }U_{0}.$
\end{theorem}

\begin{proof}
By Lemma \ref{lem26}, $w_{\pm }=w_{\pm }\left( U,U_{0},J,\mathcal{M}\right) $
exist. By Remark \ref{rem11}, for proving the theorem, we need to establish
the existence of the weak wave operators $\widetilde{W}_{\pm }\left(
U,U_{0},J,\mathcal{M}\right) $. Set $T_{n}=U^{n}JU_{0}^{-n}$ for $n\in 
\mathbb{Z}.$ Using the notation in the proof of Lemma \ref{lem6}, for $x\in
\cup P_{m}\mathcal{H}$ and $y\in \cup Q_{m}\mathcal{H},$ we choose $m$ large
enough that $P_{m}x=x$ and $Q_{m}y=y.$ Then by the proof of Lemma \ref{lem6}%
, 
\begin{eqnarray*}
&&\left\langle T_{n}x,y\right\rangle =\left\langle P_{ac}^{\infty }\left( U,%
\mathcal{M}\right) T_{n}P_{ac}^{\infty }\left( U_{0},\mathcal{M}\right)
x,y\right\rangle \\
&=&\left\langle Q_{m}T_{n}P_{m}x,y\right\rangle =\left\langle
U^{n}J_{m}U_{0}^{-n}x,y\right\rangle \\
&\rightarrow &\left\langle W_{\pm ,m}x,y\right\rangle \text{ as }%
n\rightarrow \pm \infty .
\end{eqnarray*}%
Since $\left\Vert P_{ac}^{\infty }\left( U,\mathcal{M}\right)
T_{n}P_{ac}^{\infty }\left( U_{0},\mathcal{M}\right) \right\Vert \leq 1,$
the density of $\cup P_{m}\mathcal{H}$ and $\cup Q_{m}\mathcal{H}$ gives the
existence of 
\begin{equation*}
\widetilde{W}_{\pm }:=\widetilde{W}_{\pm }\left( U,U_{0},J,\mathcal{M}%
\right) =wot\text{-}\lim_{n\rightarrow \pm \infty }P_{ac}^{\infty }\left( U,%
\mathcal{M}\right) T_{n}P_{ac}^{\infty }\left( U_{0},\mathcal{M}\right) .
\end{equation*}%
Then for any given $\varepsilon >0$ and $x,y\in \mathcal{H},$ there is an
integer $N,$ such that for $n\geq N,$ we have 
\begin{equation*}
\left\vert \left\langle \left( P_{ac}^{\infty }\left( U,\mathcal{M}\right)
T_{n}P_{ac}^{\infty }\left( U_{0},\mathcal{M}\right) -\widetilde{W}_{\pm
}\right) x,y\right\rangle \right\vert <\varepsilon .
\end{equation*}%
By (\ref{16}), for such $N,$ we obtain 
\begin{eqnarray*}
&&\left\vert \left\langle \left( P_{ac}^{\infty }\left( U,\mathcal{M}\right)
w_{+}\left( U,U_{0},J,r,\mathcal{M}\right) -\widetilde{W}_{+}\right)
x,y\right\rangle \right\vert \\
&\leq &\left( 1-r^{2}\right) \sum_{n\geq 0}r^{2n}\left\langle \left(
P_{ac}^{\infty }\left( U,\mathcal{M}\right) T_{n}P_{ac}^{\infty }\left(
U_{0},\mathcal{M}\right) -\widetilde{W}_{+}\right) x,y\right\rangle \\
&\leq &\left( 1-r^{2}\right) \sum_{n=0}^{N-1}\left\vert \left\langle \left(
P_{ac}^{\infty }\left( U,\mathcal{M}\right) T_{n}P_{ac}^{\infty }\left(
U_{0},\mathcal{M}\right) -\widetilde{W}_{+}\right) x,y\right\rangle
\right\vert +\varepsilon r^{2N}.
\end{eqnarray*}%
By a similar argument for $\widetilde{W}_{-},$ we conclude that 
\begin{equation*}
wot\text{-}\lim_{\underset{r\in (0,1)}{r\rightarrow 1}}P_{ac}^{\infty
}\left( U,\mathcal{M}\right) w_{\pm }\left( U,U_{0},r,J,\mathcal{M}\right) =%
\widetilde{W}_{\pm }\left( U,U_{0},J,\mathcal{M}\right) .
\end{equation*}%
Hence by Lemma \ref{lem6} 
\begin{equation*}
\widetilde{W}_{\pm }\left( U,U_{0},J,\mathcal{M}\right) =\widetilde{w}_{\pm
}\left( U,U_{0},J,\mathcal{M}\right) .
\end{equation*}%
Applying Lemma \ref{lem16},%
\begin{equation*}
\widetilde{W}_{\pm }\left( U,U_{0},J,\mathcal{M}\right) ^{\ast }\widetilde{W}%
_{\pm }\left( U,U_{0},J,\mathcal{M}\right) =P_{ac}^{\infty }\left( U_{0},%
\mathcal{M}\right) .
\end{equation*}%
Since $J$ is a unitary, $\widetilde{W}_{\pm }\left( U_{0},U_{0},J^{\ast }J,%
\mathcal{M}\right) =P_{ac}^{\infty }\left( U_{0},\mathcal{M}\right) .$ It
follows that 
\begin{equation*}
\widetilde{W}_{\pm }\left( U,U_{0},J,\mathcal{M}\right) ^{\ast }\widetilde{W}%
_{\pm }\left( U,U_{0},J,\mathcal{M}\right) =\widetilde{W}_{\pm }\left(
U_{0},U_{0},J^{\ast }J,\mathcal{M}\right) .
\end{equation*}%
Then by Remark \ref{rem11}, Lemma \ref{lem16} and Lemma \ref{lem26}, we
obtain that%
\begin{equation*}
W_{\pm }\left( U,U_{0},J,\mathcal{M}\right) =\widetilde{W}_{\pm }\left(
U,U_{0},J,\mathcal{M}\right) =w_{\pm }\left( U,U_{0},J,\mathcal{M}\right)
\end{equation*}%
and 
\begin{equation*}
W_{\pm }\left( U,U_{0},J,\mathcal{M}\right) ^{\ast }W_{\pm }\left( U,U_{0},J,%
\mathcal{M}\right) =P_{ac}^{\infty }\left( U_{0},\mathcal{M}\right) .
\end{equation*}%
Since $J^{\ast }U-U_{0}J^{\ast }=U_{0}V^{\ast }U\in \mathcal{L}^{1}\left( 
\mathcal{M},\tau \right) \mathcal{\cap M},$ by a similar discussion, we know 
$W_{\pm }\left( U_{0},U,J^{\ast },\mathcal{M}\right) $ exist and 
\begin{equation*}
W_{\pm }\left( U_{0},U,J^{\ast },\mathcal{M}\right) ^{\ast }W_{\pm }\left(
U_{0},U,J^{\ast },\mathcal{M}\right) =P_{ac}^{\infty }\left( U,\mathcal{M}%
\right) .
\end{equation*}%
We also have 
\begin{equation*}
W_{\pm }\left( U_{0},U,J^{\ast },\mathcal{M}\right) =\widetilde{W}_{\pm
}\left( U_{0},U,J^{\ast },\mathcal{M}\right) =\widetilde{W}_{\pm }\left(
U,U_{0},J,\mathcal{M}\right) ^{\ast }=W_{\pm }\left( U,U_{0},J,\mathcal{M}%
\right) ^{\ast },
\end{equation*}%
By Lemma \ref{lem16}, we also conclude that $UW_{\pm }=W_{\pm }U_{0}.$ This
completes the proof.
\end{proof}

Next we recall the definition of diagonal unitary operators in $\mathcal{M}.$

\begin{definition}
A unitary $U$ in $\mathcal{M}$ is called diagonal if there exists a family $%
\left\{ P_{n}\right\} $ of orthogonal projections in $\mathcal{M}$ and a
family $\left\{ \lambda _{n}\right\} _{n=1}^{\infty }$ of real numbers such
that $U=\sum_{n=1}^{\infty }e^{i\lambda _{n}}P_{n.}.$
\end{definition}

Note that for a diagonal unitary $U,$ we always have $P_{ac}\left( U\right)
=0$. Then by Corollary (\ref{cor2}), we obtain that $P_{ac}^{\infty }\left(
U,\mathcal{M}\right) =0$ for such $U.$ The following corollary is an
application of Theorem \ref{thm16}.

\begin{corollary}
A unitary $U$ in $\mathcal{M}$ with $P_{ac}^{\infty }\left( U,\mathcal{M}%
\right) \neq 0$ is not a perturbation by $\mathcal{L}^{1}\left( \mathcal{M}%
,\tau \right) \mathcal{\cap M}$ of a diagonal unitary $U$ in $\mathcal{M.}$
\end{corollary}

\begin{example}
\bigskip We will give an example such that a unitary $U$ in $\mathcal{M}$
with $P_{ac}\left( U\right) \neq 0$ is a perturbation by $\mathcal{L}%
^{1}\left( \mathcal{M},\tau \right) \mathcal{\cap M}$ of a diagonal unitary $%
D$ in $\mathcal{M.}$

Suppose $\mathcal{M=L}^{\infty }\left( \left[ 0,1\right] \right) \otimes 
\mathcal{B}\left( l^{2}\left( \mathbb{Z}\right) \right) $, $\mathcal{H}=%
\mathcal{L}^{2}\left( \left[ 0,1\right] \right) \otimes l^{2}\left( \mathbb{Z%
}\right) $ and $\tau \left( T\right) =\int_{0}^{1}Tr\left( T\left( t\right)
\right) dt$ for $0\leq T\in \mathcal{M}$ where $Tr$ is the normal trace on $%
\mathcal{B}\left( l^{2}\left( \mathbb{Z}\right) \right) .$ Let $%
Q=e_{0}\otimes e_{0}$ be the rank-one projection onto $\mathbb{C}e_{0}$ and
let $E=1\otimes Q.$ Then $\tau \left( E\right) =1.$ We partition $(0,1]$
into the intervals $I_{n}=(2^{-n},2^{1-n}],$ $n\geq 1$ and set 
\begin{equation*}
P_{n}=\chi _{I_{n}}\otimes Q\text{ and }\lambda _{n}=2\pi 2^{-n}.
\end{equation*}%
The projections $P_{n}$ are pairwise orthogonal, and $\sum_{n\geq 1}P_{n}=E.$
Define the unitaries $U$ and $D$ by 
\begin{equation*}
D=\left( I-E\right) +\sum_{n=1}^{\infty }e^{i\lambda _{n}}P_{n}\text{ }
\end{equation*}%
and 
\begin{equation*}
U=1\otimes \left( I-Q\right) +\lambda \left( t\right) \otimes Q\text{ where }%
\lambda \left( t\right) =e^{2\pi it}.
\end{equation*}%
It is clear that $P_{ac}\left( D\right) =0$ and $P_{ac}\left( U\right) =E.$
For $t\in I_{n},$%
\begin{equation*}
\left( U-D\right) \left( t\right) =\left( e^{2\pi it}-e^{i\lambda
_{n}}\right) Q,
\end{equation*}%
it follows that 
\begin{equation*}
\left\Vert U-D\right\Vert _{\mathcal{L}^{1}\left( \mathcal{M},\tau \right)
}=\sum_{n=1}^{\infty }\int_{I_{n}}\left\vert e^{2\pi it}-e^{i\lambda
_{n}}\right\vert dt\leq 2.
\end{equation*}%
Thus $U-D\in \mathcal{L}^{1}\left( \mathcal{M},\tau \right) \mathcal{\cap M}%
. $ It implies that 
\begin{equation*}
U=D+K
\end{equation*}%
where $K\in \mathcal{L}^{1}\left( \mathcal{M},\tau \right) \mathcal{\cap M.}$
\end{example}

\begin{equation*}
\end{equation*}

\textbf{AI statement:} During the preparation of this manuscript, the author
used GPT 6 Sol for searching necessary references, examining mathematical
arguments, and language editing. The author takes full responsibility for
the final version.

\bigskip 
\begin{equation*}
\end{equation*}

\vspace{1cm}

\end{document}